\documentclass[a4paper,11pt]{amsart}
\usepackage{amssymb}
\usepackage{amscd}

\def\Box{\square}

\def\frak{\mathfrak}
\def\Bbb{\mathbb}
\def\Cal{\mathcal}

\def\prettylongrightarrow{\relbar\joinrel\relbar\joinrel
                          \relbar\joinrel\relbar\joinrel\longrightarrow}

\def\kth{k^{\underline{{\rm th}}}}

\let\suchthat=\mid
\def\endrk{\hbox{$|\!\!|\!\!|\!\!|\!\!|\!\!|\!\!|$}}

\def\bbW{{\mathbb{W}}}
\def\bbX{{\mathbb{X}}}
\def\bbY{{\mathbb{Y}}}
\def\bbZ{{\mathbb{Z}}}

\newcommand{\hook}{\raisebox{-0.35ex}{\makebox[0.6em][r]
{\scriptsize $-$}}\hspace{-0.15em}\raisebox{0.25ex}{\makebox[0.4em][l]{\tiny
 $|$}}}

\let\hash=\sharp

\let\d=\delta
\let\e=\varepsilon
\let\f=\varphi
\let\i=\iota
\let\N=\nabla
\let\w=\omega
\let\r=\rho

\def\uw{\underline{\wedge}}

\def\AMB{{\Cal{M}}}

\def\fsl{{\frak{sl}}}
\def\fg{{\mathfrak{g}}} 
\def\fsu{{\mathfrak{su}}}
\def\fu{{\mathfrak{u}}}
\newcommand{\bh}{\mbox{\boldmath{$h$}}}
\newcommand{\qq}{\mbox{\boldmath{$q$}}}

\def\LOT{{\rm LOT}}

\def\ix{{\bf p}}

\newcommand{\miniX}{\mbox{\boldmath{$\scriptstyle{X}$}}}
\newcommand{\miniY}{\mbox{\boldmath{$\scriptstyle{Y}$}}}
\newcommand{\LX}{{\Cal L}_{\!\miniX}}
\newcommand{\KX}{{\Cal K}_{\!\miniX}}

\newcommand{\NX}{{\mbox{\boldmath$ \nabla$}}_{\!\miniX}}

\newcommand{\ip}[2]{\langle{#1},{#2}\rangle}
\let\Euler=\NX
\newcommand{\IT}[1]{{\rm(}{\it{\!#1}}{\rm)}}

\newcommand{\newc}{\newcommand}

\renewcommand{\Im}{\operatorname{Im}}

\newcommand{\x}{\times}

\let\ccdot\cdot
\def\cdot{\hbox to 2.5pt{\hss$\ccdot$\hss}}

\newcommand{\om}{\omega}
\renewcommand{\phi}{\varphi}
\newcommand{\ph}{\varphi}

\newcommand{\si}{\sigma}

\newcommand{\Om}{\Omega}

\newc{\aI}{\mbox{\boldmath{$ I$}}}
\newc{\aR}{\mbox{\boldmath{$ R$}}}
\newc{\aDeR}{\mbox{\boldmath{$ U$}}_B{}^P{}_C{}^Q}
\newc{\al}{\mbox{\boldmath$ \Delta$}}             
\newc{\nda}{\mbox{\boldmath$ \nabla$}}
\newc{\ad}{\mbox{\boldmath$ d$}}
\newc{\da}{\mbox{\boldmath$ \delta$}}
\newc{\aK}{\mbox{\boldmath{$ K$}}}
\newc{\aL}{\mbox{\boldmath{$ L$}}}

\newtheorem{theorem}{Theorem}[section]
\newtheorem{lemma}[theorem]{Lemma}
\newtheorem{proposition}[theorem]{Proposition}
\newtheorem{corollary}[theorem]{Corollary}

\newcommand{\cB}{{\Cal B}}
\newcommand{\cC}{{\Cal C}}
\newcommand{\cf}{{\Cal F}}
\newcommand{\cH}{{\Cal H}}
\newcommand{\cg}{{\Cal G}}
\newcommand{\cN}{{\Cal N}}
\newcommand{\cV}{{\Cal V}}

\newcommand{\ce}{{\Cal E}}
\newcommand{\cE}{{\Cal E}}
\newcommand{\cq}{{\Cal Q}}
\newcommand{\cQ}{{\Cal Q}}
\newcommand{\cK}{{\Cal K}}
\newcommand{\cL}{{\Cal L}}
\newcommand{\cT}{{\Cal T}}
\newcommand{\cR}{{\Cal R}}

\newcommand{\ct}{{\Cal T}}

\newcommand{\bW}{{\Bbb W}}

\newcommand{\bX}{{\Bbb X}}
\newcommand{\bZ}{{\Bbb Z}}
\newcommand{\bV}{{\Bbb V}}
\newcommand{\bL}{{\Bbb L}}
\newcommand{\bG}{{\Bbb G}}
\newcommand{\bQ}{{\Bbb Q}}
\newcommand{\bS}{{\Bbb S}}

\newcommand{\nd}{\nabla}

\newcommand{\Rho}{{\mbox{\sf P}}}
\newcommand{\Up}{\Upsilon}

\newcommand{\End}{\operatorname{End}}

\newcommand{\Ric}{\operatorname{Ric}}

\newcommand{\vol}{\mbox{\large\boldmath $ \epsilon$}}
\newcommand{\tstar}{\mbox{\Large $ \star$}}

\newcommand{\astar}{\mbox{\Large\boldmath $ \star$}}

\newcommand{\fb}  {\mbox{$                                     
\begin{picture}(9,8)(1.6,0.15)
\put(1,0.2){\mbox{$ \Box \hspace{-7.8pt} /$}}
\end{picture}$}} 

\newcommand{\afl}{\mbox{$
\begin{picture}(9,8)(1.6,0.15)
\put(1,0.2){\mbox{$ \al \hspace{-7.8pt} /$}}
\end{picture}$}}

\newcommand{\fD}{\mbox{$
\begin{picture}(9,8)(1.6,0.15)
\put(1,0.2){\mbox{$ D \hspace{-7.8pt} /$}}
\end{picture}$}}

\newcommand{\afD}                         
{\mbox{$
\begin{picture}(9,8)(1.6,0.15)
\put(1,0.2){\mbox{$ \D \hspace{-7.8pt} /$}}
\end{picture}$}}

\newcommand{\fl}{\mbox{$                                     
\begin{picture}(9,8)(1.6,0.15)
\put(1,0.2){\mbox{$ \Delta \hspace{-7.8pt} /$}}
\end{picture}$}}

\newcommand{\act}{\mbox{\boldmath{$\ct$}}}

\newcommand{\acf}{\mbox{\boldmath{$\Cal F$}}}
\newcommand{\acg}{\mbox{\boldmath{$\Cal G$}}}

\newcommand{\td}{\tilde{d}}
\newcommand{\dt}{\tilde{\d}}

\newcommand{\obL}{\overline{\Bbb L}}
\newcommand{\bK}{{\Bbb K}}

\newcommand{\bN}{{\Bbb N}}
\newcommand{\obN}{\overline{\bN}}

\def\eye{{\rm i}}
\let\i=\iota
\let\x=\xi

\newcommand{\nn}[1]{(\ref{#1})}

\newcommand{\tD}{\tilde{D}}

\newcommand{\D}{\mbox{\boldmath{$ D$}}}
\newcommand{\btD}{\tilde{\mbox{\boldmath{$ D$}}}}

\newcommand{\bF}{\mbox{{$\Bbb F$}}}
\newcommand{\bT}{\mbox{{$\Bbb T$}}}

\newcommand{\X}{\mbox{\boldmath{$ X$}}}

\newcommand{\Y}{\mbox{\boldmath{$ Y$}}}

\newcommand{\sX}{\mbox{\scriptsize\boldmath{$X$}}}        
                                                          
\newcommand{\h}{\mbox{\boldmath{$ h$}}}
\newcommand{\bg}{\mbox{\boldmath{$ g$}}}

\newcommand{\cce}{\tilde{\ce}}                         
\newcommand{\aM}{\tilde{M}}
\newcommand{\tF}{\tilde{F}}

\newcommand{\tV}{\tilde{V}}
\newcommand{\tU}{\tilde{U}}

\newcommand{\sbg}{\mbox{\scriptsize\boldmath{$g$}}}

\let\m=\mu

\let\G=\Gamma

\newcommand{\V}{{\mbox{\sf P}}}                   
\newcommand{\J}{{\mbox{\sf J}}}
\newcommand{\U}{{\mbox{\sf U}}}

\newc{\strutdd}{\rule{0mm}{5mm}}

\newcommand{\lpl}                         
{\mbox{$
\begin{picture}(12.7,8)(-.5,-1)
\put(2,0.2){$+$}
\put(6.2,2.8){\oval(8,8)[l]}
\end{picture}$}}

\usepackage{ifthen}

\newc{\tensor}[1]{#1}

\newc{\Mvariable}[1]{\mbox{#1}}

\newc{\down}[1]{{}_{
\ifthenelse{\equal{#1}{;}}{|}{#1}}}

\newc{\up}[1]{{}^{#1}}
\newc{\C}{C}

\newc{\JulyStrut}{\rule{0mm}{6mm}}
\newc{\midtenPan}{\mbox{\sf S}}
\newc{\midten}{\mbox{\sf T}}
\newc{\midtenEi}{\mbox{\sf U}}
\newc{\ATen}{\mbox{\sf E}}
\newc{\BTen}{\mbox{\sf F}}
\newc{\CTen}{\mbox{\sf G}}

\def\sideremark#1{\ifvmode\leavevmode\fi\vadjust{\vbox to0pt{\vss
 \hbox to 0pt{\hskip\hsize\hskip1em
 \vbox{\hsize3cm\tiny\raggedright\pretolerance10000
 \noindent #1\hfill}\hss}\vbox to8pt{\vfil}\vss}}}%

\begin{document}
\renewcommand{\today}{}
\title{Conformally invariant operators, differential forms, 
cohomology and a generalisation of Q-curvature}
\author{Thomas Branson and A. Rod Gover}

\address{Department of Mathematics\\
  The University of Iowa\\
  Iowa City IA 52242 USA}\email{branson@math.uiowa.edu}
\address{Department of Mathematics\\
  The University of Auckland\\
  Private Bag 92019\\
  Auckland 1\\
  New Zealand} \email{gover@math.auckland.ac.nz}

\vspace{10pt}

\renewcommand{\arraystretch}{1}
\maketitle
\renewcommand{\arraystretch}{1.5}

\pagestyle{myheadings}
\markboth{Branson \& Gover}{Conformal operators, forms, cohomology and $Q$}

\begin{abstract}
On conformal manifolds of even dimension $n\geq 4$ we construct a family
of new conformally invariant differential complexes, each containing
one coboundary operator of order greater than 1.  
Each bundle in each of these complexes appears either in the de Rham complex
or in its dual (which is a different complex in the non-orientable
case). 
Each of the new 
complexes is elliptic in case the conformal structure has
Riemannian signature.  We also construct gauge companion operators
which (for differential forms of order $k\leq n/2$) complete the
exterior derivative to a conformally invariant and (in the case of
Riemannian signature) elliptically coercive system.  These
(operator,gauge) pairs are used to define finite dimensional
conformally stable form subspaces which are are candidates for spaces
of conformal harmonics.  This generalises the $n/2$-form and $0$-form
cases, in which the harmonics are given by conformally invariant
systems.  These constructions are based on a family of operators on
closed forms which generalise in a natural way Branson's
Q-curvature. We give a universal construction of these new operators
and show that they yield new conformally invariant global pairings
between differential form bundles.  Finally we give a geometric
construction of a family of conformally invariant differential
operators between density-valued differential form bundles and develop
their properties (including their ellipticity type in the case of
definite conformal signature).  The construction is based on the
ambient metric of Fefferman and Graham, and its relationship to the
tractor bundles for the Cartan normal conformal connection.  For each
form order, our derivation yields an operator of every even order in
odd dimensions, and even order operators up to order $n$ in even
dimension $n$.  In the case of unweighted (or {\em true}) forms as
domain, these operators are the natural form analogues of the 
critical order conformal Laplacian of Graham et al., and are
key ingredients in the new differential
complexes mentioned above.
\end{abstract}
\noindent\thanks{TB gratefully acknowledges support from US NSF grant INT-9724781.}
\thanks{ARG gratefully acknowledges support from the Royal Society of 
New Zealand via Marsden Grant no.\ 02-UOA-108,
and to the New Zealand Institute of Mathematics and its Applications
for support via a Maclaurin Fellowship.}
\section{Introduction}

Conformal structure on manifolds is the natural setting for the study
of massless particles in Physics.
It also plays a role in curvature
prescription, in extremal problems for metrics in Riemannian geometry,
and in string and brane theories.  Via the Fefferman bundle and
metric, conformal structure also makes its presence felt in complex
and CR geometry.

The main point of this paper is the discovery and
construction of a host of new local and global conformally invariant
objects associated with what is perhaps the most fundamental
domain for conformal geometry, namely the {\em true} differential forms; that is, 
exterior powers $E^k$ of $T^*M$ endowed with its natural conformal weight.

One of the important concepts is that of a {\em conformally invariant
differential operator}, like the conformal Laplacian (sometimes called
the Yamabe operator), or the Maxwell operator.  
Such operators act on sections of vector bundles natural for conformal
structure, may be defined by universal natural formulae, and depend
only on the conformal structure (and not on any choice of metric
tensor from the conformal class).  Among the most basic bundles are
weighted differential forms.  In this paper, we give a geometric
construction of a family of conformally invariant differential
operators between weighted (i.e.\ density-valued) differential forms
on pseudo-Riemannian manifolds $M$ of arbitrary conformal curvature.
Our construction gives all operators of this type that are known by
abstract methods to exist in this general setting (and generalises the
construction of Graham et al.\ \cite{GJMS} which gives all
conformally invariant operators between scalar densities, the so-called 
{\em GJMS operators}).  (Conjecturally, no other invariant
form-density operators exist in the arbitrarily conformally curved setting;
some evidence in this direction is given in 
\cite{Grnon} and \cite{GoH}.)
While we feel that this result, giving a complete picture of form-density
operators, is an important aspect of the current work, 
we focus most of our attention on the 
deeper, hitherto unexpected structure associated to the true form operators,
which we believe will have long-term resonance.

What is more significant is that our construction gives a preferred
family of such operators and this is especially evident in the
subfamily of operators $L_k$ which act on true forms.  From previously
known results and general reasoning one knows that these operators
exist in even dimensions $n$, carry $k$-forms to $k$-forms of a
certain weight (or equivalently, in the oriented case, to true
$(n-k)$-forms), and have principal part $(\delta d)^{n/2-k}$, where
$\delta$ is the formal adjoint of the exterior derivative $d$.  The
$L_k$ we construct are formally self-adjoint and have factorisations
of the form $\d Md$ (or in more detail \nn{rarefac} below), where $d$
is the exterior derivative and $\d$ is its formal adjoint with respect
to the conformal structure. Thus these operators generalise the
Maxwell operator $\delta d$ on $(n/2-1)$-forms and give a family of
complexes, that we introduce for the first time here, the {\em detour
complexes} \nn{detourdetail}. In the case of Riemannian signature
these complexes are elliptic.  As we point out below the factorisation
of the $L_k$ is subtle and unexpected. This property is absolutely
crucial to the other constructions we present.

In the Riemannian signature case, the operators $L_k$ are evidently
non-elliptic, having only positive semidefinite leading symbol.  For
the Maxwell operator, this state of affairs is tied to the {\em gauge
fixing} problem; roughly speaking, the search for a suitable operator
$G$ for which the system $\delta d\f=0$, $G\f=0$ is elliptically
coercive.  The choices of $G$ usually employed for the Maxwell
operator, however (notably $G=\delta$, the {\em Coulomb gauge}), are
not friendly to the formulation of the problem in terms of conformal
structure alone, since the coupled system is not conformally
invariant.  In this paper, we give a geometric construction of a {\em
gauge companion} for each $L_k$.  This is an operator $G_k$ of order
$n-2k+1$ on $k$-forms for which the system $(L_k,G_k)$ is conformally
invariant and elliptically coercive.  In fact, $G_k$ lands in a bundle
of weighted $(k-1)$-forms, and has principal part
$\delta(d\delta)^{n/2-k}$.  $G_k$ is not itself conformally invariant
on arbitrary $k$-forms, but {\em is} invariant on the forms
annihilated by $L_k$. In particular the $G_k$ are conformally invariant
on the
subspace (recall the factorisation of the $L_k$) of closed forms. This
leads to a type of conformal de Rham Hodge theory that we describe
below.  We show that $L_k$ and
$G_k$ are aspects of a single conformally invariant  operator, valued in a bundle that is
reducible but indecomposable for conformal structure.

The pairs $(L_k,G_k)$ harbor further, still deeper structure,
generalising the {\em Q-curvature}, an object that has inspired much
recent activity; see 
\cite{CGY,CQY,CYAnnals,CY,FGrQ,FeffHir,RobinKengo,GrZ,gursk,gurskyv} 
and references later in this
paragraph.  Beckner's generalization to $S^n$
of the Moser-Trudinger inequality \cite{beckner,carlenl} has a natural
statement in terms of Q-curvature; see \cite{tbsrni2004}.
The Q-curvature was first
defined in \cite{tbkorea,tomsharp} in arbitrary even dimensions,
generalising the 4-dimensional construction of \cite{tbprsnl,tbbo91}.
$Q$ is a local scalar invariant that appears naturally in formulas for
quotients of functional determinants for pairs of conformal metrics
\cite{tbbo91,tbkorea,tomsharp,CYAnnals}.  It also has a natural
relation to the {\em Fefferman-Graham ambient construction}
\cite{FGast}, which imbeds a conformal manifold of dimension $n$ into
a pseudo-Riemannian manifold of dimension $n+2$ (to a certain finite
order in even dimensions), by formally solving the Goursat problem for
the Einstein equation.  In \cite{RobinKengo}, Graham and Hirachi show
that the total metric variation of $\int Q$ is the {\em
Fefferman-Graham tensor}; i.e.\ the obstruction, at the appropriate
order, to the power series solution for the ambient metric in even
dimensions.  This in turn makes the Q-curvature of interest in the
study of the AdS/CFT correspondence \cite{FGrQ} and in scattering
theory \cite{GrZ}.  The Q-curvature has analogues in other parabolic
geometries, for example CR geometry; see \cite{FeffHir}, and in 
dimension 3, the earlier work \cite{KH}.

One of the salient features of the Q-curvature is its
conformal deformation law.  Given metrics $g$ and $\hat g=
e^{2\omega}g$ with $\omega$ a smooth function,
\begin{equation}\label{confchgQ}
\hat Q=Q+L_0\w,
\end{equation}
where the convention is that hatted (resp.\ unhatted) quantities
are computed in $\hat g$ (resp.\ $g$).  Thus $Q$ is not a conformal
invariant, but rather an invariant with a {\em linear} conformal 
change law.    
Generically, the conformal change law
for an invariant density of the same weight as $Q$
contains differential expressions
of homogeneities $1,2,\cdots,n$ in $\omega$.
(If $Q$ is viewed as a function rather than a density, the 
conformal change law reads $\hat Qe^{n\w}=Q+L_0\w$.  
As a nonlinear curvature prescription law, this equation has analytic
behavior similar to its 2-dimensional special case, the
Gauss curvature prescription equation.)
Note that one of the operators from our $L_k$
series appears in \nn{confchgQ}; in fact, this is the 
{\em critical GJMS operator} constructed in \cite{GJMS}.
It is evident from \nn{confchgQ} that the critical GJMS operator,
itself a delicate and celebrated object, may be reconstructed
from a knowledge of $Q$; the original construction of $Q$, on
the other hand, made essential use of the whole series of GJMS operators.

Since $L_0$ generalises to $L_k$, it is plausible that $Q$ generalises to
some form-valued object $Q_k$.  In this paper, we give a natural
geometric construction of this generalisation.
These form analogues of $Q$ are not
form-densities, but rather differential operators that act between certain
invariant subquotients of section spaces of form-density bundles.
(The source space is a true form subquotient, and in the orientable
case, the target space may also be realised as such.)
The analogue of \nn{confchgQ} is
\begin{equation}\label{confchgQk}
\hat Q_k\,u=Q_ku+L_k(\w u)\ \ \mbox{for }u\mbox{ a {\em closed} }k\mbox{-form}.
\end{equation}
Equation \nn{confchgQk} hints at a possible role for $Q_k$ as a
cohomology map.  This role, in fact, materialises in our work below,
which, as mentioned, could be described as a conformal Hodge theory.
The relevant cohomologies and harmonic spaces are related to the
elliptic detour complexes, mentioned above; and to the operators $L_k$
and their gauge companions $G_k$.  This is explained in somewhat more
detail just below, in our itemised list of results.

Our construction of $G_k$ generalises and is inspired by the special
case of the Maxwell operator in dimension 4, for which Eastwood and
Singer \cite{EastSin1,EastSin2} constructed the corresponding gauge
companion.  We are quick to note, however, that our construction of
the $G_k$, and even of the $L_k$, requires more powerful techniques,
since constructing these as classical tensor formulas is not an
option: the size of such formulas would grow rapidly 
with the order $n-2k$, and the
number of invariant expressions that could possibly appear undergoes a
combinatorial explosion.

One of the devices that allows us to work in such generality here is
the Fefferman-Graham ambient metric construction mentioned above.  
All of our main
results, save for some of the operators of order $n$, are obtained
from a single uniform construction based on the ambient metric, and
its relation, as exposed in \cite{CapGoFG,GoPet}, to a class of vector
bundles natural for conformal structure, the so-called {\em tractor}
bundles.  Tractor bundles and their normal connections may be viewed
as structures associated to the Cartan normal connection
\cite{CapGotrans} but may also be constructed directly \cite{BEGo} by
an idea which, in the conformal setting, dates back to Thomas
\cite{T}.  Penrose's local twistor bundle is an example of a tractor
bundle, as are the {\em spannor} and {\em plyor} bundles of Irving
Segal and his collaborators \cite{plyors,pilot} (though these
references work only in the conformally flat case).  Here, to avoid
unnecessary background, we use the ambient manifold to actually define
the tractor bundles required.

Even with these tools, extracting information from the ambient
construction is not necessarily a straightforward process.  
One needs a conceptual and 
detailed understanding of how  tractor and form operators  on 
the underlying conformal manifold $M$ arise from ambient operators.
A key part of the work we do here is to extend the results in  
\cite{CapGoFG,GoPet} and 
set up a calculus which is capable of restricting pseudo-Riemannian
information in ambient space to conformal information on $M$.
This calculus of this paper involves commutation and anticommutation
relations for certain relevant operators in ambient space, which in
turn put us in contact with a naturally occurring copy of the 
8-dimensional Lie
superalgebra ${\frak{sl}}(2|1)$.  (Work of Holland and Sparling
\cite{HS} on powers of the ambient Dirac operator has turned up a 
5-dimensional superalgebra isomorphic to the orthosymplectic
algebra ${\frak{osp}(2|1)}$, which may be realized as a subalgebra
of ${\frak{sl}}(2|1)$.)
Mediating between ambient space and
the form bundles on $M$ are {\em form-tractor} bundles; viewed from
$\AMB$, these are essentially restrictions of form bundles; viewed
from $M$, they are semidirect sums of form bundles.  They may also be
productively viewed as bundles over the $(n+1)$-dimensional conformal
metric bundle $\cQ$.  

The following is a precise, compact guide to the principal objects that
we construct on the conformal manifold $M$, and 
(along with \nn{confchgQk} above) asserts their main
properties:
\begin{itemize}
\item There are natural (built polynomially 
from $\nabla,R$), formally self-adjoint
differential operators $L_k:\cE^k\to\cE^k[2k-n]$, which at each
choice of metric have the factorisation
\begin{equation}\label{rarefac}
L_k=\underbrace{\delta\Bigl\{\overbrace{(d\delta)^{n/2-k-1}+\LOT}^{Q_{k+1}}
\Bigr\}}_{G_{k+1}}\,d,
\end{equation}
up to a nonzero constant factor that depends on $n$ and $k$,
and which are conformally invariant: $\hat L_k=L_k$.
Here $\cE^k$ denotes the smooth $k$-forms, and $\cE^k[w]$
the smooth $k$-forms of conformal weight $w$.  (Our normalization
of the conformal weight is uniquely determined by the fact that 
$TM=\cE^1[2]$.)  Here $\LOT$ stands for 
``lower order terms'', and the hat has the same meaning as above:
$\hat g=e^{2\omega}g$ with $\omega$ a smooth function.
$Q_{k+1}$ is a universal (but not conformally invariant)
expression in the covariant derivative and curvature.
\item The $Q_{k+1}$ are formally self-adjoint.
\item The $L_k$ are not elliptic, but the system 
$(L_k,G_k)$
is graded injectively elliptic, and conformally invariant in the
sense that up to a nonzero constant multiple,
$$
\hat G_k-G_k=d\omega\wedge L_k.
$$
(The sense of the last expression, of course, is $\varphi\mapsto d\omega\wedge
L_k\varphi$.)
\item Let $\cC^k$ denote the closed $k$-forms.
Each operator in the diagram 
\begin{equation}\label{ciQ}
\cC^k\stackrel{Q_k}{\longrightarrow}\cE^k[2k-n]/\cR(L_k)\stackrel
{{\rm quotient}}{\prettylongrightarrow}\cE^k[2k-n]/\cR(\d)
\end{equation}
is conformally invariant.  Thus $Q_k$ gives rise to operators
from closed $k$-forms to either of the quotients in the diagram.
\item $Q_k:\cN(L_k)\to\cE^k[2k-n]/\cN(\delta)$ is conformally invariant.
\item $Q_01$ is the ``classical'' Q-curvature.
\item The diagram
\begin{equation}\label{detourdetail}
\cdots\stackrel{d}{\longrightarrow}\cE^{k-1}\stackrel{d}{\longrightarrow}
\cE^k
\stackrel{L_k}{\longrightarrow}\cE^k[2k-n]\stackrel{\delta}{\longrightarrow}
\cE^{k-1}[2(k-1)-n]\stackrel{\delta}{\longrightarrow}\cdots
\end{equation}
is an elliptic complex, the {\em detour complex}.
As a result, in the orientable case,
$$
\cdots\stackrel{d}{\longrightarrow}\cE^{k-1}\stackrel{d}{\longrightarrow}
\cE^k
\stackrel{\star L_k}{\longrightarrow}\cE^{n-k}\stackrel{d}{\longrightarrow}
\cE^{n-k+1}\stackrel{d}{\longrightarrow}\cdots,
$$
where $\star$ is the Hodge star operator,
is an elliptic complex.
\item In the complex \nn{detourdetail}, let $H^{n-k}_L$ be the
cohomology at $\cE^k[2k-n]$. Note from \nn{rarefac} we have the conformally
invariant surjection $H^{n-k}_L \to H^{n-k}$.  Then the conformally
invariant operators \nn{ciQ} compress to conformally invariant
operators acting betweeen finite-dimensional conformally invariant
vector spaces according to
$$
Q_k:\cH^k\to H_L^{n-k}\to H^{n-k},
$$
where $\cH^k$ is the null space of $G_k$ within the closed forms $\cC^k$.
The space $\cH^k$ is conformally invariant
because $G_k$ is 
invariant on $\cN(L_k)\supset\cN(d)$; it is finite dimensional
because $(L_k,G_k)$
elliptically coercive.
\end{itemize}

Note that the first two points above indicate that the $(L_k,G_k)$ for
various $k$ are interlocked in an interesting way: the first statement
relates $L_k$ to $G_{k+1}$, while the second relates $L_k$ to $G_k$.
An equation that evokes this interlocking quite readily follows from
\nn{confchgQk} and \nn{rarefac}: up to a nonzero constant factor that
depends on $n$ and $k$, we have
$$
\hat Q_ku=Q_ku+\d Q_{k+1}d(\w u)\ \mbox{ for }u\in\cC^k.
$$
A special case of this interlocking was obtained by Eastwood and
Singer in \cite{EastSin1}, where it was shown that (in the current
language, and up to nonzero constant factors) 
$L_0=G_1d$ and $\hat G_1-G_1=d\omega\wedge L_1$ in dimension
4.  The last point indicates that $\cH^k$ is a candidate for a space
of {\em conformal harmonics}; under mild restrictions, ${\rm
dim}\,\cH^k$ recovers the $\kth$ Betti number.  More generally, we
give estimates bounding the size of each $\cH^k$ in terms of the de
Rham cohomology and the cohomology of the detour complexes mentioned
above.

In fact, our results point to what may be a better generalisation of
the Maxwell equations than is provided by the $(n/2-1)$-forms.  In arbitrary
dimension, we may take U$(1)$-connections, represented by one-forms $A$,
and take the corresponding curvature $F=dA$.  A natural
conformally invariant system of equations on $F$ in even dimensions is then 
\begin{equation}\label{newmax}
dF=0,\ \ G_2F=0,
\end{equation}
and this specialises to the usual Maxwell equations $dF=0$, $\delta F=0$
in dimension 4.  Just as $F=dA$, 
$\delta dA=0$ implies the Maxwell equations on $F$
in dimension 4, the system $F=dA$, $0=L_1A=G_2dA$ implies the Maxwell-like 
system \nn{newmax}.  The interlocking of different orders is illustrated by
the fact that $G_1$ is the natural gauge companion operator for $L_1$, so 
that both $G_1$ and $G_2$ are involved in the problem.  

Despite early glimpses of the factorisation \nn{rarefac} at low order
in \cite{tbms}, Theorem 2.10, such factorisations of invariant
operators are rare and surprising, and not at all to be expected (see
\cite{Gosrni04})
from the curved translation principle of Eastwood and Rice
\cite{EastRice,Esrni}.  Via this factorisation and the formal
self-adjointness of $Q_1$, we immediately have a constructive proof of
the existence of a version of the critical GJMS operator $L_0$ that is
formally self-adjoint and annihilates constants; this was an issue in
the earlier construction of the Q-curvature ($Q_01$ in the present
language).  Earlier proofs of the existence of an $L_0$ of this form
were given in \cite{GrZ,FGrQ}.

The current work arose in part from a desire to extend the
Eastwood-Singer gauge companion idea mentioned above, and partly from
a desire to have a differential form generalization of the GJMS
construction with an optimally clean idea of extension to and
restriction from ambient space.  One of the motivations for the latter
desideratum is a need to clarify the conformal geometric meaning of
Q-curvature.  In the process of carrying this out, we have observed
unexpected and, in our opinion, exciting further structure,
culminating in the $Q_k$ (in their various incarnations as
differential operators on closed forms, and as cohomology maps).  To
explain briefly our use of the ambient construction, powers of the
ambient form Laplacian are shown to descend to conformally invariant
operators on form tractor bundles on $M$; these bundles are exterior
powers of the standard tractor bundle. These descended operators are
then composed fore and aft with certain tractor operators to yield
invariant operators on weighted forms over the underlying conformal
manifold $M$.  These tractor operators, which we anticipate will be of
independent interest, are also defined via the ambient metric, and an
extensive ambient form calculus is developed to establish the relevant
properties of all the tractor operators involved in our compositions.
Since the tractor bundle and connection are well understood in terms
of the underlying (pseudo-)Riemannian structures (for metrics from the
conformal class) there is a straightforward algorithm for expression
of our operators in terms of the Levi-Civita connection and its
curvature; this provides considerable scope for future worthwhile
work.

In the next section we present the main results. Proofs are included
in Section \ref{mainthmsappls} only if they are accessible given
results already presented, Otherwise the proofs are delayed until
Section \ref{mainproofs}. The arguments of that section take advantage
of our ambient calculus, presented in Section \ref{ambient}, and its
interpretation in terms of tractor bundles given in Section
\ref{tractor}. In fact, Section \ref{mainproofs} gives theorems
generalising many of the results of Section \ref{mainthmsappls}, as
well as other theorems of independent value. However these require the
technical background of the earlier sections in order to be stated.
Section \ref{Qsect} is devoted to defining the operators $Q_k$, and
deriving their main properties (Theorem \ref{Qthm}). In Section
\ref{var}, we show that these operators generalise and relate the
definitions of the Q-curvature given recently in \cite{GoPet} and
\cite{FeffHir}, while in Section \ref{nont} we give a nontriviality
result for the maps $\cH^k\to H_k(M)$ that they determine.  In Section
\ref{var} we also describe ways of proliferating other operators with
transformation laws similar to that of $Q_k$.

It is a pleasure to thank Michael Eastwood, Charlie Frohman and Ruibin
Zhang for helpful discussions, and Gregg Zuckerman for posing some
questions that motivated parts of this investigation.  We would also like
to thank the American Institute of Mathematics for sponsoring a workshop 
on Conformal Structure in Geometry Analysis, and Physics, during which
this paper was completed.

\section{The main theorems} \label{mainthmsappls}

Let $M$ be a smooth manifold of dimension $n\geq 3$. A 
{\em conformal structure\/} on $M$ of signature $(p,q)$ (with $p+q=n$) is an
equivalence class $[g]$ of smooth pseudo--Riemannian metrics of signature
$(p,q)$ on $M$, with two
metrics being equivalent if and only if one is obtained from the other
by multiplication with a positive smooth function. Equivalently a
conformal structure is a smooth ray subbundle $\cq\subset S^2T^*M$,
whose fibre over $x$ consists of the values of $g_x$ for all metrics
$g$ in the conformal class. A metric in the conformal class is a section 
of $\cq$. A conformal structure of signature $ (n,0)$ is termed {\em
Riemannian}.

We can view $ \cq$ as the total space of a principal bundle
$\pi:\cq\to M$ with structure group $\Bbb R_+$ and so there are
natural line bundles on $ (M,[g])$ induced from the
irreducible representations of ${\Bbb R}_+$. For $w\in {\Bbb R}$, we
write $E[w]$ for the line bundle induced from the representation of
weight $-w/2$ on ${\Bbb R}$ (that is ${\Bbb R}_+ \ni x\mapsto
x^{-w/2}\in {\rm End}(\Bbb R)$).  Thus a section of $E[w]$
corresponds to a real-valued function $ f$ on $\Cal Q$ with the
homogeneity property $f(x, \Omega^2 g)=\Omega^w f(x,g)$, where
$\Omega$ is a positive function on $M$, $x\in M$, and $g$ is a metric
from the conformal class $[g]$. We shall write $ \ce[w]$ for the
space of smooth sections of this bundle. 

Note that there is a tautological function $\bg$ on $\cq$ taking
values in $S^2T^*M$, namely the function which assigns to the point
$(x,g_x)\in \cq$ the metric $g_x$ at $x$.  
This is homogeneous of
degree 2 since $\bg (x,s^2 g_x) =s^2 g_x$. If $\si$ is any
non-vanishing function on $\cq$ homogeneous of degree $+1$ then
$\si^{-2} \bg$ is independent of the action of ${\Bbb R}_+$ on the
fibres of $\cq$, and so $\si^{-2} \bg$ descends to give a metric from
the conformal class. Thus $\bg$ determines and is equivalent to a
canonical section of $S^2T^*M\otimes E[2]$ (called the {\em conformal
metric}) that we also denote $\bg$. 
This in turn determines a
canonical section $\bg^{-1}$ of $S^2T^*M\otimes E[-2]$ with the
property that a single contraction between these gives the 
identity endomorphism
of $TM$.  The conformal metric gives an isomorphism of $TM$ with
$T^*M[2]$ that we will view as an identification.
We usually also denote by $\si\in \ce[1]$ the density on $M$
equivalent to the homogeneous function $\si$ on $\cq$. Since
$\si^{-2}\bg$ is a metric from the conformal class we term $\si$ a
choice of {\em conformal scale}.

We will use the notation $\ce^k$ for the space of smooth sections of
$\wedge^kT^*M$, for which we shall sometimes use the alternative
notation $E^k$, and we write $\ce^k[w]$ for the smooth sections of the
tensor product $\wedge^kT^*M \otimes E[w]$.  Some statements about
forms or form-densities admit simpler formulations if we allow values
of $k$ falling outside the range $0\le k\le n$; by convention these
$\wedge^k$ are zero bundles.  
We also write $
E_k[w]$ (with section space $\ce_k[w]$) as a shorthand for $
E^k[w+2k-n]$.  This notation is suggested by the duality between the
section spaces $\ce^k$ and $\ce_k $ as follows. For $\phi\in \ce^k$
and $\psi\in \ce_k $, with one of these compactly supported, there is
the natural conformally invariant global pairing
\begin{equation}\label{pair}
\phi,\psi \mapsto \langle \phi,\psi \rangle :=\int_M \phi\cdot \psi\,
d\mu_{\sbg},
\end{equation}
where $\ph\cdot \psi\in\ce[-n]$
denotes a complete contraction between $\phi $ and $\psi$. This scales
so that when $M$ is orientable we have
$$
\langle \phi,\psi \rangle =\int_M \phi \wedge \star \psi 
$$
where $\star $ is the conformal Hodge star operator.
(On orientable conformal manifolds the bundle of
volume densities can be canonically identified with $ \ce[-n]$ and so
the Hodge star operator for each metric from the conformal class
induces an isomorphism that we shall also term the Hodge star
operator: $\star: \ce^k \cong \ce^{n-k}[n-2k]$.) 
The integral is also well-defined if instead $\phi\in \ce_k$ and
$\psi\in \ce^k $.  If we also denote this pairing by $\langle\cdot,\cdot
\rangle$, then $\langle\phi,\psi\rangle=\langle\psi,\phi\rangle$.

We write $\delta$ for the formal adjoint of the exterior derivative
with respect to 
our pairings.
That is, for $\phi\in\cE^k$ and $\mu \in
\ce_{k+1}$, at least 
one having compact support, we have $ \langle d\phi,\mu \rangle  
=\langle \phi,\delta \mu \rangle $. 
The notation
$\cC^k$ is used for the space of closed $ k$-forms and $ \cC_k$ denotes 
the formal dual space 
$\ce_k/{\Cal R}(\delta)$.

Our first main theorem concerns the construction and form of a family
of natural conformally invariant operators between density valued
differential form bundles in dimension $n\geq 3$.  We say that $P$ is
a {\em natural
differential operator} if $P_g$ can be written as a universal
polynomial in covariant derivatives with coefficients depending
polynomially on the conformal metric, its
inverse, the curvature tensor and its covariant derivatives.  The
coefficients of natural operators are called {\em natural tensors}. In
the case that they are scalar they are often also called {\em
Riemannian invariants}.  We say $P$ is a {\em conformally invariant
differential
operator} if it is well defined on conformal structures (i.e.\ is
independent of a choice of conformal scale).

\begin{theorem}\label{long} 
For each choice of $ k\in \{0,1,\cdots , n\}$, let 
$\ell\in\{0,1,2,\ldots\}$ if $n$ is odd and, if $n$ is even 
let $\ell\in \{0,1,\ldots, n/2\}$.
Let
$$
w=k+\ell-n/2.
$$
On conformal $n$-manifolds,
there is a formally self-adjoint conformally invariant natural
differential operator
$$
L_k^\ell:\ce^k[w]\to \ce_k[-w]
$$
of order $2\ell$. 
If $n$ is even, $k\leq n/2$,
and $w=0$, we have
\begin{equation}\label{mk}
L_k:=
L_k^{n/2-k}=\delta M_k d,
\end{equation}
where $M_k:\cC^{k+1}\to \cC_{k+1}$
is a conformally invariant operator.

On Riemannian conformal manifolds the following holds.
The operator $L_k^\ell$ is elliptic
if and only if $k\neq n/2\pm \ell$,
and it is positively elliptic if and only if $k\notin[n/2-\ell,n/2+\ell]$.
For each $k$ the differential operator sequence
\begin{equation}\label{detour}
\ce^0\stackrel{d}{\to}\cdots\stackrel{d}{\to}\ce^{k-1}\stackrel{d}{\to}
\ce^k\stackrel{L_k}{\to}\ce_k\stackrel{\delta}{\to}\ce_{k-1}
\stackrel{\delta}{\to}
\cdots\stackrel{\delta}{\to}\ce_0
\end{equation}
is an elliptic complex.
\end{theorem}

The operator $L^1_{n/2-1}$ is, up to a non-zero multiple, the Maxwell
operator in even dimension $n$. $L^0_{n/2}$ is zero.  For $k\neq 0$
the existence of conformally invariant differential operators between
the spaces $\ce^k[w]$ and $\ce_k[-w]$ as in the theorem follows from
the general theory in \cite{EastSlo}. The main point here is the
special form \nn{mk}, together with an explicit construction of these
operators.  This construction employs tractor calculus and the
Fefferman-Graham ambient construction in tandem to generalise the
scalar density results in \cite{GJMS} to density valued differential
forms.  For almost all cases this construction is described in
expression \nn{bLdef} based on operators given in \nn{bGdef},
\nn{average} and Section \ref{ambext}.  The main results concerning
order and ellipticity are the subject of Proposition \ref{Lelliptic}.
(Since we are generally only concerned with operators up to a
non-vanishing constant multiple, here and throughout the article we
say an operator is {\em positively elliptic} if it has positive {\em
  or} negative leading symbol.)  Naturality and the form \nn{mk} are
established in Theorem \ref{Main} parts \IT{i} and \IT{ii}. There are
two classes of exceptional operators not given by expression
\nn{bLdef}: (a) the operators satisfying $k=\ell+n/2$; (b) The
operators of order $n$ when $k\geq 1$.  (These classes share the
operator of order $n$ on $n$-forms of weight $n$.)  The operators of
type (a) are exactly those which have the spaces $\ce_k$, with $k\geq
n/2$, as domain, and are treated in Theorem \ref{colong} below. On $k$
forms with $k\geq n/2$, they are given by $\star L_{n-k} \star $,
where $L_{n-k}$ is as in \nn{mk} above (and so are never elliptic). In
particular the operator of order $n$ on weighted $n$-forms arises this
way.  Proposition \ref{ordern} gives a construction of the other
operators of order $n$. The anomalous behaviour of the operators in
the two exceptional classes is not unexpected. In the case of
operators of type (a) it arises because $\d$ acts invariantly on the
domain bundles $\ce^k[2k-n]$. This means a certain differential
splitting operator (see Proposition \ref{splitop}) involved in the
general construction \nn{bLdef} fails for these bundles. The failure
of the operators of type (b) to arise from \nn{bLdef} is a reflection
of the conformal invariance of the Fefferman-Graham obstruction
tensor. While we do not elaborate on this point in the current
article, in dimension 4 this claim is clear from section 3 of
\cite{Grnon}.

The complex \nn{detour} will be referred to as the {\em $k^{\rm th}$
de Rham detour complex}.  We write $H^k_L(M)$ for the cohomology of
this complex at the point $\stackrel{d}{\to} \ce^k\stackrel{L_k}{\to}$
and $H^L_k(M)$ for the cohomology at
$\stackrel{L_k}{\to}\ce_k\stackrel{\delta}{\to}$.  In view of the
factorisation $L_k=\d M_k d$, it is immediate that there is a
canonical conformally invariant invariant injection $H^k(M)\to
H^k_L(M)$ and similarly a canonical surjection $H_k^L(M)\to H_k(M)$.
Here, in accordance with our other conventions, by $H_k(M)$ we mean
$\cN(\d)/\cR(\d)$ at $\ce_k$. (So on oriented manifolds $H_k(M)\cong
H^{n-k}(M)$.)  

On compact Riemannian conformal manifolds, 
Hodge theory shows that 
\nn{pair} determines a perfect pairing between the
standard de Rham cohomology $H^k(M)$ and $H_k(M)$. That is, via
\nn{pair} $H_k(M)$ is just the vector space dual of $H^k(M)$. Next
observe that since $L_k$ is formally self-adjoint it follows easily
from standard Hodge theory that $\dim (H^k_L(M))$ is finite and $\dim
(H^k_L(M))=\dim (H_k^L(M))$. On the other hand it is easily verified
that the pairing \nn{pair} descends to a well defined pairing between
$H^k_L(M)$ and $H_k^L(M)$.
We state this in a proposition.

\begin{proposition}\label{Hpairing}
On compact Riemannian conformal manifolds {\rm(}\ref{pair}{\rm)} induces 
an invariant perfect pairing between  $H^k_L(M)$ and $ H_k^L(M)$. 
\end{proposition}
\noindent{\bf Proof:} The invariance assertion is clear 
by construction.  If we fix an arbitrary choice of conformal scale
then, for each $k$, the spaces $\ce_k$ can be identified with the spaces
$\ce^k$ and the pairing \nn{pair} gives an inner product on each of
the spaces $\ce^k$.  In this setting, using that $L_k$ is formally
self-adjoint, the standard Hodge theory of the complex \nn{detour}
gives
\begin{equation}\label{HajAtScale}
\cR (d)\oplus \cR (L_k) \oplus (\cN(\d )\cap \cN (L_k))
\end{equation}
as the Hodge decomposition of the space $\ce^k$ and also 
of $\ce_k$. 
For any class $[\phi']\in
H^k_L(M)$ we can find a representative $\phi\in \cN(\d )\cap \cN (L_k)$
and, via the identification of $\ce^k$ with $\ce_k$, this is also the
preferred representative of a class in $H_k^L(M)$. 
The pairing of these classes produces $(\phi,\phi)$, which is 
positive 
if $[\phi']$ is non-zero; this establishes nondegeneracy.
\quad $\Box$ \\

\noindent{\bf Remark:} The pairing in the proposition induces and is
 equivalent 
 a symplectic inner
 product on $H_L^k(M)\oplus H_k^L(M)$, via
 $\langle(v,v'),(w,w')\rangle=\langle v,w'\rangle -\langle w,v'\rangle $. 
$\qquad\endrk$

To state the next theorem we need one more result.

\begin{proposition}\label{Gksemidef}
On a conformal manifold of dimension $ n$, for each $ k\in 0,1,\cdots ,
n+1$ there is a natural indecomposable bundle $\bG_k$ with a natural
subbundle isomorphic to $E_{k-1}$, and corresponding quotient
isomorphic to $E_k$.
\end{proposition} 
The bundle $\bG_k$, its dual and their weighted
variants are defined in expression \nn{cgdefs} of Section
\ref{ambext}.  In particular 
$\bG_k$ is a subbundle of a certain tractor bundle (see
Section \ref{tractor}) and arises naturally from the ambient
construction. From either picture the properties described in the
proposition are immediate.
To summarise the composition series of $ \bG_k$ we will
often use the semi-direct sum notation
$$
\bG_k=E_k \lpl E_{k-1},
$$
or, on the level of section spaces,
$$
\cg_k=\ce_k \lpl \ce_{k-1}.
$$
For the natural quotient bundle map onto $E_k $ we shall write $
q^k:\bG_k\to E_k$.

{}From Theorem \ref{long}, on Riemannian conformal manifolds the
operators $L_k^\ell$ are elliptic except when operating on unweighted
forms, that is when the dimension is even and $\ell=n/2-k$. The next
theorem asserts that these operators $L_k:\ce^k\to\ce_k$ admit special
{\em gauge companion} operators $G_k$ so that the pairs $(L_k,G_k)$
are graded injectively elliptic and, in an appropriate
sense, conformally invariant. 
A general definition of graded injective ellipticity is possible along
the lines of \cite{DougNir}.  A definition more focused on our present
purposes is as follows. Let $P:V\to W$ be a natural differential operator
between bundles which are natural for conformal structure,
and suppose that 
a choice of
scale $g$ naturally splits $V$ as $V_1\oplus\cdots
\oplus V_r$ and $W$ as $W_1\oplus\cdots\oplus W_s$.  Let $P^j_i:V_j\to W_i$
be the block decomposition with respect to this splitting.  Then
$P$ is 
{\em graded injectively elliptic} if there is a positive integer $m$ and 
there are differential operators
$\overline{P}^i_k:W_i\to V_k$, natural for Riemannian structure, with
\begin{equation}\label{injeelli}
\sum_{i=1}^s\overline{P}^i_kP^j_i=\delta^j_k\Delta^m+({\rm order}\,<2m)
\end{equation}
at any Riemannian metric.  
The various $P^i_j$ will generally have different orders.
Graded injective ellipticity in the above sense implies, 
for an appropriately natural operator,
any reasonable graded injective ellipticity property that 
might be formulated in the setting of more general partial differential operators.
An injectively elliptic $P$ has finite-dimensional kernel on a compact
manifold, since the operator $\overline{P}P$ described in \nn{injeelli} does.

By the same token, we can speak of {\em graded surjectively elliptic
operators}.  In the notation above, $P$ is such if there are natural
differential operators $\underline{P}^k_j:W_k\to V_j$ with
\begin{equation}\label{surjell}
\sum_{j=1}^rP^j_i\underline{P}^k_j=\d^k_i\Delta^p+({\rm order}\,<2p)
\end{equation}
for some $p$.

In other signatures, we can still consider the 
properties \nn{injeelli} and \nn{surjell}, but
now the right hand side
 will involve powers of the (pseudo-)Laplacian for that signature.
 In
general signature, we will say that an operator $P$ is
 {\em quasi-Laplacian} if it is a right factor of a power Laplacian in
 the sense of \nn{injeelli}.
Of course, for indefinite metric
 signatures, the quasi-Laplacian property does not naturally relate to
 ellipticity considerations, but rather to properties of hyperbolicity
 or ultra-hyperbolicity.  
Most of the following theorem takes place in the setting of arbitrary metric
signature.

\begin{theorem}\label{main}
On conformal manifolds of even dimension $n$, for each $k=0,1,\cdots ,n/2$
there is a natural
conformally invariant differential operator
$$
\bL_k:\ce^k\to \cg_k
$$ with the following properties.\newline 
\IT{i}\ $q^k\bL_k=L_k$.  In particular $q^k\bL_k$ is 
trivial on the null
space of $L_k$.\newline
\IT{ii}\ $\bL_k$ determines a conformally invariant
operator $G_k:{\Cal N}(L_k)\to \ce_{k-1}$,
which satisfies
$$
(n-2k+4)G_kd=L_{k-1}  .  
$$ 
\IT{iii}\ $G_k=\delta \tilde{M}_{k-1}$,
where
$\tilde{M}_{k-1} : {\Cal N}(L_k)\to \ce_{k}/{\Cal N}(\d)$ is
conformally invariant.    \newline 
\IT{iv}\ For each $k\leq n/2-1$ the operator
$\bL_k$ is quasi-Laplacian, and thus in the case of Riemannian signature
it is a graded injectively elliptic operator.
\end{theorem}

The operator $ \bL_k$ is the $\ell=n/2-k$ special case of
the operator $ \bL_k^\ell$ defined by equation (\ref{bLdef})
below. This will actually make the equation
$q^k \bL_k^\ell = L_k^\ell$ the definition of
$L_k^\ell$.
That these are conformally invariant operators is
established in part \IT{i} of Theorem \ref{Main}.  Part \IT{ii} above
is proved in part \IT{iii} of Theorem \ref{Main}.  This and the fact
that the $ L_k$ are formally self-adjoint enable us to conclude the
result (\ref{mk}) in Theorem \ref{long}. Part \IT{iii} above is part
\IT{iv} of Theorem \ref{Main}, and finally the result \IT{iv} above is
exactly Proposition \ref{bLelliptic}. Note that, from parts \IT{iv,v} 
of Theorem \ref{Main}, $L_{n/2}=0$ and $G_{n/2}$ is a non-zero multiple of $\d$.
At the other extreme of $k$ we note that 
$q^0$ is the identity on $\cg_0\cong \ce_0$ so 
$\bL_0=L_0$ and $G_0$ is the zero operator. 

Note that the injective ellipticity of $\bL_k$ for conformal Riemannian
 structures implies that it has a finite-dimensional conformally
 invariant null space ${\Cal H}^k_\bL $ for compact $M$, and thus this 
null space is a
 candidate for a space of ``conformal harmonics''.  Since the exterior
 derivative is a right factor of $L_k$ {\em ab initio}, 
the space ${\Cal H}^k:=\cN(G_k:\cC^k\to \ce_{k-1})$ is contained in
 ${\Cal H}^k_\bL $ and leads to simpler results which we discuss
 first.  First we state a proposition giving a relationship between
 $\cH^k$ and $H^k(M)$ in the general case.  

\begin{proposition} \label{ESTimate}
On even dimensional conformal manifolds there
is a canonical exact sequence of vector space homomorphisms
$$
0\to H^{k-1}(M)\to H^{k-1}_L(M) \to \cH^k \to H^k(M) 
\quad \mbox{ for } k=1,\cdots ,n/2-1, 
$$
where $\cH^k\to H^k(M)$ is the map taking $w \in \cH^k\subset \cC^k$ to its 
equivalence class 
in $H^k(M)$. Thus in the compact Riemannian case we have
$$
\dim H^{k-1}(M)+\dim \cH^k \leq \dim H^k(M)+\dim H^{k-1}_L(M).
$$
\end{proposition}

\noindent{\bf Proof:} The theorem is clear for $k=0$ since 
by their definitions both
$\cH^0$ and $H^0(M)$
are the space of locally constant functions. Suppose now $k\geq 1$.    
If $ w \in \cH^k$ is mapped to the class of 0 in $H^k(M)$ then 
$w $ is exact. Since in addition $G_k w=0$, it follows from 
Theorem \ref{main} part \IT{ii} that 
$w = d u$ for $u\in \cN(L_{k-1})$. On the other hand recall that by 
Theorem \ref{long},
$d$ is a right factor of $L_{k-1}$ and so $\cC^{k-1}\subseteq \cN(L_{k-1})$. 
Thus there is an exact sequence 
$$
0\to \cC^{k-1}\to \cN(L_{k-1})  \to \cH^k \to H^k(M) 
$$ 
where the map $\cN(L_{k-1}) \to \cH^k$ is the restriction of
exterior differentiation.  Since $H^{k-1}(M)$ is the image of
$\cC^{k-1} $ under the composition $\cC^{k-1}\to \cN(L_{k-1})\to
H^{k-1}_L(M) $, the sequence in the lemma is constructed. 
\quad $\Box$\\

\noindent{\bf Remark:} By a very similar argument one shows that there
 is also a canonical exact sequence of vector space homomorphisms
$$
0\to H^{k-1}(M)\to H^{k-1}_L(M) \to \cH^k_\bL \to H^k_L(M) 
$$
giving $\dim H^{k-1}(M)+\dim \cH^k_\bL \leq \dim H^k_L(M)+\dim H^{k-1}_L(M).$
$\qquad\endrk$

By the proposition, the map $\cH^k \to H^k(M) $ is injective if and only if 
$H^{k-1}(M) =H^{k-1}_L(M)$. In fact in this case it is an isomorphism. 

\begin{theorem} \label{strong}
On compact conformal Riemannian manifolds of even dimension $n$, 
$\cH^{n/2}$ is isomorphic to $H^{n/2}(M)$.  In addition, for each
$k=0,1,\cdots ,n/2-1$, if $H^{k-1}(M)=H^{k-1}_L(M)$ then the
conformally invariant null space $\cH^k$ of $G_k$, acting on $\cC^k$,
is naturally isomorphic, as a vector space, to $H^k(M)$. 
\end{theorem}

\noindent{\bf Proof:} 
The first statement is immediate, 
since $G_{n/2}$ is a non-zero constant multiple of $\d$.

For $k-1<n/2$, as in the proof of Proposition \ref{Hpairing},
at an arbitrary choice of conformal scale we have the Hodge decomposition
\nn{HajAtScale} of $\ce^{k-1}$ or $\ce_{k-1}$.
But under the
assumption of the theorem we have $H^{k-1}(M)=H^{k-1}_L(M)$ and so $\cN(\d
)\cap \cN (L_{k-1}) $ is the usual space of de Rham harmonics $ \cN(\d
)\cap \cN (d)$. Since  by the usual de Rham Hodge decomposition 
$\cR(\d)$ is the 
intersection of $\ce^k$ with the $L^2$ 
orthocomplement of 
$ \cR(d)\oplus (\cN(\d )\cap \cN (d))$, it follows that
$$
\cR (L_{k-1}) =\cR(\d).
$$

Now recall 
from Theorem \ref{main} part \IT{iii} that $G_k=\d \tilde{M}_{k-1}$ on $
\cC^k$. From this and the last display it is immediate that given any
equivalence class $[w]$ in $H^k(M)$ there exists $u\in \ce^{k-1}$
solving the equation
$$
(n-2k+4)G_k w+  L_{k-1} u =0 .
$$ 
By Theorem \ref{main} part \IT{ii} this gives the solution $w'=w+d
u$ to the problem of finding $w'\in [w]$ satisfying $G_k w'=0$.
Thus the injective map $\cH^k\to H^k(M)$ is also surjective.  \quad
$\Box$ 

Since in general we would expect that $H^k_L(M)=H^k(M)$, 
it is worth noting the obvious corollary.
 
\begin{corollary}\label{weak} For each
$k=0,1,\cdots ,n/2-1$,
if $H^{k-1}(M)=H^{k-1}_L(M)$ and $H^{k}(M)=H^{k}_L(M)$, then
$\cH^k_\bL:= \cN(\bL_k)$ 
is naturally isomorphic, as a vector space, to $H^k(M)$.
\end{corollary} 

Note that if $k=n/2-1$ then $L_k$ is the usual Maxwell operator; 
thus the condition 
$H^{k}(M)=H^{k}_L(M)$ is automatically satisfied and $\cH^k_\bL= \cH^k_L(M)$.

Finally in this section we show the operators $M_k$ and $\tilde{M_k}$
above are related to an operator on forms which in an appropriate
sense generalises Branson's Q-curvature.  
Each operator $L_{k-1}$ of Theorem \ref{long} evidently only
determines $M_{k-1}$ as an operator $\cC^{k}\to \ce_{k}/\cN(\d)$,
while $\tilde{M}_{k-1}$ is similarly fixed only as an operator
$\cN(L_k)\to \ce_{k}/\cN(\d)$. One might hope that these conformally
invariant operators are the restrictions of some conformally invariant
operator or operators $\ce^k\to \ce_k$. In fact this is impossible. 
The invariant
differential operators on the standard conformally flat model are
classified via the bijective relationship with generalised Verma
module homomorphisms and the classification of the latter in
\cite{BoeColl2}. (See \cite{EastSlo} and references therein for a
summary of the relevant representation theoretic results and details on how
these correspond dually to a classification of 
invariant differential operators.)  
{}From this well known classification, 
it is clear that on such structures the $L_k$ are (up to
constant multiples) unique. Since $L_k d=0=\d L_k$ these are not
candidates for the $M_{k-1}$ or the $\tilde{M}_{k-1}$.  Expression
\nn{Qdef} of Section \ref{Qsect} defines, on even dimensional
conformal manifolds, for each choice of conformal scale $\si$ and for
$k=0,1,\cdots n/2$, an operator
$$
Q^\si_k:\ce^k\to \ce_k .
$$ Parts \IT{i,ii} of the next theorem assert that the
operators $M_{k-1}$ and $\tilde{M}_{k-1}$ are (up to a multiple)
restrictions of $Q^\si_k$. Parts \IT{iii--v} show
that as operators on closed forms the $Q^\si_k$ generalise
Branson's Q-curvature.

We should point out that the construction \nn{Qdef} defines each
$Q^\si_k$ as a composition of tractor operators arising from the
ambient construction, but via \nn{bGdef} and the theory of Section
\ref{ambext}, this may be readily re-expressed as a universal polynomial
expression in the covariant derivative $\N$ and the Riemann curvature
$R$.

\begin{theorem}\label{Qthm} 
\noindent \IT{i} In a conformal scale $\si$ the operator 
$$
Q^\si_k:\ce^k\to \ce_k
$$
is formally self-adjoint. \newline
\noindent  \IT{ii} As an operator on $\cN(L_k)$,
$
\d Q_k^{{\si}}
$
is conformally invariant and
$$
\d Q_k^{{\si}}= G_k.
$$
\noindent  \IT{iii} Operating on
$
\ce^{k-1}
$,
we have
$$
(n-2k+4)\d Q_k^\si d= L_{k-1} .
$$
\noindent  \IT{iv}  As an operator on closed $k$-forms
$$
Q_k^{\si}:\cC^k \to \ce_k,
$$
$Q_k^{\si} $ has the conformal transformation law
$$
Q_k^{\hat{\si}}u=Q_k^\si u + L_k (\Upsilon u)
$$
where $ \Upsilon$ is a smooth function and 
$ \hat{\si}=e^{-\Upsilon}\si$. \\
\noindent  \IT{v}
$ Q^{\si}_0 1$ is the Branson Q-curvature.
\end{theorem}

This theorem is proved in the last part of Section \ref{Qsect}. Note
that from part \IT{ii} above we see that acting on $\cN(L_k)$,
$Q^\si_k$ gives the operator $\tilde{M}_{k-1}$ of Theorem \ref{main}.
On $\cC^k\subseteq \cN (L_k)$, we have $(n-2k+4) Q^\si_{k} = M_{k-1}
$.  That this is conformally invariant as an operator $\cC^k\to \cC_k$
is immediate from part \IT{iv} since $L_k: \ce^k\to \cR(\d)\subset
\ce_k $.  From our observations concerning $M_{k}$ after 
Corollary \ref{weak}
it is clear that, as an
operator on $\ce^k$ with $k<n/2$, $Q_k^\si$ is {\bf not} of the form
$(\mbox{conformally invariant operator})+ \d U +V d $, where $U$ and
$V$ are differential operators.  Using the tools we
develop below, it is easy to verify that $Q^\si_{n/2}$ is a multiple
of the identity.

A celebrated property of the Q-curvature is that its
integral is conformally invariant. In the next result we observe that
there is a somewhat stronger invariance result, in that one
can integrate invariantly 
against the null space $\cN(L_0)$.  This property
is generalised by the operators $Q^\si_k$.

\begin{theorem} \label{pairings}
\IT{i} As an operator between $\cC^k$ and $\ce_k/\cR (L_k) $,
$Q_k^{\si}$ is conformally invariant.  Further restricting to
$\cH^k\subset \cC^k$, we obtain a conformally invariant operator 
$$
Q_k: {\Cal H}^k\to H_{k}^L(M) , 
$$ 
where $ {\Cal H}^k:=\cN(G_k: \cC^k
\to \ce_{k-1})$. \\ \
\IT{ii} On compact manifolds,
$Q_k^{\si}$ gives a conformally invariant pairing between ${\Cal
N}(L_k)$ and $\cC^k$ by 
$$
(u,w)\mapsto \langle u, Q_k w \rangle
$$
for $w\in\cC^k $ and $u\in \cN(L_k)$. 
On compact conformal manifolds, 
the same integral formula
determines a pairing between $H^k_L(M)$ and $\cH^k$
by taking $w\in\cH^k $ and $u$ any representative of $[u]\in H_{k}^L(M)$.\\
\end{theorem}

\noindent{\bf Proof:} The first statement is a trivial
consequence of part \IT{iii} of the previous theorem.  Now suppose
that $u\in {\Cal H}^k$.  By part \IT{ii} of Theorem \ref{Qthm},
$\d Q^\si_k u=G_k u =0 $.  
So the conformally invariant map $Q_k^{\si}:\cC^k \to \ce_k/\cR
(L_k) $ descends to a well-defined map $Q_k: {\Cal H}^k\to
H_{k}^L(M)$. (In view of the conformal invariance  we
omit the argument $\si$ from $Q_k$ here.) This establishes
part \IT{i}.  For part \IT{ii},
the first statement follows from the pairing \nn{pair}, the
first statement in part \IT{i} and the fact that $L_k$ is formally
self-adjoint. The second can be deduced from this, or follows
immediately from part \IT{i} and the earlier observation that
\nn{pair} induces a conformally invariant pairing between $H^k_L(M)$ and
$H_k^L(M)$.
\quad $\Box$\\

Recall that, from the factorisation $L_k=\d M_k d$, there is a natural
 (conformally invariant) surjection $\ce_k/\cR(L_k)\to \cC_k$ inducing
 the map $H^L_k(M)\to H_k(M)$. Thus from part \IT{i} above, $Q^\si_k$
 induces a conformally invariant map $Q_k:\cC^k\to \cC_k$. Since
 $\cC^k\subseteq \cN(L_k)$, $Q^\si_k$ induces a conformally invariant
 pairing of $\cC^k$ with itself by restriction of the pairing in part
 \IT{ii} above.  (In the compact Riemannian setting the latter is
 equivalent to the map $Q_k:\cC^k\to \cC_k$.)  Similarly the
 composition of the displayed map in part \IT{i} with $H^L_k(M)\to
 H_k(M)$ gives a conformally invariant map into de Rham cohomology,
 $\cH^k\to H_k(M)$, and the pairing just described descends to a
 conformally invariant pairing between $H^k(M)$ and $\cH^k$.  We
 summarise this in:

\begin{corollary}\label{corpairings} $Q_k$ induces conformally invariant maps
$\cC^k\to\cC_k$ and $\cH^k\to H_k(M)$. 
In the compact setting, $Q_k$ induces conformally invariant
 pairings of $\cC^k$ with itself, and of $H^k(M)$
and $\cH^k$.\end{corollary}

As an application of these results we can now show that on 
compact Riemannian manifolds 
$\dim
\cH^k\geq \dim H^k(M)$. 
Combining the map $\cH^k\to H_k(M)$ just
discussed with the composition of $H^k_L(M)\to H^k(M)$ and the map $Q_k:
\cH^k\to H^L_k(M)$ of Theorem \ref{pairings} part \IT{i}, we obtain a
conformally invariant map $ I:\cH^k\to H^k(M)\oplus H_k(M) $ given by
$$
\w \mapsto ([\w],[Q^\si_k \w]) .
$$
Clearly $\w$ is killed by this map if and only if $\w$ is 
both exact and $Q^\si_k \w$ is in 
$\cR(\d)$. Thus the space 
$$
\cB^k:=\{d\ph\suchthat Q_k^\si d\ph\in \cR(\d) \}
$$
is conformally invariant and
$$
I:\cH^k/\cB^k \to H^k(M)\oplus H_k(M)
$$ is injective. It turns out that the domain space here has the same
dimension as $H^k(M)$. To show this we explicitly identify a space which,
via the pairing \nn{pair}, is its vector space dual. Let
$$
\cH_k:= \{ \xi\in \ce_k \suchthat \d \xi=\d Q_k^\si \eta 
~\mbox{for some } \eta\in \cC^k \},
$$
and
$$
\cB_k:= 
\{ \xi \in \ce_k \suchthat \xi-Q_k^\si d \ph\in \cR(\d) 
~\mbox{for some }\ph\in \ce^{k-1}\}.
$$ The space $\cH_k$ is conformally stable since from Theorem
 \ref{Qthm} part \IT{ii}, $\d Q^\si_k = G_k$ is conformally invariant
 on $\cN(L_k)\supseteq \cC^k$.  The conformal invariance of $\cB_k$
 is immediate from part \IT{iii} of the same theorem.

Next let us fix a conformal scale $\si$. Then we have the map
$$
P: \cC^k\oplus \cN(\d)\to \cH_k\  
\mbox{ given by }\ 
(\eta,\xi)\mapsto \xi - Q_k^\si \eta .
$$ 
This is surjective since for $\xi\in \cH^k$ there is $\eta\in
\cC^k$ such that $\d \xi=\d Q_k^\si \eta$, and so
$(-\eta,\xi-Q^\si_k \eta)$ is a pre-image of $\xi$. Now if $\eta= d
\ph$ and $\xi =\d \rho$ then $P(\eta,\xi)= \d \rho -Q^\si_k d \phi $
which is in $\cB_k$. Thus $P$ descends to a well-defined surjective map 
$H^k(M)\oplus H_k(M)\to \cH_k/\cB_k$.
 It is clear that $\ker(P)=\Im (I)$ so, in summary, in a conformal
scale we have an exact sequence
\begin{equation} \label{exact}
0\to \cH^k/\cB^k \stackrel{I}{\to} H^k(M)\oplus H_k(M) \stackrel{P}{\to} 
\cH_k/\cB_k \to 0,
\end{equation}
from which it is clear that $\dim (\cH_k/\cB_k)$ is finite. 
The following result shows that 
$\dim (\cH_k/\cB_k)=\dim (\cH^k/\cB^k)=b^k:=\dim H^k(M)$.

\begin{theorem}\label{perfecttoo}
For compact even dimensional Riemannian conformal manifolds and each 
$k=0,1,\cdots ,n/2-1$,
the conformally invariant pairing between $\cH^k$ and $\cH_k(M)$
given  by the restriction of \IT{\ref{pair}}
descends to a well-defined conformally
invariant perfect pairing of $\cH^k/\cB^k$ with $\cH_k/\cB_k$.
\end{theorem} 
\noindent{\bf Proof:} The conformal invariance of the pairing between
$\cH^k$ and $\cH_k$ is immediate from the invariance of the pairing
\nn{pair}. It is clear that if this descends as claimed then the
result is an invariant pairing. Note that for $k=0$ we have by
construction $\cH^0=H^0(M)$ and $\cH_0=H_0(M)$ 
and so the result holds trivially.  
Henceforth we assume $k\geq 1$ and fix a conformal scale
$\si$.

Note that if 
$d\ph\in \cB^k$ and
 $\xi\in \cH_k$ we have
$$ 
\langle d\ph , \xi \rangle= \langle \ph,\d
\xi \rangle= \langle \ph,\d Q_k^\si \eta \rangle= \langle Q_k^\si d
\ph , \eta \rangle= \langle \d \mu,\eta \rangle= \langle\mu,d \eta
\rangle= 0,
$$ 
where we have used that $Q_k^\si$ is formally self-adjoint.
Thus via the pairing, $\cB^k$ annihilates $\cH_k$. Now
suppose that $\w\in \cH^k$ and $\xi\in \cB_k$. Then there is a pair 
$(\ph,\rho)\in\ce^{k-1}\oplus \ce_{k+1}$ so that
$$ \langle \w, \xi \rangle= \langle \w, Q_k d\ph -\d \rho \rangle =
 \langle \d Q^\si_k\w, \ph\rangle - \langle d \w,\rho \rangle=0,
$$ 
since $w\in \cH^k$ means that $d\w=0$ and $\d
Q^\si_k\w=G_k \w=0$. So $\cB_k$ annihilates $\cH^k$ and 
$\langle\cdot,\cdot \rangle$ descends to a bilinear function
on $(\cH^k/\cB^k)\times(\cH_k/\cB_k)$ as claimed.
It remains to show the pairing
is perfect. Suppose that $\w\in \cH^k$ satisfies $\langle
\w,\xi\rangle=0$ for all $\xi\in \cH_k$. Then in particular $\langle
\w, \xi \rangle=0$ for any $\xi\in \cN(\d)\subset \ce_k$. Thus by
Poincar\'{e} duality $\w=d \ph$ for some $\ph\in \ce^{k-1}$. Now
consider $\langle \eta,Q_k \w \rangle$ where $\eta\in \cC^k$.  Note
that $Q_k \eta \in \cH_k$ and so, by the assumption, $\langle
\w,Q^\si_k \eta \rangle=0$. But $Q^\si_k$ is formally self-adjoint and
this implies $\langle \eta, Q^\si_k \w \rangle=0$. Again by
Poincar\'{e} duality, since $\eta$ is an arbitrary element of $\cC^k$,
this implies $Q_k^\si \w\in \cR(\d)$. So the pair $(\w,Q^\si_k \w)$
vanishes in $H^k(M)\oplus H_k(M)$.
 
Finally suppose that $\xi\in \cH_k$ satisfies $\langle
\w,\xi\rangle=0$ for all $\w\in \cH^k$. Observe that by standard de
Rham Hodge theory $d \w =0 \Leftrightarrow \fl^{\ell}d\w=0$ where
$\fl$ is the form Laplacian $\d d+d \d$ and 
$\ell= n/2-k$. This gives an alternative description of $\cH^k$, 
$$
\cH^k=\{\w\in \ce^k \suchthat \fl^{\ell}d\w=0 
~\mbox{ and  }~\d Q^\si_k \w=0 \} .
$$ So by the assumption $\xi$ is orthogonal to the kernel in $\ce^k$
of the operator pair $(\d Q^\si_k,\fl^{\ell}
d):\cE^k\to\cE_{k-1}\oplus \cE_{k+1}$ 
(ignoring conformal weights as we may since we have
fixed a conformal scale).  
Put another way, $\xi$ is
orthogonal to kernel of the adjoint of the operator
\begin{equation}\label{surj}
\begin{array}{r}
(Q^\si_k d,\d\fl^\ell):\ce^{k-1}\oplus\ce^{k+1} \to \ce^k, \\  
\renewcommand{\arraystretch}{0.8}
\left(\begin{array}{c}\ph \\ \mu\end{array}\right) 
\mapsto Q^\si_k d \ph +\d\fl^{\ell}\mu .
\end{array}
\end{equation}
\renewcommand{\arraystretch}{1.1}
But from Theorem \ref{Qthm} part \IT{iii} it follows that the
leading term of $ Q_k^\si d$ is, up to a non-vanishing constant
multiple, $(d\d)^\ell d=\fl^\ell d$. Thus there is a number $\alpha$
so that 
$$
\left(Q^\si_kd , \d\fl^{\ell}\right)\left(\begin{array}{c}
                                                \alpha \d\\
                                                 d 
\end{array}\right) =\fl^{\ell+1} +\LOT: \ce^k\to \ce^k.
$$ 
(Here and below ``$\LOT$'' is an abbreviation for ``lower order terms''.)
Thus the operator \nn{surj} is surjectively elliptic, and in the
resulting Hodge decomposition
\renewcommand{\arraystretch}{1.3}
$$
\cE^k=\cR(Q^\si_kd,\d\fl^\ell)\oplus\cN
\left(\begin{array}{c}\d Q^\si_k\\ \fl^\ell d
\end{array}\right),
$$ 
\renewcommand{\arraystretch}{1.5} 
$\xi$ must be in the first
summand. Thus $\xi \in \cB_k$.  \quad $\Box$ 

By the finite dimensionality of the spaces $\cH^k/\cB^k$ and
$\cH_k/\cB_k$, the duality established above and the exact sequence
\nn{exact} we have the following.

\begin{corollary} \label{bk}
For any compact even dimensional Riemannian conformal 
manifold and $k= 0,1,\cdots ,n/2-1$ we have 
$$
\dim (\cH^k/\cB^k)=b^k=\dim (\cH_k/\cB_k). 
$$
Thus 
we have
$$
\dim\cH^k=b^k+\dim\cB^k\quad\mbox{ and }\quad\dim\cB^k\le
\dim(H^{k-1}_L(M)/H^{k-1}(M)).
$$
\end{corollary}

Note that the final conclusion is clear, since from the definition
of $\cB^k$ and part \IT{iii} of Theorem \ref{Qthm},
$d\phi\in\cB^k\Rightarrow L_{k-1}\phi=0$.

{}From Theorem \ref{perfecttoo} we see that via \nn{pair} $\cH_k/\cB_k$
may identified with the vector space dual to $ \cH^k/\cB^k$. 
It is straightforward to verify that via \nn{pair} $H_k(M)$ 
may be similarly identified 
with the vector space dual to $H^k(M)$. Thus we have the following result. 

\begin{corollary}\label{regular} 
For any even
dimensional compact Riemannian conformal manifold and 
$k= 0,1,\cdots ,n/2-1$, we have
$$
\begin{array}{rrcl}
&\cH^k/\cB^k\to H^k(M) \mbox{ is injective } & \Leftrightarrow & 
\cH^k/\cB^k\to H^k(M) \mbox{ is surjective }\\
\Leftrightarrow & H_k(M)\to \cH_k/\cB_k  \mbox{ is injective } &   
\Leftrightarrow & H_k(M)\to \cH_k/\cB_k  \mbox{ is surjective.} 
\end{array}
$$
\end{corollary}

Theorem \ref{perfecttoo}, Corollary \ref{bk} and Corollary \ref{regular}
generalise, respectively, Theorem 4.1, Corollary 4.2 and Corollary
4.3, of Eastwood and Singer \cite{EastSin2}, which deal with 1-forms
in dimension 4. Our treatment of these last three results has been
heavily influenced by their development of that case. For each $k$,
the (equivalent) conditions of Corollary \ref{regular} constitute some
conformally invariant condition on the Riemannian conformal manifold, or
Riemannian manifold, that we term $(k-1)$-{\em regularity}. 
This generalises the notion of `regular' for $k=1$, $n=4$ discussed in
\cite{EastSin2}. Similarly the conformally invariant hypothesis
$H^{k}(M)=H^{k}_L(M)$ of Theorem \ref{strong}, which we term {\em
strong $k$-regularity}, generalises the dimension 4 notion of `strong
regularity' in \cite{Sin}.

Of course if $H^{k+1}(M)$ 
vanishes, then clearly $M$ is $k$-regular; this is the trivial
case. 
We expect that $k$-regularity, for each $k$, 
should hold generically in some appropriate sense for 
compact conformal Riemannian
manifolds. 
Note that in
$\ce^k$ we have the subspace inequality 
$$
\cB^{k+1}=\{d\ph\suchthat Q_{k+1} d\ph\in \cR(\d) \} \subseteq
\{d\ph\suchthat L_{k}\ph =(n-2k+2) \d Q_{k+1} d\ph =0 \} .
$$ 
By 
Proposition \ref{ESTimate}, $k$-regularity is the case of equality 
in this inequality,
while strong $k$-regularity is the assertion that the
space on the right side vanishes.  In particular, strong $k$-regularity
implies $k$-regularity.  In dimension 4, \cite{EastSin2} shows that
Einstein manifolds are $0$-regular,
and also gives an example of a
manifold which fails to be strongly $0$-regular.

\subsection{Extensions to the theory and non-orientable manifolds} 
\label{codetoursect}

No assumptions have been made above concerning the orientability of
$M$. In the general case that $M $ may have non-orientable components
there is further information to be extracted via an extension to the
theory.  We present this rather concisely since at one level this
extension arises rather simply from the machinery described above and
given in the following sections.

In the case that $M$ is orientable the conformal  
Hodge star operator gives an isomorphism
$\star: \ce^k[w] \stackrel{\cong}{\to} \ce^{n-k}[w+n-2k]=\ce_{n-k}[w]$. 
Since up to a sign 
$\star \star$ is the identity, it follows that there are operators
$
L^{n-k}_{\ell,\star}:\ce_{n-k}[w]\to \ce^{n-k}[-w]
$
where 
\begin{equation}\label{Lstar}
L^{n-k}_{\ell,\star} = \star L_k^{\ell} \star  \quad  \mbox{locally}.
\end{equation}
We have written ``locally'' since we also define the operators
$L^{n-k}_{\ell,\star} $ on non-orientable manifolds by requiring that
\nn{Lstar} hold on every orientable neighbourhood, for any choice of
orientation on a given neighbourhood. Since the local choice of
orientation only affects the sign of $\star$ it is clear that each
$L^{n-k}_{\ell,\star}$ is well defined. When $M$ is orientable the
operators $L^{n-k}_{\ell,\star} $ are, by construction, equivalent to
the $L^\ell_k$. Otherwise there need not be an isomorphism 
between $\ce^k$
and $\ce_{n-k}$ and so these are new formally self-adjoint conformally
invariant operators. These are non-trivial and natural for $k$ and
$\ell$ as in Theorem \ref{long}.  When $w:=k+\ell-n/2=0$ we use the
alternative notation $L^{n-k}_{\star}:= L^{n-k}_{\ell,\star}$ and we have
the following result.

\begin{theorem}\label{colong}
If $n$ is even and $w=0$, we have that 
\begin{equation*}
L^{n-k}_{\star} =d M^{n-k}_{\star} \d,
\end{equation*}
where, up to a constant multiple, $M^{n-k}_{\star}$ is given locally 
{\rm(}and in a choice
of conformal scale $\si${\rm)} by $\star Q_{k+1}^\si \star $.

On Riemannian conformal manifolds the following holds.
The operator $L^{n-k}_{\ell,\star}$ is elliptic
if and only if  $k\neq n/2\pm \ell$,
and it is positively elliptic if and only if $k\notin[n/2-\ell,n/2+\ell]$.
For each $k$ the differential operator sequence
\begin{equation}\label{codetour}
\ce_n\stackrel{\d}{\to}\cdots\stackrel{\d}{\to}\ce_{n-k+1}\stackrel{\d}{\to}
\ce_{n-k}\stackrel{L^{n-k}_{\star}}{\to}\ce^{n-k}\stackrel{d}{\to}\ce^{n-k+1}
\stackrel{d}{\to}
\cdots\stackrel{d}{\to}\ce^n
\end{equation}
is an elliptic complex.
\end{theorem}

\noindent{\bf Proof:} On orientable neighbourhoods we have
$L^{n-k}_{\star} = 2(\ell+1)\star \d Q^\si_{k+1} d \star$, from Theorem
\ref{Qthm}.  Recalling also that as operators on ${k'}$-forms we have
$\star d=(-1)^{{k'}+1}\d \star$ and $d\star =(-1)^{k'}\star \d$ the first
result is proved. 
The remaining facts follow immediately from the corresponding results
for the $L_k$.  \quad $\Box$

That the $k^{th}$ {\em de Rham
codetour complex} \nn{codetour} is not equivalent to the detour
complex \nn{detour} is already clear by taking the cohomology at $0\to
\ce_n\stackrel{\d}{\to} $. 

\begin{proposition}\label{Hn}
On compact Riemannian conformal manifolds,\newline
$\dim H_n(M)$ is the number of oriented components in $M$.
\end{proposition}

\noindent{\bf Proof:} Suppose
that, on a connected component $M'$ of $M$, 
$\phi$ is a non-zero section of $\ce_n$ annihilated by
$\d$. Pick a metric $g$ from the conformal class and identify $\ce_n$
with $\ce^n$ on $M'$.  Consider an arbitrary orientable neighbourhood
$U$
in $M'$, let $V^g$ be a choice of volume form on $U$ consistent with
$g$, and observe that $\phi = a V^g$ for some function $a$. 
Since $\d V^g=0$ and $V^g$ is co-closed, it follows that $\i(d a)V^g=0$. 
Thus $a$ is constant on $U$.
Now $V^g$ is preserved by the Levi-Civita connection
and so $\phi$ is covariantly constant on $U$.
Since $U$ was arbitrary, it follows that $\phi$
is covariantly constant and so nowhere vanishing on $M'$. Thus $M'$ is
orientable.  \quad $\Box$ 

As a result, $\dim H_n(M)$ is in general less than
$\dim H^0(M)$. Note that from the perfect pairing between $H_n(M)$ and
$H^n(M)$ the proof above recovers the well known result that $ \dim
H^n(M)$ is the number of connected components.

The main point now is that there is a theory for the operators
$L^{n-k}_\star$ which closely parallels the theory for the operators
$L_k$. For example for $k\leq n/2-1$ each $L^{n-k}_\star$ 
has an extension to a quasi-Laplacian operator
$\bL^{n-k}_\star:\ce_{n-k}\to
\cg^{n-k}_\star$ where $\cg^{n-k}_\star$ is the section space of a
bundle with a composition series $E^{n-k} \lpl E^{n-k+1}$.  There is
an analogue for Theorem \ref{main}. In this, for example, (ignoring
non-zero constant scalar multiples) the analogue of $G_k$ is $d
Q_{n-k}^{\si,\star}$ where, for each choice of conformal scale $\si$,
$Q_{n-k}^{\si,\star}$ is the unique operator given on oriented
neighbourhoods by $\star Q^\si_k \star$. $Q_{n-k}^{\si,\star}$ is
another type of Q-operator and satisfies the obvious analogue of
Theorem \ref{Qthm} parts \IT{i--iv}. The existence of the
bundle $\cg^{n-k}_\star$, and the other local results we have
mentioned here, follow trivially from earlier results since the
operators $L^{n-k}_\star$ and $L_k$ are, by construction, locally
equivalent.  On the other hand it is also straightforward to directly
develop these results using the calculus in the later sections and
some straightforward variations on the constructions there. See in particular
the remarks on pages \pageref{cosubtractors} and \pageref{coK} where
the key tools are described.

Of course when $M$ is not orientable there is no result corresponding
to part \IT{v} of Theorem \ref{Qthm} and this is an important
distinction between the cases at a global level.  There are obvious
analogues for all the cohomological theorems. We may define
$H_{n-k}^{L_\star}(M)$ and $H^{n-k}_{L_\star}(M)$ for the cohomology of
\nn{codetour} at, respectively, the bundles $\ce_{n-k}$ and
$\ce^{n-k}$. On the other hand we may define the space $\cH_{n-k}^\star$
of conformal harmonics to be the null space of $d
Q_{n-k}^{\si,\star}$ as an operator on $\cN(\d:\ce_{n-k}\to
\ce_{n-k+1})$. These satisfy analogues of Proposition \ref{Hpairing},
Proposition \ref{ESTimate}, Theorem \ref{strong}, Theorem
\ref{perfecttoo} and the corollaries of the latter. The Q-operator
$Q_{n-k}^{\si,\star}$ satisfies the analogue of Theorem
\ref{pairings}.

\section{The ambient construction and tractor calculus}\label{ambient}

The basic relationship between the Fefferman-Graham ambient metric
construction and tractor calculus is described in \cite{CapGoFG}. We
review this briefly and establish our notation before developing an
exterior calculus for the ambient manifold and for tractor fields.

Let $\pi:\cq\to M$ be a conformal structure of signature $ (p,q)$.
Let us use $\rho $ to denote the ${\Bbb R}_+$ action on $ \cq$ given
by $\rho(s) (x,g_x)=(x,s^2g_x)$.  An {\em ambient manifold\/} is a
smooth $(n+2)$-manifold $\aM$ endowed with a free $\Bbb R_+$--action
$\rho$ and an $\Bbb R_+$--equivariant embedding
$i:\cq\to\aM$.  We write $\X\in\frak X(\aM)$ for the fundamental field
generating the $\Bbb R_+$--action, that is for $f\in C^\infty(\aM)$
and $ u\in \aM$ we have $\X f(u)=(d/dt)f(\rho(e^t)u)|_{t=0}$.

If $i:\cq\to\aM$ is an ambient manifold, then an {\em ambient
metric\/} is a pseudo--Riemannian metric $\h$ of signature $(p+1,q+1)$
on $\aM$ such that the following conditions hold:

\smallskip
\noindent
(i) The metric $\h$ is homogeneous of degree 2 with respect to the
$\Bbb R_+$--action, i.e.\ if $\Cal L_{\sX}$
denotes the Lie derivative by $\X$, then we have $\Cal L_{\sX}\h=2\h$.
(I.e.\ $\X$ is a homothetic vector field for $h$.)

\noindent
(ii) For $u=(x,g_x)\in \cq$ and $\xi,\eta\in T_u\cq$, we have
$\h(i_*\xi,i_*\eta)=g_x(\pi_*\xi,\pi_*\eta)$. 

\noindent
To simplify the notation  we will usually identify $\cq$ 
with its image in $\aM$ and suppress
the embedding map $i$.
\smallskip
To link the geometry of the ambient manifold to the underlying
conformal structure on $M $ one requires  further conditions. In
\cite{FGast} Fefferman and Graham treat the problem of constructing a formal power
series solution along $ \Cal Q$ for the Goursat problem of finding
an ambient
metric $ \h$ satisfying (i) and (ii) and the condition that it be
Ricci flat, i.e.\ Ric$(\h)=0$. A key result is Theorem~2.1 of their
paper: If $n$ is odd, then up to a $ \Bbb R_+$-equivariant
diffeomorphism fixing $ \Cal Q$, there is a unique power series
solution for $ \h$ satisfying (i), (ii) and Ric$(\h)=0$.  If $ n$ is
even, then up to a $ \Bbb R_+$-equivariant diffeomorphism fixing $
\Cal Q$ and the addition of terms vanishing to order $ n/2$, there is
a unique power series solution for $ \h$ satisfying 
\renewcommand{\arraystretch}{1}
$$
\begin{array}{l}
\mbox{(i), (ii);} \\
\mbox{Ric}(\h)\mbox{ vanishes to order }n/2-2\mbox{ along }\Cal Q; \\
\mbox{tangential components of Ric}(\h)\mbox{ vanish to order }n/2-1
\mbox{ along }\Cal Q.
\end{array}
$$
\renewcommand{\arraystretch}{1.5}
It turns out that in metrics satisfying these conditions $ Q:=\h(\X,\X)$
is a defining function for $ \cq$ and $2\h(\X,\cdot)=dQ$ to all orders in
odd dimensions and up to the addition of terms vanishing to order
$n/2$ in even dimensions. It is straightforward to show
\cite{GoPet,GJMS} that one can extend the solution slightly in even
dimensions to obtain
\renewcommand{\arraystretch}{1}
$$
\noindent {\rm (iii)}  \quad 
\Ric(\h)=0\; \left\{\begin{array}{l} \mbox{to all orders if $n$ is odd,} \\
                                 \mbox{up to the addition of terms vanishing} \\
                                 \qquad\mbox{to order $n/2-1$ if $ n$ is even,}
\end{array} 
\right.
$$ 
with (i), (ii) and $\h(\X,\cdot)=\frac{1}{2}d Q$ to all orders in
both dimension parities.  Henceforth, unless otherwise indicated, the
term ambient metric will mean an ambient manifold with metric
satisfying all these conditions.  We should point out that we only use
the existence part of the Fefferman-Graham construction. The
uniqueness of the operators we will construct is a consequence of the
fact that they can be uniquely expressed in terms of the underlying
conformal structure as we shall later explain. Finally we note that if
 $M$ is locally conformally flat then there is a canonical
solution, to all orders, to the ambient metric problem. This is the
flat ambient metric. This is forced by
(i--iii) in odd dimensions (see \nn{Wform} and the proof of 
Lemma \ref{GoPet4.4}). 
But in even dimensions this extends the solution.
When discussing the conformally flat case we assume this solution. 

We write $ \nda $ for the ambient Levi-Civita connection determined by
$ \h$ and use upper case abstract indices $A,B,\cdots $ for tensors on
$ \aM$. For example, if $ v^B$ is a vector field on $\aM $, then the
ambient Riemann tensor will be denoted $\aR_{AB}{}^C{}_{D}$ and defined
by $ [\nda_A,\nda_B]v^C=\aR_{AB}{}^C{}_{D}v^D$. In this notation the
ambient metric is denoted $ \h_{AB}$ and with its inverse this is used
to raise and lower indices in the usual way. We will not normally
distinguish tensors related in this way even in index free
notation; the meaning should be clear from the context. Thus for example
we shall use $\X$ to mean both the Euler vector field $\X^A$ and the
1-form $ \X_A=\h_{AB}\X^B$.

The condition $\LX\h=2\h$ is equivalent to
the statement that the symmetric part of $\nda\X$ is $\h$.  On the other
hand, since $\X$ is exact, $\nda\X$ is symmetric.  Thus
\renewcommand{\arraystretch}{1.5}
\begin{equation} \label{ndaX} 
\nda\X =\h,
\end{equation} 
which in turn implies 
\begin{equation} 
\label{XRrel} 
\X\hook\aR=0. 
\end{equation} 
Equalities without qualification, as here, indicate that the results
hold to all orders or identically on the ambient manifold.

\subsection{Tractor bundles}\label{tractor}

Let $ \cce(w)$ denote the space of functions on $ \aM$ which are
homogeneous of degree $ w\in {\Bbb R}$ with respect to the action $
\rho$. That is $f\in \cce(w)$ means that $ \X f=wf$. Similarly a
tensor field $F$ on $ \aM$ is said to be {\em homogeneous of degree} 
$w$ if $\rho(s)^* F= s^w F$ or equivalently $ \cL_{\miniX} F=w F$.  Just
as sections of $ \ce[w]$ are equivalent to functions in $
\cce(w)|_\cq$ we will see that the restriction of homogeneous tensor
fields to $\cq$ have interpretations on $M$.

 On the ambient tangent bundle $T\aM$ we define an action of $\Bbb
R_+$ by $s\cdot \xi:=s^{-1}\rho(s)_\ast \xi$. The sections of $ T\aM$
which are fixed by this action are those which are homogeneous of
degree $ -1$. Let us denote by $ \act$ the space of such sections and
write $\act(w)$ for sections in $\act\otimes \cce(w)$, where the
$\otimes$ here indicates a tensor product over $\cce(0)$.  
Along $ \cq$ the
$\Bbb R_+$ action on $T\aM$ is compatible with the ${\Bbb R}_+$ action
on $\cq$, so defining $\bT$ to be the quotient $(T\aM|_\cq)/\Bbb R_+$,
yields a rank $ n+2$ vector bundle $ \bT$ over $\cq/\Bbb R_+=M$.  By
construction, sections of $p:\bT \to M$ are equivalent to sections
from $\act|_\cq$. We write $\ct$ to denote the space of such sections.

Since the ambient metric $\h$ is homogeneous of degree $2$
it follows  that for vector fields $\xi$ and $\eta$ on $\aM$
which are homogeneous of degree $-1$, the
function $\h(\xi,\eta)$ is homogeneous of degree $0$ 
and thus descends to a smooth function on
$M$. Hence $\h$ descends to a smooth bundle metric $h$ of signature
$(p+1,q+1)$ on $\bT$.

Next we show that the space $\act$ has a filtration reflecting the geometry of
$\aM$.  First observe that for $\phi\in \cce(-1)$, $ \phi \X\in
\act$. Restricting to $\cq$ this determines a canonical inclusion $
E[-1] \hookrightarrow \bT $ with image denoted by $\bV$. 
Since $ \X$ generates the fibres of $\pi:\cq\to M$ the smooth
distinguished line subbundle $ \bV\subset \bT$ reflects the
inclusion of the vertical bundle in $ T\aM|_\cq$. 
We write $X$ for the canonical section in $\ct[1]$ giving this
inclusion.  We define $ \bF$ to be the orthogonal complement of $ \bV$
with respect to $ h$.  Since $Q=\h(\X,\X)$ is a defining function for
$\cq $ it follows that $ X$ is null and so $\bV\subset \bF$. Clearly $
\bF$ is a smooth rank $ n+1$ subbundle of $ \bT$. Thus $ \bT/ \bF$ is
a line bundle and it is immediate from the definition of $ \bF$ that
there is a canonical isomorphism $E[1]\cong \bT/ \bF $ arising from
the map $\bT\to E[1]$ given by $V\mapsto h(X,V)$. Now recall
$2\h(\X,\cdot)=d Q$, so the sections of $\act$ corresponding to
sections of $\bF$ are just those that take values in $ T\cq\subset
T\aM|_\cq$. 
Finally we note that if
$\tilde \xi$ and $\tilde{\xi}'$ are two lifts to $ \cq$ of $\xi\in
{\frak X}(M)$ then they are sections of $T\cq$ which are homogeneous
of degree 0 and with difference $\tilde \xi - \tilde{\xi}' $ taking
values in the vertical subbundle.  Since $\pi: \cq \to M$ is a
submersion it follows immediately that $\bF[1]/\bV[1]\cong TM\cong
T^*M[2]$ 
(where recall by our conventions $\bF[1]$ means $\bF\otimes E[1]$ etc.\/). 
Tensoring this with $E[-1]$ and combining this observation 
with our earlier results we can summarise the filtration of $\bT$ by
the composition series
\begin{equation} \label{trcomp} 
\bT= E[1]\lpl T^*M[1] \lpl E[-1] .
\end{equation}

Next we show that the Levi-Civita connection $\nda$ of $\h$ determines
a linear connection on $ \bT$. Since $\nda$ preserves $\h$ it follows
easily that if $U\in \act(w)$ and $ V\in \act(w')$ then $ \nda_U V \in
\act(w+w'-1)$.  The connection $ \nda $ is torsion free so $ \NX U -
\nda_U \X-[\X,U]=0$ for any tangent vector field $ U$.  Now $\nda_U \X
=U$, so this simplifies to $ \NX U = [\X,U]+U$.  Thus if $ U\in \act$,
or equivalently $[\X, U] =-U$, then $\NX U =0$. The converse is clear
and it follows that sections of $ \act$ may be characterised as those
sections of $ T\aM$ which are covariantly parallel along the integral
curves of $ \X$ (which on $ \cq$ are exactly the fibres of $\pi$).
These two results imply that $ \nda $ determines a connection $ \nd$
on $\bT$. For $ U\in \ct$, let $ \tU $ be the corresponding section of
$ \act|_{\cq}$. Similarly a tangent vector field $ \xi$ on $M$ has a
lift to a field $ \tilde{\xi}\in \act(1)$, on $ \Cal Q$, which is
everywhere tangent to ${\Cal Q}$.  This is unique up to adding $f\X$,
where $ f\in \cce(0)$. We extend 
$ \tU$ and $\tilde{\xi} $ smoothly and homogeneously to
fields on $ \aM$.  Then we can form $ \nda_{\tilde{\xi}} \tU$; this is
clearly independent of the extensions. Since $ \NX \tU =0$, 
the section $ \nda_{\tilde{\xi}} \tU$ is also
independent of the choice of $ \tilde{\xi}$ as a lift of $ \xi$. Finally, 
$ \nda_{\tilde{\xi}} \tU$ is
a section of $ \act(0)$ and so determines a section $ \nd_\xi U $
of $ \bT$ which only depends on $ U$ and $ \xi$. It is easily verified
that this defines a covariant derivative on $ \bT$ which, by
construction, is compatible with the bundle metric $ h$.

The ambient metric is conformally invariant; no choice of metric from
the conformal class on $M$ is involved in solving the ambient metric
problem.  Thus the bundle, metric and connection $(\bT,h,\nd)$ are by
construction conformally invariant.  On the other hand the ambient
metric is not unique. Nevertheless it is straightforward to verify that $
\nd$ satisfies the required non-degeneracy condition and curvature
normalisation condition that lead to the following result.

\begin{proposition}
The bundle and connection pair $(\bT,\nd)$, induced by $\h$, is a
normal standard (tractor bundle, connection) pair.
\end{proposition}

This is proved in \cite{CapGoFG}. (In fact it is shown there that to
obtain the normal standard tractor bundle and connection it is
sufficient to replace property (iii) of the ambient metric with the
weaker condition that the tangential components of $\Ric(\h)$ vanish
along $ \cq$.) From a standard tractor bundle and connection it is
straightforward to construct a Cartan bundle $ \cg$ and connection
$\omega$ from which the tractor bundle and connection arise as
associated structures (via the defining representation of
SO$(p+1,q+1)$).  The notion of normality of a tractor connection is
equivalent to that on Cartan structure (see \cite{CapGotrans}). So
although the ambient metric is not unique the induced tractor bundle
structure $(\bT,h,\nd)$ is equivalent to a normal Cartan connection,
and so
is unique up to bundle isomorphisms preserving the filtration
structure of $\bT$, and preserving $h$ and $\nd$.

In particular this means that given a choice of metric $g$ from the
conformal class the structure $(\bT,h,\nd)$ can be expressed in terms
of $T^*M$, $g$ and the Levi-Civita connection for $g$ (which is also
denoted $\nd$) by explicit formulae which we give below.  In an
abstract index notation $TM$ is denoted $E^a$ and $E_a$ means $T^*M$;
we write $\ce^a$ and $\ce_a$ for the corresponding section
spaces. (We use the early part of the alphabet for abstract
indices. In view of this and context $\ce^a$ should not be confused
with the space of $k$-forms $\ce^k$).  
Similarly the section
spaces of the tractor bundle and its dual can also be denoted $\ct^A$
and $\ct_A$.  It is often convenient choose a metric $g$ from the
conformal class which determines \cite{BEGo,CapGoluminy} a canonical
splitting of the composition series \nn{trcomp}. Via this the
semi-direct sums $\lpl$ in that series get replaced by direct sums
$\oplus$, and we introduce $g$-dependent sections
$Z^A{}^b\in\ct^{Ab}[-1]$ and $Y^A\in\ct^A[-1]$ that describe this
decomposition of $\bT$ into the direct sum $ \bT^A= E[1]\oplus E_a[1]
\oplus E[-1]$.  A section $V\in\ct$ then corresponds to a triple 
$(\sigma,\mu,\rho)$ of
sections from the direct sum according to $V^A=
Y^A\sigma+Z^{Ab}\mu_b+X^A\rho $, and in these terms
the tractor metric is given by $h(V,V)=\bg^{ab}\mu_a\mu_b +2 \si
\rho$.
The sections $Y$ and $Z$ are defined in terms of the
Levi-Civita connection, and have ambient space equivalents which will be partially
described below. 
If $\hat{Y}^A$ and $ \hat{Z}^A{}_b$ are the
corresponding quantities in terms of the metric $ \hat{g}=e^{2\om}g$
then we have
\begin{equation*}
\textstyle
\begin{array}{rl}
\hat Z^{Ab}=Z^{Ab}+\Up^bX^A, &
\hat Y^A=Y^A-\Up_bZ^{Ab}-\frac12\Up_b\Up^bX^A,
\end{array}
\end{equation*}
where $\Up:=d\om$.
In terms of this splitting for $g$ the tractor connection is given by  
\begin{equation}\label{connids}
\begin{array}{rcl}
\nd_aX_A=Z_{Aa}\,, &
\nd_aZ_{Ab}=-\V_{ab}X_A-Y_A\bg_{ab}\,, & \nd_aY_A=\V_{ab}Z_A{}^b ,
\end{array}
\end{equation}
(see \cite{BEGo,GoPet}) where $\V_{ab}$ is a trace adjustment of a
constant multiple of $\Ric(g)$ known as the Schouten (or Rho)
tensor. (Note that in \nn{connids}, $\nd$ is the coupled
tractor--Levi-Civita connection.) 

The bundle of $k$-form tractors $\bT^k$ is the $k^{\underline{\rm
th}}$ exterior power of the bundle of standard tractors. This has a
composition series which, in terms of section spaces, is given by
\begin{equation} \label{formtractorcomp}
\cT^k=\Lambda^k\cT\cong\ce^{k-1}[k]\lpl\{\ce^k[k]\oplus
\ce^{k-2}[k-2]\}\lpl\ce^{k-1}[k-2].
\end{equation}
Given a choice of
metric $g$ from the conformal class there is a splitting 
of this composition series 
corresponding to the splitting of $\ct$ as mentioned above.   
Relative to this, a typical $k$-form tractor field $F$
corresponds to a 4-tuple $(\si,\m,\f,\r)$ of sections of the direct sum 
(obtained 
by replacing each $\lpl$ with $\oplus$ in \nn{formtractorcomp}) and we write
$$
F=\bbY^k\cdot \si +\bbZ^k \cdot \m + \bbW^k \cdot \f+ \bbX^k\cdot \r,
$$
where `${\cdot}$' is the usual pointwise form inner product in the
tensor arguments,
$$
\f\cdot\psi=\dfrac1{p!}\f^{a_1\cdots a_p}\psi_{a_1\cdots a_p}\ 
\mbox{ for $p$-forms},
$$
and for $k>1$, if $\uw$ is the wedge product in the tractor arguments,
\begin{equation}\label{parts}
\bbZ^k=Z\uw\bbZ^{k-1},\ \ \bbX^k=X\uw\bbZ^{k-1},\ \
\bbY^k=Y\uw\bbZ^{k-1},\ \
\bbW^k=Y\uw X\uw\bbZ^{k-2}.
\end{equation}
(By convention, $\bbZ^0=1$ and $\bbZ^{-1}=0$.)
Note that because $Z$ is vector valued,
$Z\uw Z$ does not vanish, though expressions like $X\uw X$ and $Y\uw Y$
do vanish.  The form tractor bundles $\bT^k$ are non-zero for $k=0,\ldots,
n+2$; $\bbZ^k$ vanishes for $k\ge n+1$; $\bbW$ vanishes for $k\le 1$;
and $\bbX^k$, $\bbY^k$ vanish for $k=0,n+2$.
$\bbX^k$ is an invariant section, while $\bbY^k$ depends on
a choice
of scale.  $\bbZ^k$ depends on a choice of scale unless $k=0$;
$\bbW$ depends on a choice of scale unless $k=n+2$.

An invariant metric on $\bT^k$ is
\begin{equation}\label{formtracmet}
\ip{(\nu,\m,\f,\r)}{(\tilde \nu,\tilde\m,\tilde\f,\tilde\r)}=
\nu\cdot\tilde\r+\r\cdot\tilde\nu+\m\cdot\tilde\m-\f\cdot\tilde\f.
\end{equation}
In fact, this is the restriction of the ambient $k$-form metric 
\begin{equation}\label{ambimet}
\Phi\bullet\Psi:= \frac{1}{k!}\Phi^{A_1\cdots A_k}\Psi_{A_1\cdots A_k},
\end{equation}
in the sense that after restricting to homogeneous ambient $k$-forms
along $\cQ$ and identifying these with $k$-form tractors we obtain a
$\bullet$ on the latter. Choosing a conformal scale we observe that 
$$
(\bbY^k\cdot \nu +\bbZ^k \cdot \m + \bbW^k \cdot \f+ \bbX^k\cdot \r)
\bullet
(\bbY^k\cdot \tilde\nu +\bbZ^k \cdot \tilde\m + \bbW^k \cdot \tilde\f
+ \bbX^k\cdot \tilde\r)
$$
reduces to the expression in \nn{formtracmet}.
This follows in turn from the formulae
$$
Z^{Ab}Z_{Ac}=\delta^b{}_c,\ \ X^AY_A=1,
$$
with all other quadratic contractions of $X,Y,Z$ vanishing, together
with formula \nn{wedgie} below 
for the $\uw$ with a tractor-one-form.

If $\alpha$ and $\beta$ are one-forms and $\f$ is a form
(or if these objects are form-densities), let
$$
E(\alpha\otimes\beta)\f:=\alpha\otimes\e(\beta)\f,\ \
I(\alpha\otimes\beta)\f:=\alpha\otimes\i(\beta)\f,
$$
and extend from simple tensors $\alpha\otimes\beta$
to arbitrary 2-tensors by linearity.  
The formulae for the covariant derivatives of $X,Y,Z$ at a scale imply
that
\begin{equation}\label{trconnform}
\begin{array}{rlrl}
\nd\bbX&=-E(\bg)\bbW+I(\bg)\bbZ,\qquad\qquad&
\nd\bbZ&=-E(\V)\bbX-E(\bg)\bbY, \\
\nd\bbW&=I(\V)\bbX-I(\bg)\bbY,\qquad\qquad&
\nd\bbY&=I(\V)\bbZ+E(\V)\bbW,
\end{array}
\end{equation}
where we have suppressed
the superscript $k$.
Under a change of scale $\hat g=e^{2\om}g$, the behaviour of
$X,Y,Z$ gives
\renewcommand{\arraystretch}{1.1}
$$
\begin{array}{rl}
\hat \bbX&=\bbX, \\
\hat \bbZ&=\bbZ+\e(\Up)\bbX, \\
\hat \bbW&=\bbW-\i(\Up)\bbX, \\
\hat\bbY&=\bbY-\i(\Up)\bbZ-\e(\Up)\bbW
+\frac12(\e(\Up)\i(\Up)-\i(\Up)\e(\Up))\bbX,
\end{array}
$$
\renewcommand{\arraystretch}{1.1}
where again $\Up=d\omega$.

\subsection{Exterior calculus on the ambient manifold} \label{ambext}

Let $\ad$ be the exterior derivative on the ambient manifold
$\aM$, and let $\da$ be its formal adjoint with respect to the
usual form metric, which is derived in turn from the ambient metric $\bh$.
If $u$ is a one-form on $\aM$, we have exterior multiplication
$\e(u)$ and its formal adjoint,
the interior multiplication $\i(u)$. Using the ambient metric 
\nn{ambimet} to raise and 
lower indices
the conventions are as follows. 
Exterior and interior multiplication by a 1-form $\w$ are given by
\begin{equation}\label{wedgie}
\begin{array}{rl}
(\e(\w)\f)_{A_0\cdots A_k}&=(k+1)\w_{[A_0}\f_{A_1\cdots A_k]}\,, \\
(\i(\w)\f)_{A_2\cdots A_k}&=\w^{A_1}\f_{A_1\cdots A_k}\,.
\end{array}
\end{equation}
We extend the notation for interior and exterior multiplication in
an obvious way to operators which increase the rank by one. For
example since the ambient connection is symmetric we have
$\ad\f=\e(\nda)\f$ and $\da \f = -\i(\nda) \f$. These notations and
conventions are also used for form tractors, and for forms and form-densities
on the underlying conformal manifold.

Building polynomially on $\ad,\da,\e(\X),\i(\X)$, we get several
more differential operators.  In particular, we get the
{\em form Laplacian} $\afl:=\da\ad+\ad\da$, and
$Q=\i(\X)\e(\X)+\e(\X)\i(\X)$.
We also obtain the Lie derivative with respect to $\X$, and its formal
adjoint as operators on forms:
\begin{equation}\label{LXdie}
\LX=\i(\X)\ad+\ad\i(\X),\ \ \LX^*=\da\e(\X)+\e(\X)\da.
\end{equation}

In general, given a bundle, we shall use the notation $\G(\cdot)$ for its
smooth section space, if this space has not been given another name.
The subspace of $\Gamma(\wedge^k T^*\aM)$ consisting of ambient
$k$-forms $F$ satisfying $\NX F =w F$ for a given $w\in {\Bbb R}$ will be
denoted $\act^k(w)$. We say such forms are (homogeneous) of {\em
weight} $ w$. From the definitions in the previous section, it is
straightforward to verify that
$$
\act^k(w) =(\Lambda^k \act)\otimes \cce(w)
$$ where the tensor (and exterior) products are over
$\cce(0)$. 
Given its weight, the
degree of $F$ is dependent on its order. In general from \nn{ndaX} we
have the identities
\begin{equation}\label{LXNXp}
\begin{array}{rl}
\LX&=\NX+{\ix}, \\
\LX-\LX^*&=2\NX+n+2, 
\end{array}
\end{equation}
where $ \ix$ is the operator that multiplies by $p-q$ the part of a
tensor with rank $ (q,p)$ (i.e.\ $q$ indices up and $p$ down).

Considering further commutators and anticommutators we note that the
eight operators in Tables 1 and 2
generate an isomorphic copy
$\fg$ of the linear Lie superalgebra $\fsl (2|1)$.
This decomposes into  the $-1$, $0$ and $1$ eigenspaces of 
(the bracket with) $Z:=(\LX+\LX^*)/2$:
$\fg=\fg_1^{-}\oplus \fg_0\oplus \fg_1^+$.  The
odd part of $\fg$ is $\fg_1= \fg_1^{-}\oplus \fg_1^+$; the subspace 
$\fg_1^-$ is
spanned by $\da$ and $\i(\X)$; the subspace $\fg_1^+$ is spanned by 
$\ad$ and $\e(\X)$; and via
anticommutators these generate $\fg_0$, that is $\{\fg_1,\fg_1\}=\fg_0$.
Denoting by $E_{ij}$ ($i,j=1,2,3$) the standard matrix units in ${\Bbb
C}^3\times {\Bbb C}^3$ (or ${\Bbb R}^3\times {\Bbb R}^3$) 
one family of Lie superalgebra 
isomorphisms from $\fg$ to the defining representation of 
$\fsl (2|1)$ is given by
$$
\begin{array}{llll}
 \i(\X) \mapsto -\frac{2}{c}E_{31}, & \da \mapsto \frac{2}{c}E_{32}, 
& \ad \mapsto c E_{13}, & \e(\X) \mapsto -c E_{23},\\
Q\mapsto 2 E_{21},& \afl \mapsto 2 E_{12}, 
& \LX \mapsto  -2 E_{11} -2 E_{33}, & \LX^* \mapsto -2 E_{22}-2 E_{33}, \\
\end{array}
$$ 
where $c$ is a non-zero complex number (or real number if we work over
$\Bbb R$).  Observe that the even part $\fg_0$ of $\fg$ is
isomorphic to $\fu(2)$. Since $Z$ is central in $\fg_0$, it is clear
that the decomposition of $\fg_1^{-}\oplus \fg_0\oplus \fg_1^+$ is
$\fg_0$-equivariant and we note from the isomorphism that $\fg_1^-$
and $\fg_1^+$ are dual $\fg_0$-modules.  In the $0$-form-density case,
the $\fsu(2)$ subalgebra played a role in \cite{GJMS}.  
\begin{table}[ht]
\begin{center}
\begin{tabular}{c|c|c|c|c|}
\{\cdot,\cdot\} & $\ad$ &    $\da$ &    $\e(\X)$ &    $\i(\X)$ \\
\hline
$\ad$   &   $0$       &    $\afl$   &    $0$         &   $\LX$ \\
\hline
$\da$ & $\afl$ & $0$ & $\LX^*$ & $0$ \\
\hline
$\e(\X)$ & $0$ & $\LX^*$ & $0$ & $Q$  \\
\hline
$\i(\X)$ & $\LX$ & $0$ & $Q$ & $0$  \\
\hline
\end{tabular}\end{center}
\caption{Anticommutators $\{\mathfrak{g}_1,\mathfrak{g}_1\}$}
\label{anticomms}
\end{table}
\begin{table}[ht]
\begin{center}
\begin{tabular}{c|c|c|c|c|c|c|c|c|c|}
$[\cdot,\cdot]$ & $\ad$ & $\da$ & $\e(\X)$ & $\i(\X)$ & & $\afl$ & $\LX$ & 
$\LX^*$ & $Q$ \\
\hline
$\afl$ & $0$ & $0$ & $-2\ad$ & $2\da$ & & $0$ & $2\afl$ & $-2\afl$ & $-2\KX$ \\
\hline
$\LX$ & $0$ & $-2\da$ & $2\e(\X)$ & $0$ & & $-2\afl$ & $0$ & $0$ & $2Q$ \\ 
\hline
$\LX^*$ & $2\ad$ & $0$ & $0$ & $-2\i(\X)$ & & $2\afl$ & $0$ & $0$ & $-2Q$ \\ 
\hline
$Q$ & $-2\e(\X)$ & $2\i(\X)$ & $0$ & $0$ & & $2\KX$ & $-2Q$ & $2Q$ & $0$ \\
\hline
\end{tabular}\end{center}
\caption{Commutators $[{\mathfrak{g}}_0,{\mathfrak{g}}_1]$ and
$[{\mathfrak{g}}_0,{\mathfrak{g}}_0]$, where $\KX:=\LX-\LX^*$}
\label{comms}
\end{table}

The relations in these tables all follow from 
$Q=\h(\X,\X)$, $dQ =2\X$, (\ref{ndaX}), and the usual
identities of exterior calculus on pseudo-Riemannian manifolds 
(e.g.\ (\ref{LXdie})).  In
particular, they hold in all dimensions and to all orders. 

Of particular interest are differential operators $P$ on ambient form
bundles, or subquotients thereof, which act {\em tangentially} along
$\cq$, in the sense that $PQ=QP'$ for some operator $P'$ (or
equivalently $[P,Q]=QP''$ for some $P''$). Note that compositions of
tangential operators are tangential. If tangential operators are
suitably homogeneous then they descend to operators on on $ M$. Such
$P$ are {\em fragile} if they are tangential only when acting on
sections $F$ of some weight $w$; a
fragile operator descends to an operator which is invariant for a
particular weight.  The {\em robust} $P$ are those which are
tangential when acting on arbitrary smooth sections; these descend to
operators which are invariant for any weight.

An example of a fragile tangential operator is given by:

\begin{proposition}\label{aflmtang}
$\afl^m:\act^k(m-n/2)\to\act^k(-m-n/2)$ is tangential.
\end{proposition}

\noindent{\bf Proof:} 
We need to calculate
$\afl^m(Qf)$ for 
$Qf$ of homogeneity $m-n/2$ (i.e.\ for $f$ of homogeneity $m-2-n/2$).
Without any homogeneity assumption, we have
\begin{equation}\label{mth}
[\afl^m,Q]=\sum_{p=0}^{m-1}\afl^{m-1-p}[\afl,Q]\afl^p,
\end{equation}
and $[\afl,Q]=2(\LX^*-\LX)$, from Table \ref{comms}.  Thus letting 
\nn{mth} act on $\act^k(w)$,
the $p^{\underline{{\rm th}}}$ term on the right acts as
$-2[2(w-2p)+n+2]\afl^{m-1}$, 
so that $[\afl^m,Q]$ acts
as $-2m(2w-2m+n+4)\afl^{m-1}$.  This vanishes identically if and only if 
$w=m-2-n/2$, so that $\afl^m$ is tangential on $\act^k(m-n/2)$, as
desired. \quad $\Box$
 
\noindent{\bf Remark:} We should point out that results along these
lines are not peculiar to forms or the form Laplacian. Recall that when
acting on $\act^k(w)$, $2(\LX^*-\LX)= -2(2\NX+n+2)$. On the other hand
from \nn{ndaX} one calculates that for the ambient Bochner Laplacian
$\al:=\nda^*\nda=-\nda^A\nda_A$ we have
\begin{equation}\label{alQ}
[\al,Q]=-2(2\NX+n+2) ,
\end{equation}
as an operator on {\em any ambient tensor}.  Thus by
essentially the same argument as above we conclude that $\al^m$ is
tangential on arbitrary tensors of weight $w$. $\qquad\endrk$

Simple examples of robust tangential operators are given in the 
following proposition:

\begin{proposition}\label{efdd}
The operators
$$
\e(\afD): = \ad(n+2\Euler-2) +\e(\X)\afl,
\quad \i(\afD):=-\da(n+2\Euler-2) + \i(\X)\afl 
$$ 
act tangentially along $ \cq$.
\end{proposition}

\noindent{\bf Proof:} Using Table \ref{comms} one calculates 
$[Q,\e(\afD)]=-4Q\ad$ and $[Q,\i(\afD)]=4Q\da$.
$\quad\Box$\\

Note that in view of the relations $[\NX,\e(\afD)]=-\e(\afD)$ and
$[\NX,\i(\afD)]=-\i(\afD)$, these operators lower weight by 1 and by 
construction $ \e(\afD)$ raises form order by 1 while $ \i(\afD)$ lowers it by 1. In summary:
$$
\e(\afD):\act^k(w)\to\act^{k+1}(w-1),\ \
\i(\afD):\act^k(w)\to\act^{k-1}(w-1).
$$
These satisfy  identities as follows.

\begin{proposition}\label{DD0}
On $\aM$,
$$
\i(\afD)\i(\afD)=0,\ \ \e(\afD)\e(\afD) =0,\ \
\i(\afD)\e(\afD)+\e(\afD)\i(\afD)=Q\afl^2.
$$
\end{proposition}

\noindent{\bf Proof:}
These formulae follow from \nn{LXNXp} and 
Tables \ref{anticomms} and \ref{comms}.
$\quad\Box$

\subsection{Form tractors and invariant operators} 

Recall that we write $\ct^k$ for the space of {\em $k$-form tractor}
 fields. That is $\ct^k$ is the space of sections of $\bT^k$.
 $\ct^k[w]$ denotes the space of {\em weighted}
 $k$-form tractors of weight $w$; this is the space of sections of
 $\bT^k\otimes E[w]$. (Naturally these are non-trivial for $0\leq k
 \leq n+2$. For $k$ outside this range we take $\bT^k$ to be the zero
 bundle.) From the relationship between $\ce[w]$ and $\cce(w)$ and the
 definition of $\cT$ in Section \ref{tractor}, it follows that the
 sections in $\ct^k[w]$ are equivalent to ambient manifold form fields, 
along $\cq$, 
which lie in $ \act^k(w)|_\cq$. 
Clearly
 the operators of Propositions \ref{aflmtang} and 
\ref{efdd} determine  operators between twisted form tractor bundles. By
 construction these tractor operators are conformally invariant but 
{\em ab initio}
 might depend on the choices in the ambient metric. Hence there is a need for
 specific results on when such operators are natural conformal
 objects.

A similar comment applies to natural tensors on the ambient
manifold.  The ambient curvature $\aR$ is in $(\otimes^4\act)(-2)$ and so
corresponds to a section of $(\otimes^4\ct)[-2]$. In dimensions other
than 4 this tractor field depends only on the conformal structure
\cite{CapGoFG,GoPet} and is $1/(n-4)$ times the tractor $W$ of
\cite{gosrni}.
In terms of a metric from the conformal class and the notation from Section 
\ref{tractor},
$W$ is given by
\begin{equation}\label{Wform}
\begin{array}{rl}
W_{ABCE}&=(n-4)\left\{\tfrac14\bbZ^{ab}_{AB}\bbZ^{ce}_{CE}C_{abce}
-\bbZ_{AB}^{ab}\bbX_{CE}^e\nd_{[a}\Rho_{b]e}\right. \\
\qquad&\left.-\bbX_{AB}^b\bbZ^{ce}_{CE}\nd_{[c}\Rho_{e]b}\right\}
+\bbX_{AB}^b\bbX_{CE}^e(2\nd^q\nd_{[q}\Rho_{b]e}+\Rho^{qp}C_{peqb}),
\end{array}
\end{equation}
where $C$ is the Weyl tensor.

Closely related to $\e(\afD)$ and $\i(\afD)$ is the operator $\D:=\nda
(n+2\NX-2) +\X \al$ of \cite{CapGoFG,GoPet} (and see e.g.\ \cite{BEGr}
where this was used in the setting of
the standard flat model).  Using $2\h(\X,\cdot)=d Q$,  \nn{alQ} and \nn{LXNXp} 
it is easily verified that $\D$ is a robust tangential operator on
ambient tensor fields of any rank or symmetry. For the class of ambient metrics
that we consider (in particular we need \nn{XRrel}) the operator $\D$,
restricted to homogeneous tensors along $ \cq$, is equivalent 
\cite{CapGoFG,GoPet} to the well-known tractor-D operator $D$ of
\cite{T,BEGo}.  $D$ is natural and so depends only on the conformal
structure. Explicitly, for a metric from the conformal class, $D$ is
given by
\begin{equation}\label{Dform}
D_A V:=(n+2w-2)w Y_A V+ (n+2w-2)Z_A{}^{a}\nabla_a V +X_A\Box V,
\end{equation} 
 where $ V$ is a section of any twisted tractor bundle of weight $w$,
 and writing $\J$ for the trace (by $\bg^{-1}$) of the Schouten tensor,
 we have 
$$
\Box V :=-(\nd_p\nd^p V+w\J V). 
$$
In these formulae $\nd$ means the coupled tractor--Levi-Civita connection.

Let us write $\hash$ ({\em hash}) for the natural tensorial action of
sections $A$ of $\End(T\aM)$ on tensors. If $A$ is skew for $\h$ then
this commutes with the raising and lowering of indices. 
As a section of the tensor square
of the $\h$-skew bundle endomorphisms of $T\tilde M$, the ambient
curvature has a double hash action on tensors, 
 and, on forms, this
is exactly the difference between the ambient Bochner and form
 Laplacians. That is we have
\begin{equation}\label{lapdiff}
\afl-\al=-\aR\hash\hash .
\end{equation}
It follows immediately that, as operators on $k$-forms,
\begin{equation*}
\e(\afD) =\e(\D) -\e(\X)\aR\hash\hash,\ \
\i(\afD) =\i(\D) -\i(\X)\aR\hash\hash .
\end{equation*}
Summarising with some additional
results we have the following.

\begin{proposition}\label{tracefdd}
The operators of Proposition \ref{efdd}
descend to natural conformally invariant differential operators 
$$
\e(\fD):\ct^k[w]\to\ct^{k+1}[w-1],\ \
\i(\fD):\ct^{k}[w]\to\ct^{k-1}[w-1] ,
$$ 
and $\i(\fD): \ct^{k+1}[1-n-w]\to\ct^{k}[-n-w]$ is the formal adjoint of
$\e(\fD): \ct^k[w]\to\ct^{k+1}[w-1]$. These satisfy
$$\textstyle
\e(\fD) =\e(D)- \e(X) \Omega \hash \hash ,\ \
\i(\fD) =\i(D)- \i(X) \Omega \hash\hash ,
$$ 
where $\Omega \hash \hash$ is a curvature action {\rm(}and so
has order 0 as a differential operator{\rm)}.
\end{proposition} 

\noindent{\bf Proof:} It is immediate from Proposition \ref{efdd}
that the operators there descend to conformally invariant
operators. 

{}From the discussion above, in dimensions other than 4,  
$(n-4)\aR|_\cq$ is equivalent to the tractor
field $W$, where $W$ is given explicitly in
\nn{Wform}. Viewing $W$ as a section of $\wedge^2 \bT \otimes \End(\bT)$, 
 we note that from \nn{Wform} 
$\e(X)W=(n-4)\e(X)\Omega$,
where (we also view $\Om$ as a section of $\wedge^2 \bT \otimes \End(\bT)$ and)
$\Om$ is given by
$$
\tfrac14\bbZ_{AB}^{ab}\bbZ_{CE}^{ce}C_{abce}-\bbZ_{AB}^{ab}\bbX_{CE}^e
\nd_{[a}\Rho_{b]e}
-\bbX_{AB}^b\bbZ_{CE}^{ce}\nd_{[c}\Rho_{e]b}.
$$
It is easily verified that $\Om$ has Weyl tensor type symmetries.
In dimension 4, $\aR$ is not equivalent to a natural tractor but, as
explained in Section 3.2 of \cite{CapGoFG}, $\e(X)\aR|_\cq$ is
determined by the conformal structure and it follows easily from the
discussion there that it is equivalent to the tractor field
$\e(X)\Omega$.  Next, from the above, $\D$ along $\cq$ is equivalent to the
operator $D$ of \nn{Dform} on tractors. Thus in all dimensions we have
$$\textstyle
\e(\fD) =\e(D)-\e(X) \Om \hash \hash ,\ \ \ 
\i(\fD) =\i(D)-\i(X) \Om \hash\hash ,
$$ 
where the interpretation of the action $\Om\hash\hash$ on form
tractors is obvious from the corresponding action $\aR\hash\hash$ on
ambient form fields. 

{}From the explicit formulae for $\Om$ and $D$ it is immediately clear
that the operators $\e(\fD)$ and $\i(\fD)$ are natural and differential.
That $\i(D): \ct^{k+1}[1-n-w]\to\ct^{k}[-n-w]$ is the formal adjoint
of $\e(D): \ct^k[w]\to\ct^{k+1}[w-1]$ is a special case of a more
general result in Sec.\ 7 of \cite{BrGoPacific}. 
It follows easily from this, the Weyl-tensor-type
symmetries of $\Om$, and the fact that $\Om$ is clearly annihilated by
contraction with $X$, that $\e(\fD)$ and $\i(\fD)$ are mutual formal
adjoints.  $\quad\Box$

\noindent{\bf Remark:} 
{}From Proposition \ref{DD0} we immediately have the identities 
$$
\i(\fD)\i(\fD)=0,\ \ \e(\fD)\e(\fD) =0,\ \
\i(\fD)\e(\fD)+\e(\fD)\i(\fD)=0.
$$ 
Observe that if $k=0$ the $\aR \hash\hash$ action is trivial by
definition while if $k=1$ it amounts to an action of $\Ric(\h)$ and so
vanishes along $\cq$.  In either case we have $\e(\afD) =\e(\D)$, and
$\i(\afD)=\i(\D)$ along $\cq$. Thus acting on form tractors of rank
$k\leq 1$ we have $\e(\fD)=\e(D)$ and $\i(\fD)=\i(D)$.
$\qquad\endrk$

Before we continue with the main theme let us digress briefly, to give
a direct formula for the scale dependent tractor $Y$ in terms of the
invariant natural operator $\e(\fD)$ and the choice of scale.  Let
$\si\in \ce[1]$ be a choice of conformal scale. Then by the definition
of the Levi-Civita connection on densities we have $\nd \si=0$, and, from
\nn{Dform} we have $\si^{-1}D\si =\si^{-1}\e(\fD)\si= n Y -X {\sf J}$.
Let us write $I_\si$ for $\frac{1}{n}\si^{-1}\e(\fD)\si $. Thus we have
$I_\si\bullet I_\si=-2 {\sf J} /n$, $X\bullet I_\si=1$ and so
\begin{equation}\label{Yform}
Y=I_\si-\tfrac{1}{2}(I_\si\bullet I_\si) X.
\end{equation}
This formula plays an important role in later calculations.

Returning to our programme of constructing natural invariant operators,
we note an extension of our observation above, to the effect 
that when $n\neq 4$, the ambient
curvature is a natural conformal invariant.

\begin{lemma}\label{GoPet4.4}
The ambient tensors $\nda^s\al^t \aR|_\cq$ with
  $s,t\in \{0,1,2,\ldots\}$ are equivalent to 
conformally invariant tractor fields.
In odd dimensions these tractor fields are natural.  
 In even dimensions $n\neq 4$, the
same is true with the restriction $s+t\leq n/2-3$. 
\end{lemma}

\noindent{\bf Proof:} The first claim is clear since the tensors are homogeneous
and the ambient metric is conformally invariant.
Lemma 4.4 of \cite{GoPet} establishes that in 
odd dimensions, or when $s+t\leq n/2-3$,
$\nda^s\al^t \aR$ can be expressed as a partial contraction polynomial
in $\D$, $ \aR$, $ \X$, $\h$, and $ \h^{-1}$.  It follows that the
corresponding conformally invariant tractor sections are given by same
formal expression with the respective replacements $D$,
$W/(n-4)$, $ X$, $h$, and $ h^{-1}$. The claims concerning naturality are now
immediate as each of these is natural.  
$\quad\Box$

{\bf Remarks:} Some related results are in Theorem 3.4 of
\cite{CapGoFG}.  It should also be pointed out that $\al^{n/2-2}\aR$
corresponds to a natural conformal invariant related to the
ambient obstruction of \cite{FGast} -- see \cite{GoPetprogress}.  In
even dimensions the remaining ambient tensors $\nda^s\al^t \aR|_\cq$
with $s+t>n/2-3$  correspond to tractor fields which depend
on the choices involved in extending the ambient metric to, 
and beyond, order $n/2$. $\qquad\endrk$

\begin{proposition}\label{powerslap}
The operators of Proposition \ref{aflmtang}
descend to conformally invariant differential operators
$$
\fl_m:\ct^k[m-n/2]\to\ct^k[-m-n/2] \quad m=0,1,\cdots,
$$
where by convention $\afl^0$ and $\fl_0$ are identity operators.
In odd dimensions these are natural operators.
In even dimensions the same is true with the restrictions that 
either 
$m\leq n/2-2$; or
$m\leq n/2-1$ and $k=1$; or
$m\leq n/2 $  and $ k=0 $.
In the conformally flat case the operators are natural with no
restrictions on $m\in \{0,1,2,\ldots\}$.
In every case $\fl_\ell$ has the same principal part as
$\Delta^\ell$.
\end{proposition}
Note that the $k=0$ cases of the above theorem are by construction
exactly the GJMS operators of \cite{GJMS}.

\noindent{\bf Proof: } From the conformal invariance of the ambient
construction, and the relationship between $\act^k(w)$ and $\ct^k[w]$
it is clear that the operators of Proposition \ref{aflmtang} determine
conformally invariant operators $ \fl_m:\ct^k[m-n/2]\to\ct^k[-m-n/2] $
for all $m\in\{0,1,2,\ldots\}$.

It remains to establish that these are differential, natural and with
leading term as claimed.  For odd dimensions and even dimensions, up
to order $n$, the $k=0$ cases are dealt with in \cite{GJMS} and
\cite{GoPet}.  The arguments of the latter adapt easily to the more
general setting here.  The first observation in \cite{GoPet} is that
for a function on $\aM$ homogeneous of weight $m-n/2$, $\al^m f$ (or
rather $(-1)^{m-1}\X^{m-1}\al^{m} f$) is the leading term of $\al
\D^{m-1} f$. The difference between these terms involves 
$\nda$-derivatives of $f$ and $\aR$ and the main part of the argument is that
these can be re-expressed as $\D$-derivatives of $f$ and $\aR$ . (See
in particular Section 4 and the proof of Theorem 2.5).  It is easily
verified that if we instead just re-express the $\nda$-derivatives of
$f$ in this way but do not re-express the $\nda^\ell \aR$, then in all
dimensions the argument works for all $m\in {\Bbb Z}^+$. With this
variant the proof goes through, and is essentially unchanged, if we
replace $f$ with an ambient tensor field homogeneous of weight
$m-n/2$. (The only difference is that more curvature terms turn up.) A
similar argument can be applied to powers of $\afl$ by using
\nn{lapdiff} to first re-express $\afl^m$ in terms of $\al$, $\nda$
and $\aR$.  Thus for $U\in \act^k(m-n/2)$, it follows that $\afl^m U$
has an expression which is polynomial in $\h$, $\h^{-1}$,
$\nda$-derivatives of $\aR$, and is linear in $\D$-derivatives of $U$.  It
follows immediately that the operators $\fl^m$ here are differential,
since $\D$ descends to the natural differential operator $D$. The
leading term is $\Delta^m$, since from \nn{Dform} it follows easily
this is the leading term of $\Delta D^{m-1}$. Counting the powers of
$\nda$ and $\al$ that can act on $\aR$ in the expansion discussed, the
statements on naturality in odd dimensions and in even dimensions for
$m\leq n/2-2$ are now immediate from Lemma \ref{GoPet4.4}.

The improved result when $k=1$ in even dimensions is a consequence of
the observation already made that the $\aR\hash\hash$ action reduces
to a $\Ric(\h)$ action on 1-forms while by property (iii) of the
ambient metric, \nn{ndaX} and \nn{alQ} it follows that $\nda^s\al^t
\Ric(\h)|_\cq=0$ if $s+t\leq n/2-2$.  On conformally flat manifolds we
take the flat ambient metric and so the  claim for this setting is clear.
\quad $\Box $ 
 
Next we will construct natural conformally invariant operators 
$$
\begin{array}{ll}
\td : \cg^k[w]\to \cg^{k+1}[w]\quad &\mbox{and}\\
\dt : \cg_{k+1}[w]\to \cg_k[w] \quad &\mbox{for } k=0,1,\cdots ,n+1,
\end{array}
$$ 
which in an appropriate sense generalise $d$ and $\d $.
The operators act between section spaces that we define as
follows. First denote by $\bV^k[w]$ the subbundle of $\bT^k[w]$
consisting (pointwise) of $k$-form tractors annihilated by $\e(X)$.
Denote by $\bF^k[w]$ the subbundle of $\bT^k[w]$ consisting of
$k$-form tractors annihilated by $\i(X)$.  (Equivalently $\bV^k[w]$ is
the subbundle of form tractors of the form $\e(X)S$ for some $(k-1)$
form tractor $S$ and similarly $\bF^k[w]$ is the subbundle of tractors
of the form $\i(X)F$ for some $(k+1)$ form tractor $F$. Also note that
the notation for the bundles defined here is consistent with Section
\ref{tractor} in the sense that $\bV^1=\bV$ and $\bF^1=\bF$.) Then
for $k=0,1,\cdots ,n+1 $ we have the definitions:
\begin{equation}\label{cgdefs}
\begin{array}{lcl}
 \bG^k[w]:=\bT^k[w-k]/\bV^k[w-k] &\quad &  E^{k-1}[w]\lpl E^k[w]  \\
\bG_{k}[w]:=\bF^k[w+k-n] &\quad &  E_k[w]\lpl E_{k-1}[w].
\end{array}
\end{equation}
The second column here gives the composition series of the bundle
defined (which follow at once from \nn{formtractorcomp} and the
definitions here). We use the following notations for the section spaces:
$$
\begin{array}{ll}
\cg^k[w] :=\Gamma(\bG^k[w]), \qquad & \cg_k[w] :=\Gamma(\bG_k[w]), \\
\cf^k[w] :=\Gamma(\bF^k[w]), \quad & \cV^k[w] := \Gamma(\bV^k[w]).
\end{array}
$$ 
We note that via the form tractor metric $\bG_k[n-w]$ is identified
 with the bundle dual to $\bG^k[w]$, and so there is an integral pairing
 between $\cg^k[w]$ and $\cg_k[-w]$. The latter is part of the reason
 for using $\bG_k[w]$ as an alternative notation for
 $\bF^k[w+k-n]$. Along the lines of conventions for other spaces we
 write $\bG^k:=\bG^k[0]$ and $\bG_k:=\bG_k[0]$. 
Note that the wedge product $\uw$ of \nn{parts} induces a wedge
product carrying $\bG^k[w]\times\bG^{k'}[w']$ to $\bG^{k+k'}[w+w']$,
and similarly for the weighted $\bG_k$ bundles.

Since $X$ and $Y$ are null tractor fields, we have the well-defined operations
$$
\begin{array}{lcl}
\i(X):\cg^{k+1}[w]\to \cg^k[w], &\quad & \e(Y):\cg^k[w]\to \cg^{k+1}[w]  \\
\i(Y):\cg_{k+1}[w]\to \cg_k[w], &\quad & \e(X):\cg_k[w]\to \cg_{k+1}[w],
\end{array}
$$ 
arising from interior and exterior multiplication on
$\ct^k[w]$. Thus we have that $\{\i(X),\e(Y)\}$ acts as the identity
on $\cg^k[w]$ and $\{\i(Y),\e(X)\}$ acts as the identity on
$\cg_k[w]$.  Also $\e(X):\cg^k[w]\to \cV^{k+1}[w-k+1]$ is clearly well
defined (and in fact injective since $\{\i(Y),\e(X)\}$ is the identity
on $\ct^k[w-k]$) while $\i(X)$ acts as 0 on $\cg_k[w]$.  

\noindent{\bf Remark:} \label{cosubtractors} 
There are other ``subtractor'' bundles similar
 to those defined in \nn{cgdefs}. These are used for the constructions
 outlined in Section \ref{codetoursect}.  We outline briefly their
 relationship to $\bG^k[w]$ and $\bG_k[w]$. 

 Suppose $M$ is orientable, with an orientation given by the
 conformal volume form $\vol$. Recall that $\vol$ is a section of $\ce^n[n]$
with the property that, for each 
choice of conformal scale $\si$, $\si^{-n}\vol$ is
 the unique volume form compatible with the orientation and the metric
 $g=\si^{-2}\bg$.  Then
$$
\bW^{n+2}\cdot \vol
$$ 
is a conformally invariant canonical section of $\ct^{n+2}$. In an
obvious way one can use this and the tractor metric to define a Hodge
star operator, which we denote $\tstar$ , for form tractor bundles. 

The usual Hodge relations apply to $\tstar$ and in particular we have 
$$ 
\begin{array}{l}
\tstar \e(X)= (-1)^k \i(X)\tstar,\qquad \e(X)\star
=(-1)^{k-1} \tstar \i(X), \\
\mbox{and }\ \tstar \tstar = (-1)^{k(n+2-k)+q+1}.
\end{array}
$$
It follows that on oriented manifolds we have an isomorphisms
$$
\tstar: \bG^k[w]\to \bG_{n-k}^{\star}[w], 
\quad \tstar: \bG_k[w]\to \bG^{n-k}_{\star}[w] ,
$$ 
where $\bG_{n-k}^{\star}[w]$ is the quotient of $\bT^{n+2-k}[w-k]$
with composition series $E_{n-k+1}[w]\lpl E_{n-k}[w]$ and
$\bG^{n-k}_{\star}[w]$ is subbundle of $\bT^{n+2-k}[w+k-n]$ with
composition series $E^{n-k}[w] \lpl E^{n-k+1}[w]$.  Put another way,
$\bG_{n-k}^{\star}[w]$ is the pointwise quotient of $\bT^{n+2-k}[w-k]$
by form tractors of the form $\i(X)F$, and $\bG^{n-k}_{\star}[w]$ is
the subbundle of $\bT^{n+2-k}[w+k-n]$ consisting of forms annihilated
by $\e(X)$.  We write $\cg_{n-k}^\star[w]$ and $\cg^{n-k}_\star[w]$
for the section spaces of, respectively, $\bG_{n-k}^{\star}[w]$ and
$\bG^{n-k}_{\star}[w]$. $\qquad\endrk$

The operators $\td$ and $\dt$, to be constructed,
are simply restrictions of the exterior derivative and its formal
adjoint on $\cq$.  For our current purposes, it is useful to see how
they arise from tangential operators in the ambient picture.
First observe that it is clear that, when acting on
homogeneous tensors, $\ad$ preserves homogeneous degree while
increasing rank by 1, and so lowers the homogeneity weight by 1.
Since the ambient metric is homogeneous of weight 0 it follows that
$\da$ also lowers weight by 1. In summary 
$$
\ad: \act^k(w)\to \act^{k+1}(w-1), \qquad \da: \act^{k+1}(w)\to \act^k(w-1).
$$ 
Next observe that for any form field $F$, $\ad \e(\X)F=-\e(\X) \ad
 F$ so $\ad$ preserves the space of forms of the form $\e(\X)F$. Once
 again using the Tables \ref{anticomms} and \ref{comms}, 
we also note that $\ad Q F =2 \e(\X) F+ Q\ad
 F$.  Let us use the informal notation $\Lambda^kT^* \aM/\e(\X) $ for
 the quotient bundle which is the pointwise quotient of $\Lambda^kT^*
 \aM$ by forms of the form $\e(\X)F$. The space of smooth sections of
 this quotient bundle is identified with the space
 $\Gamma(\Lambda^kT^* \aM)$ modulo the subspace of smooth sections of
 the form $\e(\X)F$. Via this and our observations just above 
it is clear that $\ad$ determines a robust
tangential operator
(also to be denoted $\ad$) on $\Lambda^kT^* \aM/\e(\X) $.
 For each $w\in {\Bbb R}$ we
have an inclusion $\act^k(w) \hookrightarrow \Gamma(\Lambda^kT^*\aM)$ and so
$\act^k(w)$ has an image in this quotient of section spaces. Let us
use $\acg^k(w+k)$ to denote this. Once again using the notation $\ad$
for the restriction to this subspace, it follows that
$$
\ad: \acg^k(w)\to \acg^{k+1}(w)  
$$ 
is a tangential operator. 

On $\aM$ there is a natural integral pairing between the sections of\newline
$\Lambda^kT^* \aM/\e(\X) $ and the space of sections of
$\Lambda^kT^*\aM$ which are annihilated by $\i(\X)$. Since $\da$ is
the formal adjoint of the exterior derivative $\ad$, it preserves this
subspace, and its restriction to this subspace (denoted $\da$) may be
viewed as a formal adjoint of $\ad$ on the quotient $\Lambda^{k-1}T^*
\aM/\e(\X) $. Since the latter commutes with $Q$ it follows at
once that, as an operator on the space of sections of
$\Lambda^kT^*\aM$ which are annihilated by $\i(\X)$, $\da$ also
commutes with $Q$. We write $\acg_k(w-k+n)$ (or alternatively
$\acf^k(w)$) for the intersection of this last space with $\act^k(w)$.

{}From the definition of the tractor bundle $\bT$ in terms of the
ambient manifold in Section \ref{tractor} and the relationship between
$X\in \ct[1]$ and $\X\in \act(1)$ it follows that sections of
$\cg^k[w]$ are equivalent to sections from the space
$\acg^k(w)|_\cq$. Similarly sections of $\cf^k[w]$ are equivalent to
sections of $\acf^k(w)$.  We now have the following.

\begin{theorem} \label{gend's}
The operators 
$$
\ad: \acg^k(w)\to \acg^{k+1}(w) \ \quad \mbox{and}\quad \ 
\da: \acg_k(w)\to \acg_{k-1}(w)  
$$ 
are tangential and satisfy $\ad^2=0=\da^2$. These
operators determine
first order conformally invariant differential operators
$$
\td: \cg^k[w]\to \cg^{k+1}[w] \ \quad \mbox{and}\quad 
\ \dt: \cg_k[w]\to \cg_{k-1}[w]
$$ 
on $M$ which satisfy $\td^2=0=\dt^2$. 
The operator $\td$ satisfies the anti-derivation rule 
$\td(\e(U)V)=\e(\td U) V+(-1)^k \e(U) \td V$ 
for $U$ in $\cg^k[w]$ and $V$ in any $\cg^{k'}[w']$. 
\end{theorem} 

\noindent{\bf Proof:} We have already shown the operators are
tangential. For $V\in \cg^k[w]$, let $\tV$ be any homogeneous smooth
extension to a section of $\acg^k(w)$ of the equivalent section from
the space $\acg^k(w)|_\cq$. Then $\ad \tV|_\cq$ is a section of
$\acg^{k+1}(w)|_\cq$ dependent only on $\tV|_\cq$ and we write $\td V$
for the equivalent section of $\cg^{k+1}[w]$. This defines the
operator $\td$. By taking coordinates on the ambient manifold $\aM$,
it is easily verified that $\td$ is differential and first order. (See
formulae (\ref{tdform},\ref{dtform}) below.) 
The result $\ad^2=0$ (resp.\ ${\td}^2=0$) on
$\acg^k(w)$ (resp.\ $\cg^k[w]$) follows from the same result for the
exterior derivative on $\Lambda^k \aM$, since if $\tV' $ is a
homogeneous section of $\Lambda^k \aM$ representing $\tV\in \acg^k(w)$
(resp.\ $V\in \cg^k(w)$), then $\ad \tV'$ is a section of
$\Lambda^{k+1} \aM$ representing $\ad \tV\in \acg^{k+1}(w)$ (resp.\
$\td V \in \cg^{k+1}[w]$).

The corresponding results for $\da$ and $\dt$ follow by an analogous
argument.  The anti-derivation rule for $\td$ is immediate from its
definition.\quad $\Box$ 

It is useful to understand the geometric origins of the results above.
Recall that we write $i:\cq \to \aM$ for the embedding of $\cq$ in
$\aM$. The identification of $\cq$ with $i(\cq)$ induces an
identification of $T^*\cq=i^*(T^*\aM)$ with $(i_*T\cq)^*$. This in
turn is canonically identified with the quotient of $T^*\aM|_\cq$ by
the conormal bundle. Since (along $\cq$) $2\X=\ad Q$ is a section of
the latter it follows immediately that this quotient is, in terms of
the notation above, precisely $\left( T^*\aM/\e(\X)\right) |_\cq $. 
Taking exterior
powers we have $\Lambda^k T^*\cq$ identified with $\left(\Lambda^k T^*\aM
/\e(\X)\right)|_\cq$.  Now for $\phi \in \Gamma(\Lambda^k T^* \cq)$ let us
say that $\tilde \phi\in \Gamma(\Lambda^k\aM)$ is an {\em extension}
of $\phi$ if $i^* \tilde \phi=\phi$. Writing $\ad$ also for the
exterior derivative on $\cq$ we have $\ad i^* \tilde \phi =i^* \ad
\tilde \phi $. So $\ad \tilde \phi$ is an extension of $\ad \phi$ and
hence $\ad$ is tangential as an operator on $\Lambda^k
T^*M /\e(\X)$. Thus $\ad:\acg^k(w)|_\cq \to \acg^{k+1}(w)|_\cq$ 
(is again shown to be 
well defined and) is really just a
restriction of the exterior derivative on $\cq$ to homogeneous
sections.

Similarly $\da:\acg_k(w)\to\acg_{k-1}(w)$ may be viewed as a restriction
of the formal adjoint of the exterior derivative on $\cq$.  Since
$\cq$ has no nondegenerate metric, or even conformal structure, we should
view the latter as an operator on weighted exterior powers of the tangent
bundle $T\cq$.  From this viewpoint the properties of $\da$ on 
$\acg_k(w)$ follow easily.  Rather than introduce new notation to explain
this explicitly, essentially the same argument may be phrased in terms
of the ambient metric as follows.  Pick (just locally if necessary)
a volume form for the ambient metric.  Then we have an ambient 
Hodge star $\astar$.  Note that $\astar\e(\X)\ad$ and $\da\astar\e(\X)$
agree up to sign.  Similarly, $\astar\e(\X)$ and $\i(\X)\star$ agree 
up to sign, while $\astar\e(\X)$ gives an isomorphism between
$\Gamma(\Lambda^kT^*\aM/\e(\X))$ and the
subspace of $\Gamma (\Lambda^{n+1-k} T^* \aM)$ consisting
of sections annihilated by $\i(\X)$. It follows now that the results in
the theorem for $\da$ are equivalent to the corresponding results for
$\ad$.

\noindent{\bf Remark:} From these observations it is clear that $\td$
 and $\dt$ do not depend on the ambient construction. In fact, in a
 sense that we will presently describe, they do not depend on the
 conformal structure either. On any $n$-manifold $M$ we may view the
 total space $\tilde{\cq}$ of $\wedge^n T^* M$, with zero section removed,
 as a principal ${\Bbb R}_{\times}$-bundle. Densities of weight $w$ on
 the underlying manifold correspond to functions $f$ on $\tilde{\cq}$
 which are homogeneous in the sense that $\rho^s_* f = |s|^{w/n} f$,
 where $\rho$ denotes the ${\Bbb R}_{\times}$ action. Let $\sim$
 denote equivalence on the total space of $\wedge^kT^*\tilde{\cq}$
 given by $U\sim V$ if $U= \rho^s_* V$ for some $s\in{\Bbb R}_{\times}
 $. Then the quotient $\wedge^kT^*\tilde{\cq}/\sim$ may be viewed as a
 vector bundle $\tilde{\bG}^k$ on $M$ and this has a composition
 series $E^{k-1}\lpl E^k$.  A section $U$ of $\wedge^kT^*\tilde{\cq}$
 which is homogeneous of degree $w$ (i.e.\ $\rho^s_* U = |s|^w U$ at
 each point of $\tilde{\cq}$) is equivalent to a section of the bundle
 $\tilde{\bG}^k[w]=\tilde{\bG}^k\otimes E[w]$,
where $E[w]$ indicates
the bundle of densities of weight $w$ on $M$. 
The exterior derivative
 on $\tilde{\cq}$ preserves homogeneity and so determines an operator
 $\td:\tilde{\bG}^k[w]\to \tilde{\bG}^{k+1}[w]$. If we write
 $\tilde{\bG}_k[w]$ for the bundle $(\tilde{\bG}^k)^*\otimes E[w-n]$
 then we may define $\dt :\tilde{\bG}_{k+1}[-w] \to \tilde{\bG}_k[-w]$
 as the formal adjoint of $\td$. (Alternatively $\dt$ may be obtained from the
 formal adjoint of the exterior derivative on $\tilde{\cq}$.) It is
 not difficult to show that when $M$ has a conformal structure these
 operators agree with the operators $\td$ and $\dt$ defined above.

Given the construction described here, it is clear that there are
 analogues of $\td$ of $\dt$ for many other geometries and situations.
 This will be taken up elsewhere. $\qquad\endrk$

Next we compare $\td$ with the exterior derivative on $M$. More
precisely writing $q_k$ for the natural injection $ \ce^k[w]\to
\cg^k[w]$, we want to compare the compositions $q_{k+1} d$ and $\td
q_k$ as operators on $\ce^k$. We write $\acf:=\acf^1$. Recall that
$2\h(\X,\cdot)=d Q$, so $\acf|_\cq$ is the subspace of $\act$
consisting of sections that take values in $ T\cq\subset
T\aM|_\cq$. The elements of $\acf$ are equivalent to sections of the
tractor subbundle $\bF$.  $T^*\aM/\e(\X)$ is the bundle dual to the
subbundle of $T\aM$ consisting of vectors annihilated by $\i(\X)$.
(We have already observed above that $T^*\aM/\e(\X)|_\cq$ is the
bundle dual to $T\cq\subset T\aM_\cq$, which is what this statement
amounts to
along $\cq$.)  
Via the ambient metric
$T^*\aM/\e(\X)|_\cq\cong T\aM/\e(\X)|_\cq$ (where, as usual, in the
index free notation we are not distinguishing forms from their
contravariant equivalents obtained by the ambient metric).  Also from
above the homogeneous sections of $T\aM/\e(\X)$ of weight $w$ are
denoted $\acg^1(w)$ (with $\acg^1:= \acg^1(0)$) and so $\acg^1(w)$ is
a natural dual space to $\acf(-w)$.
Now  dualising the discussion of Section \ref{tractor} it
is clear that $\pi^*(\ce^1[1])$ is the subspace of $\acg^1|_\cq$
consisting of sections which annihilate vertical fields in $\acf|_\cq$. 
Since the vertical vector fields are generated by $\X$,  this is the subspace of
$\acg^1(1)|_\cq$ annihilated by (the form
field) $\X$. This is well defined since $\X$ is null along $\cq$. So if $\w\in
\ce^1[1]$ then the section of $\acg^1(1)|_\cq$ equivalent to $q_1\w\in
\cg^1[1]$ is exactly $\pi^* \w$. Taking exterior powers and tensoring
with an appropriate density bundle we conclude that similarly the
section of $\acg^k|_\cq$ equivalent to $q_k \phi$, for $\phi\in
\ce^k$, is $\pi^* \phi$.  Now since, as forms on $\cq$, we have $\ad
\pi^* \phi= \pi^* d \phi$ we have the first result of the following proposition. 
The second result here follows immediately from Proposition
\ref{dfa} below.  

\begin{proposition}\label{dq}
As operators $\ce^k\to \cg^{k+1}$ we have
$$
\td q_k =q_{k+1} d .
$$
Similarly 
$$
q^k \dt = \d q^{k+1} 
$$
as operators $\cg_{k+1}\to \ce_k$.
Here $q^k:\bG_k\to E_k$ is the bundle morphism algebraically dual to  $q_k$.
\end{proposition}

The operators $\td$ and $ \dt$ are readily described in terms of
metrics from the conformal class using the machinery from Section
\ref{tractor}. For 
$\si\in \ce[1]$  a choice of conformal scale 
let $\tilde{\si}\in\cce(1)$ be any section such that
$\tilde{\si}|_\cq$ is the homogeneous function equivalent to $\si$.
Define $ \Y:=\aI_\si-\frac{1}{2}(\aI_\si\bullet \aI_\si) \X $ where
$\aI_\si:=\frac{1}{n}\tilde{\si}^{-1}\e(\afD)\tilde{\si} $.  Then $\Y$
is a section of $\act(-1)$ such that $\h(\X,\Y)=1$ and by construction
$\Y|_\cq$ is the homogeneous field along $\cq$ equivalent to the $Y$
corresponding to $\si\in \ce[1]$ (via \nn{Yform}). We make the
definition $\btD:=\nda-\X\otimes \nda_{\miniY}$, where of course $\X$
means the 1-form $dQ/2$.  Note that $\btD_{\miniY}=0$, and that $\btD
Q=0$ and so $\btD$ is a robust tangential operator.  Also if, along
$\cq$, $V\in \frak X(\aM)$ takes values in $T\cq$ then we have (along
$\cq$) $\btD_{V} =\nda_V $. Furthermore if $\tilde F$ is an ambient
tensor field homogeneous of weight $w$ then we have $\btD_{\miniX}
\tilde F=w \tilde F$. With these observations it follows that that
$\btD$ descends to the operator on weighted tractor bundles given by
\begin{equation}\label{tDform}
\tD F=w Y\otimes F+ Z\cdot \nd F
\end{equation}
for $F$ of weight $w$, where
$\nd$ means the coupled tractor--Levi-Civita connection. (This operator
has already played a role in conformal geometry \cite{gosrni,goadv}.)
To see this observe first that
if $w=0$ then, aside from straightforward details, the result is tautological,
given the definition of the tractor connection in terms of the ambient
connection. For other weights observe that upon restriction to vector fields 
with values in $T\cq$ we have $\btD \tilde{F} = \tilde{\si}^w \nda
\tilde{\si}^{-w} \tilde{F} + w \Y \otimes \tilde{F}$.  Since
$\tilde{\si}^{-w} \tilde{F}$ has weight $0$, by the result for the
weight zero case this is equivalent to $\si^{w} Z\cdot \nd \si^{-w} F
+ w Y\otimes F $. But by the definition of the coupled connection this
is exactly $\tD F$ as given.  By construction $\tD$ depends on $\si$.

Since the ambient connection is torsion free we have, for ambient
$k$-forms $F$, $ \ad F=\e(\nda) F.  $ Considering tangential components
of this it follows that for $V\in \cg^k[w]$ we have 
$$
\td V = \e(\tD) V.
$$
The right-hand side here means that we start with any section $V'\in
\ct^k[w-k]$ representing $V$, and then take the class of
$\e(\tD) V'$ in $\cg^{k+1}[w]$. The equality with the left-hand-side
guarantees that this is a well defined operation, with
a result that is independent of the choice of conformal scale $\si$. 
This may also be verified directly by
first noting that $\e(X)\e(\tD)$ is independent of $Y$, and then, via
\nn{connids}, that $Z_{[A}{}^a\nd_a X_{B]}=0$, and thus
$\{\e(\tD),\e(X)\}=\e(Y)\e(X)$.  

The sections 
$\bbY^k$ and $\bbZ^k$ may be used in the obvious way to give scale-dependent
decompositions of the $\cg^k[w]$ bundles.  For example, $\bbY^k
\cdot \alpha +\bbZ^k \cdot \m $ denotes the image in $\cg^k[w]$ of $
\bbY^k \cdot \alpha+\bbZ^k \cdot\m+ \bbW^k \cdot\f+ \bbX^k \cdot \r\in
\ct^k[w-k] $ (for any $\f$ and $\r$),
and $q_k:\ce^k[w]\to \cg^k[w]$ is given by $\mu\mapsto\bbZ^k \cdot\m $.
With these conventions and using the formulae
\nn{trconnform} to calculate $ \e(\tD) V $, we obtain the very simple
formula
\begin{equation}\label{tdform}
\td V=
\bbY^{k+1}\cdot(w\mu-\e(\nd)\alpha)+\bbZ^{k+1}\cdot \e(\nd)\m.
\end{equation}

Next recall that $\dt$ on $\cf^{k+1}[w]$ arises from the action of
$\da = -\i(\nda)$ on $\acf^{k+1}(w)$.  
Now $\i(\nda)=\i(\btD)+\i(\X)\nda_{\miniY}$. So using that
$\nda_{\miniY} \X=\Y$ we have, on $\acf^{k+1}(w)$, that
$$
\i(\nda)=\i(\btD)- \i(\Y).
$$
The operators on the right-hand side here
are both tangential on $\act^{k+1}(w)$.  Thus $\dt$ on $\cf^{k+1}[w]$ is
given by $\i(Y)-\i(\tD)= (1-w)\i(Y)-\i(Z\cdot \nd)$.  We can once
again use \nn{Dform} and \nn{trconnform} to expand this; for
$F=\bbZ^{k+1}\cdot \m+\bbX^{k+1}\cdot \rho\in \cf^{k+1}[w]$ we obtain
$\dt F = \bbZ^k \cdot ((k-n-w)\rho-\i(\nd)\m)+\bbX^k\cdot \i(\nd)
\rho$.  Recall that $\cg_k[w]$ is an alternative notation for
$\cf^k[w+k-n]$, so if now we suppose that $F\in \cg_{k+1} [-w]$, we have
\begin{equation} \label{dtform}
\dt F = \bbZ^k \cdot (w\rho-\i(\nd)\m)+\bbX^k\cdot \i(\nd) \rho. 
\end{equation}
The naturality of $\td$ and $ \dt$ (and the fact that they are first-order
differential operators) is immediate from these explicit formulae, and
using these formulae, it is now straightforward to verify that $\dt$
and $\td$ are formal adjoints with respect to the integral pairing
between $\cg_k[-w]$ and $\cg^k[w]$.  
A further inspection of these formulae also
reveals parts \IT{ii} and \IT{iii} of the next proposition.

\begin{proposition}\label{dfa}
\IT{i}
The operators $\dt: \cg_{k+1}[-w]\to \cg_{k}[-w]$ and
$\td : \cg^k[w]\to \cg^{k+1}[w]$ are natural and are mutual formal adjoints.\\
\IT{ii} When $w\neq 0$, $\td q_k$ is, up to a constant multiple, a differential 
splitting of the canonical surjection $\cg^{k+1}[w]\to \ce^k[w]$ and 
$q^k \dt$ is, up to a constant multiple, 
a differential splitting of the inclusion 
$\ce_k[-w]\to \cg_{k+1}[-w]$.\\
\IT{iii} If $w\neq 0$ then $\cN(\td : \cg^k[w]\to \cg^{k+1}[w])
=\cR(\td : \cg^{k-1}[w]\to \cg^{k}[w])$ and 
$\cN(\dt: \cg_{k}[-w]\to \cg_{k-1}[-w])= 
\cR( \dt: \cg_{k+1}[-w]\to \cg_{k}[-w]) $. 
\end{proposition}

Before Proposition \ref{dq}, we observed that the section of
$\acg^k|_\cq$ equivalent to $q_k \phi$, for $\phi\in \ce^k$, is
exactly the lift $\pi^* \phi$. As a section of $\acg^{k-1}(1)|_\cq$,
$\i(\X)\pi^* \phi$ vanishes and from the considerations above it is
clear that conversely if $\psi\in \acg^k|_\cq$ such that
$\i(\X)\psi=0\in \acg^{k-1}(1)|_\cq $ then $\psi=\pi^* \phi$ for some
$\phi\in \ce^k$. We will extend the use of the term `lift' as
follows. If the restriction of $\tF\in \act^k(w)$
to $\cq$ (i.e.\ $\tF|_\cq$) is, up to a non-zero constant multiple, 
the homogeneous section
equivalent to $F\in \ct^k[w]$ 
then $\tF$ is termed an {\em ambient lift} of $F$ while $\tF|_\cq$
will be called a {\em lift} of $F$.  
For $T\in \cg^k[w]$, the term {\em ambient lift} will be used for ${\tilde
T}\in \acg^k(w)$ with the property that ${\tilde T}|_\cq$ is, 
up to a constant non-zero scale, 
equivalent to $T$;
and also for any representative of this in $\act^k(w)$.  As for the
other cases, ${\tilde T}|_\cq$ will be called a {\em lift} of $T$.  If $T=q_k
t$ for $t\in \ce^k[w+k]$ then ${\tilde T}$ or any representative of
this in $\act^k(w)$ will also be said to be {\em ambient lifts} of
$t$. 

In each case the lift of a section is unique (up to a constant multiple), 
whereas there is choice
in an ambient lift. For most weights $w$ 
(in a sense made precise in part \IT{ii} of the proposition just
below), there is a special ambient
lift of $t\in \ce^k[w]$ so that at least the restriction of 
this to $\cq$ is unique.
First note that it is obvious that acting with $\e(\X)$ gives a 
well defined operator 
from $\Gamma(\Lambda^k T^*\aM /\e(\X)|\cq )\to\Gamma(\Lambda^{k+1} T^*\aM )$
(which is, in fact, injective) and so also 
$\e(\X):\acg^k(w)\to \act^{k+1}(w-k+1)$ 
is well defined. By construction the composition $\i(\afD)\e(\X): 
\acg^k(w)\to \act^k(w-k)$ 
acts tangentially along $\cq$.

\begin{proposition}\label{mainamb}
Let $U\in \act^k(\ell-n/2)$ be an ambient lift of $u\in \ce^k[w]$,
where $w=k+\ell-n/2$. Then \\ \IT{i} $ \i (\X) \i (\afD)\e
(\X) U|_{\cq}=0$; \\ \IT{ii} Up to scale $\i (\afD)\e (\X) U$
is an ambient lift of $ u$ if and only if  $\ell\neq -1$ and $ k\neq \ell+n/2$.
\end{proposition}
\noindent{\bf Proof:}\newline \IT{i} First note that since $U\in
\act^k(\ell-n/2)$ is an ambient lift of $u\in \ce^k[w]$ the image of
$\i(\X)U|_\cq$ in $\acg^{k-1}(w))|_\cq$ must vanish.  Thus
$\i(\X)U=\e(\X)V+QW$ for some ambient forms $V$ and
$W$. The last two form fields are not independent. Since all
forms are smooth and 
$$
0=\i(\X)\i(\X)U=\i(\X)(\e(\X)V+QW),
$$ 
a short calculation gives
$\e(\X)V +\e(X)\i(\X)W=0$, and so $\i(\X)U= \i(\X)\e(\X)W$.  Thus we may
assume without loss of generality that $V=-\i(\X)W$.

Now for an ambient $k$-form $F$ of weight $s$ we have 
$$
\begin{array}{l} 
\i(\afD)F=-(n+2s-2)\da F +\i(\X) \afl F,\qquad\mbox{and so} \\
\i(\X) \i(\afD)
F=-(n+2s-2) \i(\X) \da F = (n+2s-2) \da \i(\X) F.
\end{array}
$$
Thus $\i(\X)\i
(\afD)\e (\X) U = (n+2w-2k)\da\i(\X)\e(\X)U$.  Since $
\i(\X)\e(\X)+\e(\X)\i(\X)$ is $Q$ as a left multiplication operator,
we have $(n+2w-2k)[\da Q U -\da \e(\X)\i(\X) U]$. Recall that
$[\da,Q]=-2\i(\X)$.  Thus 
$$
\begin{array}{l}
\da Q U=-2\i(\X) U +O(Q)= -2\i(\X) \e(\X)W+
O(Q)\qquad\mbox{and} \\
-\da \e(\X)\i(\X) U = -\da \e(\X) \i(\X)\e(\X)W= -\da Q
\e(\X)W \\
\qquad =2 \i(\X) \e(\X) W +O(Q)=O(Q).
\end{array}
$$
So $\i(\X)\i
(\afD)\e (\X) U =  O(Q)$, as required.  \quad $ \Box$\\ 

\noindent\IT{ii} We must show that, up to a non-zero constant multiple, $\i
(\afD)\e (\X) U|_{\cq}$ and $U|_\cq$ represent the same section of
$\acg^k(w)|_\cq$.  Let us re-express $\i (\afD)\e (\X) U$. Using
$\{\da,\e(\X)\}=\LX^*$ and the weight of $U$, we have
$$
\da \e(\X) U=-\e(\X)\da U-(n-2k+w+2)U.
$$
Also $ [\afl,\e(\X)]U=-2\ad U$, so
$$
\begin{array}{rl}
\i(\X)\afl \e(\X)U
&=-2\i(\X)\ad U +\i(\X)\e(\X)\afl U \\
&= -2\cL_{\miniX}U +2 \ad \i (\X) U - \e (\X)\i(\X)\afl U +O(Q) \\
&=-2w U + 2\ad \i(\X)U -\e(\X)\i(\X)\afl U + O(Q).
\end{array}
$$
Combining these yields
$$
\begin{array}{rl}
\i(\afD)\e(\X)U&=
(n+w-2k)(n+2w-2k+2)U+(n+2w-2k)\e(\X)\da U
\\
&\qquad
-\e(\X)\i(\X)\afl U +2\ad\i(\X)U.
\end{array}
$$
As observed above, since $U$ is an ambient lift of $u$, 
$\i(\X)U=\e(\X)V+QW$ where $V$ is a
$(k-2)$-form and $W$ a $(k-1)$-form. Thus $\ad
\i(\X)U=-\e(\X)\ad V + 2\e(\X)W +O(Q)$. As a result, upon restriction to vectors
from $ T\cq$, $ \i(\afD)\e(\X)U$ agrees precisely with
$(n+w-2k)(n+2w-2k+2)U$. This completes the proof, since
$(n+w-2k)(n+2w-2k+2)=2(n/2+\ell-k)(\ell+1)$.
\quad $ \Box$

\noindent{\bf Remark:} The proof of part (i) above can be shortened
somewhat if one first observes that for $u\in \ce^k$ there is an
ambient lift $\tilde{U}$ such that $\i(\X)\tU=0$. For example take
$\tU=\i(\X)\e(\Y)U$ where $U$ is an ambient lift of $u$ as above and as usual
$\Y$ is a section of $ \act(-1)$ such that $\h(\X,\Y)=1$. Since 
$\i(\X)\e(\Y)+\e(\Y)\i(\X)$ is the identity on forms we have
$$
\begin{array}{rl}
\tU&=U- \e(\Y) \i(\X)U\\
&= U +\e(\X)\e(\Y)V -Q \e(\Y)W
\end{array}
$$ from which is clear that $\tU\in \act^{k}(\ell-n/2)$ and $U\in
\act^k(\ell-n/2)$ represent the same section of $\acg^k(w)|_\cq$.  On
the other hand from the last display and Proposition \ref{efdd} it is
also clear that $\i (\afD)\e (\X) U|_{\cq}=\i (\afD)\e (\X) \tU|_{\cq}$.
$\qquad\endrk$

 We observed above that the tangential operator $\i(\afD)\e(\X):
\acg^k(w)\to \act^k(w-k)$ is well defined.  By Proposition
\ref{tracefdd} this descends to a natural conformally invariant
operator $\i(\fD)\e(X):\cg^k[w]\to \ct^k[w]$ for all weights $w$.
Next note that from  part (i) of the last proposition the
composition $\i(\fD)\e(X) q_k$, acting on $\ce^k[w]$, takes values in
the subbundle $\bF^k[w-k]$ of $\bT^k[w-k]$. Part (ii) shows that for
$\phi\in \ce^k[w]$ the image of $\i(\fD)\e(X) q_k \phi $ in $\cg^k[w]$
(under that natural quotient mapping $\bT^k[w-k]\to \bG^k[w]$) is in
general a non-zero multiple of $q_k \phi$. Now recall
$\bG^k[w]=E^{k-1}[w]\lpl E^k[w]$. From the definitions of
$\bF^k[w-k]$ and $\bG^k[w]$ it follows that the composition
$\bF^k[w-k]\to \bT^k[w-k]\to \bG^k[w]$ takes values in the composition
factor isomorphic to $ E^k[w]$. So finally applying $q_k^{-1}$ to
this we obtain a non-zero map $p_k:\bF^k[w-k]\to E^k[w]$ 
and this is a constant multiple of the
canonical surjection $q^k:\bF^k[w-k]\to E^k[w]$.
Gathering these observations and results we have the following.

\begin{proposition} \label{splitop}
The composition 
$$
\i(\fD)\e(X):\cg^k[w]\to \ct^k[w]
$$ 
is a conformally invariant natural differential operator.
If $w=k+\ell -n/2$, then up
to a constant non-zero multiple,
$$
\i(\fD)\e(X) q_k: \ce^k[w]\to \cf^k[w-k]
$$ is a differential splitting of the canonical surjection
$q^k:\cf^k[w-k]\to \ce^k[w]$ if and only if $\ell\neq -1$ and $ k\neq
\ell+n/2$.  That is, $q^k\i(\fD)\e(X) q_k $ is a multiple of the
identity on $\ce^k[w]$ and the multiple is non-zero if and only if 
$\ell\neq-1$ and $k\neq \ell+n/2$ {\rm(}or equivalently $k-1-n/2\ne w
\ne 2k-n${\rm)}.
\end{proposition}   

\section{Proofs of the main theorems and their extensions} \label{mainproofs}

The following set of simple
but remarkable results are central to the constructions which follow.

\begin{lemma} \label{old-domino}
If $ V\in \act^{k}(\ell-n/2+1)$ and $ U\in \act^{k}(\ell-n/2)$ then
for $\ell =0,1,\cdots $ we have
$$
\afl^{\ell} \e(\afD)V =\e(\X) \afl^{\ell+1} V, 
\quad \afl^{\ell} \i(\afD)V =\i(\X) \afl^{\ell+1} V
$$
and
$$
\e(\afD)\afl^{\ell} U = \afl^{\ell+1} \e(\X) U, 
\quad \i(\afD)\afl^{\ell} U = \afl^{\ell+1}\i(\X) U .
$$
Here $ \afl^0$ means 1.
\end{lemma}

\noindent{\bf Proof:} We will prove the first identity; the proofs of the others
are similar.  
First observe that acting on any ambient form field, we have
$$
\afl^{\ell}(2\ell \ad+ \e(\X)\afl) 
=2\ell \afl^{\ell}\ad + \afl^{\ell}\e(\X)\afl 
= 2\ell \afl^{\ell}\ad + [\afl^{\ell},\e(\X)]\afl  + \e(\X) \afl^{\ell+1} .
$$
Now recall from the Tables \ref{anticomms} and \ref{comms}
that $ [\afl,\e(\X)]=-2\ad$, and that $\afl$ and $\ad$ commute.
Thus
$ [\afl^{\ell},\e(\X)]\afl=-2 \ell \afl^{\ell}\ad$, giving
$$
\afl^{\ell}(2\ell \ad+ \e(\X)\afl) = \e(\X) \afl^{\ell+1}.
$$ 
On the other hand, from the definition of $ \e(\afD)$, we have that
$\e(\afD) V =(2\ell \ad  +\e(\X)\afl) V$ for $ V\in
\act^{k-1}(\ell-n/2+1).\qquad\Box$

Recall from Proposition \ref{powerslap} that the powers $\afl^m$
determine
conformally invariant operators $\fl_m:\ct^k[m-n/2]\to\ct^k[-m-n/2]$
for $ m\in\{0,1,2,\ldots\}$. 
Since, via the metric on $\ct^k$, sections of
$\ct^k[m-n/2]$ and $\ct^k[-m-n/2]$ pair to give sections of $\ce[-n]$,
it follows that the formal adjoint of $\fl_m$ is a conformally
invariant differential operator between the same spaces:
$$
\fl_m^*:\ct^k[m-n/2]\to\ct^k[-m-n/2].
$$
For each $m$, let us define the operator
$\fb_m$ to be the average of these, namely
\begin{equation}\label{average}
\begin{array}{l}
 \tfrac{1}{2}(\fl_m+\fl_m^*)=:\fb_m:\ct^k[m-n/2]\to\ct^k[-m-n/2], \\
 \quad m\in 
\{0,1,2,\ldots\}.
\end{array}
\end{equation}
By construction $\fb_m$ is formally self-adjoint.
We now have the following consequence of the lemma above. 

\begin{lemma} \label{domino}
If $ V\in \ct^{k}(\ell-n/2+1)$ and $ U\in \ct^{k}(\ell-n/2)$ then
for $\ell =0,1,\cdots $ we have
$$
\fb_{\ell} \e(\fD)V =\e(X) \fb_{\ell+1} V, 
\quad \fb_{\ell} \i(\fD)V =\i(X) \fb_{\ell+1} V
$$
and
$$
\e(\fD)\fb_{\ell} U = \fb_{\ell+1} \e(X) U, \quad \i(\fD)\fb_{\ell} U 
= \fb_{\ell+1}\i(X) U .
$$
Here $ \fb_0$ means 1.
\end{lemma}

\noindent{\bf Proof:} Recall from Proposition
\ref{tracefdd} that $\i(\fD)$ and $\e(\fD)$ are formal adjoints.
Similarly $\i(X)$ and $\e(X)$ are formal adjoints.
{}From Lemma
\ref{old-domino} we have $\fl_{\ell} \e(\fD) =\e(X) \fl_{\ell+1} $ on
$\ct^{k}(\ell-n/2+1)$ and $\i(\fD)\fl_{\ell} = \fl_{\ell+1}\i(X) $ on
$\ct^{k+1}(\ell-n/2)$.  The formal adjoint of the latter gives
$\fl_{\ell}^* \e(\fD) =\e(X) \fl_{\ell+1}^* $ on
$\ct^{k}(\ell-n/2+1)$. Averaging this with the former gives the first
result. Then $\fb_{\ell} \i(\fD)V =\i(X) \fb_{\ell+1} V$ follows from
a similar argument. Taking formal adjoints on both sides of these two
results then gives the remaining identities.  \quad $\Box$

For $\ell\in {\Bbb Z}$ let us define $\bK^\ell_k:\ct^k[\ell-n/2]\to
\ct^k[-\ell-n/2]$ by 
\begin{equation}\label{bGdef}
\bK^\ell_k = \left\{\begin{array}{ll} \fb_\ell 
\i(\fD)\e(X) & {\rm if} \quad \ell\geq 0, \\
                     \i(X) \e(X)   & {\rm if} \quad \ell= -1,\\
                      0  & {\rm otherwise.} \end{array} \right.
\end{equation} 
Each $\bK^\ell_k$ is a composition of conformally invariant operators
and so is conformally invariant.  Note that for $
U\in\ct^k[\ell-n/2]$, 
Lemma \ref{domino} implies that
$$
\fb_\ell \i(\fD)\e(X) U =\i(X )\fb_{\ell+1} \e(X) U = \i(X )\e(\fD)\fb_\ell U,
$$ 
and so these are each alternative expressions for $ \bK^\ell_k
$. From the form $ \bK^\ell_k =\i(X)\fb_{\ell+1}\e(X) $, it is immediate
that $ \bK^\ell_k$ is formally self-adjoint and that $ \i(X)
\bK^\ell_k=0 $.  Thus $ \bK^\ell_k: \ct^k[\ell-n/2]\to
\ct^k[-\ell-n/2] $ takes values in $\cg_k[-w]\subset\ct^k[-\ell-n/2]
$. (Recall $w:=k+\ell-n/2$.) On the other hand it is also clear 
that the composition
$\bK^\ell_k\e(X) $ vanishes, and so we may naturally view $\bK^\ell_k$
as a formally self-adjoint operator between $ \cg^k[w]$ and $
\cg_k[-w]$. Except where otherwise mentioned, we will take this point
of view.

\begin{proposition} \label{key}
  The expressions of {\rm(}\ref{bGdef}{\rm)} define formally self-adjoint
  conformally invariant differential operators  
$$
\bK^\ell_k:\cg^k[w]\to \cg_k[-w] .
$$
These are natural when $\ell=-1$, or when $\ell$ is in the range of $m$ given 
in Proposition
\ref{powerslap}.
For $V\in \cg^{k-1}[w]$, $ U\in \cg^{k}[w]$ where $n/2+w-k=\ell\geq -1 $ we have
\begin{equation}\label{tagme}
\e(X)\bK^{\ell+1}_{k-1} V = 2(\ell+2) \bK^\ell_k \td V \quad {and} \quad
\bK^{\ell+1}_{k-1} \i(X) U = 2(\ell+2)\dt  \bK^\ell_k U .
\end{equation}
\end{proposition}

\noindent{\bf Proof:} The first statement is established above. The
naturality assertion is immediate from Propositions \ref{powerslap} and
\ref{splitop} since $\i(\fD)$, $\i(X)$ and $\e(X)$ are
natural. 

Next observe that  $\bK^{\ell+1}_{k-1} \i(X) U$ is given by the
expression 
$$
\i(X)\fb_{\ell+2} \e(X) \i(X)U= -\i(X)\fb_{\ell+2} \i(X)\e(X)U.
$$
{}From Lemma \ref{domino}, $\i(X)\fb_{\ell+2} \i(X)\e(X)U=
\i(X)\i(\fD)\fb_{\ell+1}\e(X) U$.  
Now from Proposition \ref{tracefdd}
is is clear that $\i(X)\i(\fD)=\i(X)\i(D)$. On the other hand, given a
choice of scale, we have $\i(X)\i(D)F= -(n+2w-2)\i(X)\i(\tD) F$, for
$F$ any form tractor of weight $w$. Here we have used \nn{Dform} and
\nn{tDform}. 
Thus noting that $\fb_{\ell+1}\e(X)U$ has weight $
-(n/2+\ell+1)$ we have $\i(X)\i(\fD)\fb_{\ell+1}\e(X)U=
2(\ell+2)\i(X)\i(\tD) \fb_{\ell+1}\e(X)U$.  Next, using the formulae
\nn{trconnform} for the tractor connection it is straightforward to verify
that $[\tD,X]=h-X\otimes Y$. Thus we have the identity
$\i(X)\i(\tD)=(\i(Y)-\i(\tD))\i(X) =\dt \i(X)$, and so 
$$
\i(X)\fb_{\ell+2} \e(X) \i(X)U|_{\cq}=2(\ell+2)\dt \i(X)\fb_{\ell+1}\e(X)U|_{\cq},
$$
which is the second identity of \nn{tagme}.
(Note that none of the operators composed on either side depend 
on a choice of scale.)
The first identity of \nn{tagme}
follows immediately by taking formal adjoints.
\quad  $\Box$

\noindent{\bf Remark:} One can establish the identity $\i(X)\i(\tD)=\dt
\i(X)$ without a choice of scale via the ambient metric. Simply
observe that since $-\da=\i(\nda)=\i(\btD)+\i(\X)\nda_{\miniY}$, we
have $\i(\X)\i(\btD)=-\i(\X)\da=\da \i(\X)$. $\qquad\endrk$

\noindent{\bf Remark:} \label{coK} We may also  
define the conformally invariant operators
\begin{equation*}
\bK_{\ell,\star}^{n-k} = \left\{\begin{array}{ll} \fb_\ell \e(\fD)\i(X) 
& {\rm if} \quad \ell\geq 0, \\
                     \e(X) \i(X)   & {\rm if} \quad \ell= -1,\\
                      0  & {\rm otherwise.} \end{array} \right.
\end{equation*} 
on $\ct^{n-k+2}[\ell-n/2]$. By arguments similar to those above, we find
that each $\bK_{\ell,\star}^{n-k}$ descends to a well-defined 
formally self-adjoint 
conformally invariant operator
$$
\bK_{\ell,\star}^{n-k}: \cg_{n-k}^\star[w]\to \cg^{n-k}_\star[-w],
$$
where $\cg_{n-k}^\star[w]$ and $\cg^{n-k}_\star[-w]$ are defined in the 
remark on page \pageref{cosubtractors}. 
In this context we obtain 
$$
\i(X)\bK_{\ell+1,\star}^{n-k+1}= -2(\ell+2)\bK_{\ell,\star}^{n-k}\dt_\star 
\quad \mbox{ and } \quad 
\bK_{\ell+1,\star}^{n-k+1} \e(X) =-2(\ell+2)\td^\star  \bK_{\ell,\star}^{n-k} .
$$ 
Here $\dt_\star: \cg_{n-k+1}^\star[w]\to \cg_{n-k}^\star[w]$ and
$\td^\star : \cg^{n-k}_\star[-w]\to \cg^{n-k+1}_\star[-w]$ are first
order conformally invariant operators which arise, respectively, from the ambient
$\da$ and $\ad$ by constructions parallelling the constructions of $\dt$
and $\td$. 

Suppose now that $M$ is oriented. It is easily verified that on
weighted $k$-form tractors we have
$$
\i(\fD)\tstar =(-1)^k\tstar \e(\fD) \quad \mbox{ and } \quad 
\e(\fD)\tstar =(-1)^{k-1}\tstar \i(\fD) .
$$ 
On the other hand, from the relation between $\tstar$ and the
ambient volume form one obtains that on $\ct^k[\ell-n/2]$, we have
$\tstar \fb_\ell\tstar=(-1)^{k(n+2-k)+q+1}\fb_\ell$. It follows that
up to a sign, $\bK_{\ell,\star}^{n-k}$ is exactly $\tstar \bK^k_\ell
\tstar$. So on general manifolds these operators are, in a suitable
sense, locally equivalent, and we may even 
(ignoring the issue of a possible overall
sign difference) use this to give an alternative definition of  
$\bK_{\ell,\star}^{n-k}$ 
on oriented neighbourhoods as $\tstar
\bK^k_\ell \tstar$. $\qquad\endrk$

Recall that $\cg^k[w]$ has the composition series $\ce^{k-1}[w]\lpl
\ce^k[w]$ and $q_k:\ce^k[w] \to \cg^k[w] $ is the canonical inclusion.
Dually $\cg_k[-w]$ has the composition series
$\ce_k[-w]\lpl\ce_{k-1}[-w]$ with canonical surjection $q^k: \cg_k[-w]
\to \ce_k[-w]$. Thus we can define conformally invariant differential
operators 
by compositions as follows:
\begin{equation}\label{bLdef}
\begin{array}{rl}
(\bK^\ell_kq_k =:\bL^\ell_k):&\ce^k[w]\to \cg_k[-w], \\
(q^k\bK^\ell_k =:\obL^\ell_k):&\cg^k[w]\to \ce_k[-w], \\
 (q^k\bK^\ell_kq_k =:L^\ell_k):&\ce^k[w]\to \ce_k[-w]. 
 \end{array}
\end{equation}
Note that clearly $L^\ell_k = q^k\bL^\ell_k=\obL^\ell_k q_k$. By
construction and from Theorem \ref{key} we have that $ \bL^\ell_k$ and
$\obL^\ell_k$ are 
formal adjoints and that $L^\ell_k$ is
formally self-adjoint.  To simplify the notation, when the source
bundles are true (unweighted)
forms (the case $w=0$ above), and when $\ell=n/2-k$, we shall
often omit the $\ell$ superscript.

A first concern is to verify that these operators are non-trivial. It
suffices to establish this for the family $L^\ell_k$. 

  \begin{proposition} \label{Lelliptic}
    For $ k\in \{0,1,\cdots , n\}$ and $ \ell\in \mathbb{N}$ such that
    $\ell+n/2\neq k $, the operator $ L_k^\ell:\ce^k[w]\to \ce_k[-w] $ is
  conformally invariant, formally self-adjoint, non-trivial and of order $2\ell$. 
It is quasi-Laplacian if and
    only if $k+\ell-n/2 =:w\neq 0$.
In the Riemannian signature case the operator $ L_k^\ell$ is
    elliptic if and only if  $w\neq 0$,
and it is positively elliptic if and only if $k\notin[n/2-\ell,n/2+\ell]$.
For each $k$ the differential operator sequence {\rm(}\ref{detour}{\rm)}
is an elliptic complex.
\end{proposition}

\noindent{\bf Proof:} The claims of conformal invariance and symmetry
under taking adjoints are established above.

Suppose we are in the Riemannian signature
setting.  Since $ \fb_\ell$ is elliptic it has a finite-dimensional
null space on compact manifolds. On the other hand, from Proposition
\ref{splitop}, the range of $ \i(\fD)\e(X) q_k$ is infinite
dimensional. Thus the composition $\fb_\ell \i(\fD)\e(X)
q_k=\bK^\ell_k q_k=:\bL_k^\ell$ is non-trivial in general. From
Proposition \ref{key}, this composition takes values in the subbundle $\bG_k[-w]$
of $ \bT^{k}(-n+k-w)$.

Let us consider $\bL^\ell_k: \ce^k[w]\to \cg_k[-w]$ on the flat model
$S^n$ with its standard conformal structure. Recall that $ \cg_k[-w]$
has the composition series $ \ce_k[-w]\lpl \ce_{k-1}[-w]$ and $ q^k$
is the map onto the quotient $q^k:\cg_k[-w] \to \ce_k[-w]$. If
$q^k\bL^\ell_k=:L^w_k: \ce^k[w]\to \ce_k[-w]$ is trivial then
$\bL^\ell_k$ determines a non-trivial conformally invariant operator
$\ce^k[w]\to \ce_{k-1}[-w]$. There is no such operator \cite{EastSlo}
and so $L^\ell_k$ is non-trivial. (All the invariant operators between
forms preserve $k$ except $d$ and $\d$, or restrictions or projections
of $d$ or $\d$ to $n/2$-forms of one duality.
$\d$ maps
$\ce_k=\ce^k[2k-n]$ to $\ce_{k-1}$, so the only possibility for a true form
operator is when $k=n/2$. 
but this is an $\ell=0$ case,
and so fails the assumption $k\neq \ell +n/2$.) 
By construction,
Proposition \ref{tracefdd} and Proposition \ref{key}, $L^\ell_k$ is
natural in the conformally flat case.

Next observe that it is clear from the formulae for $\bK^\ell$, $q^k$
and $q_k$ and the proof of Theorem \ref{powerslap}  that $L^\ell_k$ has
formally the same leading symbol on all structures (where we range
over both signature and curvature). Thus the $L^\ell_k$ are non-trivial.

Specialising once again to the conformally flat Riemannian case the
differential operators $L^\ell_k$ must be the unique (up to constant
multiples) operators between the bundles concerned.
For the remainder of the proof let us fix some choice of scale. Up to
a non-zero constant multiple, the operator $L^\ell_k $ has the form
\begin{equation*}
\underbrace{(n-2k+2\ell)}_{-2u}(\d d)^\ell+
\underbrace{(n-2k-2\ell)}_{-2w}
(d\d)^\ell+\LOT,
\end{equation*}
and carries $\cE^k[w]$ to $\cE^k[u]$.
This follows from the formulae
for the operators on the sphere (\cite{tbjfa}, Remark 3.30).
In particular $u\ne 0$ and 
$\ce^k[u]=\ce_k[-w]$ is the target bundle of $L^\ell_k$. Thus in
all cases the operators are of order $2\ell$.  Next note that if
$w\ne 0$ then
$$
(-u^{-1}\d d-w^{-1}d\d+\LOT) L^\ell_k = \fl^{\ell+1} +\LOT.
$$
The leading symbol of $d$ is $\eye\e(\x)$. 
Thus if, on the other hand, $w=0$, the
leading symbol of $L^\ell_k$ annihilates the range of $\e(\x)$, and so  
cannot be a right factor of the leading symbol of a power of the
Laplacian.  We conclude that $L^\ell_k$ is quasi-Laplacian if and
only if $w\neq 0$.

Specialising to the Riemannian setting, it follows that $L^\ell_k$ is
elliptic if and only if $w\neq 0$. Using the fact that the leading symbol of
$\d$ is $-\eye\i(\x)$ we have that, up to a non-zero constant multiple,
the leading symbol of $L^\ell_k$ is
$|\x|^{2\ell}(-u\i(\x')\e(\x')-w\e(\x')\i(\x'))$, where
$\x'=\x/|\x|$.  But $\i(\x')\e(\x')$ and $\e(\x')\i(\x')$ are
complementary projections on the fibre $E^k_x$.  
Thus the real linear combination
$-u\i(\x')\e(\x')-w\e(\x')\i(\x') $ is definite if and only
if $w$ and $u$ have the same sign. 
On the other hand if $w=0$ then the leading symbol
of $ L_k$ is, up to a non-vanishing scalar factor, just
$\i(\x')\e(\x')$, and so once again using the fact that $\i(\x')\e(\x')$ and
$\e(\x')\i(\x')$ are complementary projections on the fibre $E^k_x$
(and so also $E_k|_x$), it follows that the symbol sequence
is exact at $\ce^k$ and $\ce_k$. Since the adjoint de Rham sequences
are also elliptic, this shows that the sequence \nn{detour} is an
elliptic complex.  \quad $\Box$

We are now ready for one of the main results.

\begin{theorem} \label{Main}
\IT{i} For $n/2+w-k=\ell\geq 0 $, the operators $\bL^\ell_k$ and $\obL^\ell_k$  
have the
factorisations
$$
\bL^\ell_k = \dt \bN^\ell_k \qquad {\rm and} \qquad \obL^\ell_k = \obN^\ell_k \td
$$ 
where, in a choice of scale, 
$ \bN^\ell_k =2(\ell+1)\bK^{\ell-1}_{k+1}\e(Y)q_k$ and 
$\obN^\ell_k =2(\ell+1)q^k\i(Y)\bK^{\ell-1}_{k+1}$. 
For $n/2-k=\ell\geq 1$, the operator
$L_k$ has the factorisation
$$
L_k=\delta M_k d,
$$
where
\begin{equation} \label{Qtoo}
M^\si_k = -4\ell(\ell+1) q^{k+1}\i(Y)\bK^{\ell-2}_{k+2}\e(Y)q_{k+1}. 
\end{equation}

\IT{ii} The operators $\bL^\ell_k$, $ \obL^\ell_k$ and $ L^\ell_k$ are natural
as follows: in odd dimensions for integers $-1\leq \ell $; in even
dimensions for integers $-1\leq \ell \leq n/2-1$ and for $\ell=n/2$ if
$k=0$. 

\IT{iii} The differential operator $ G^\si_k:\ce^k\to \ce_{k-1}$,
defined on even dimensional manifolds for $ k\leq n/2+1$ by $
G^\si_k:= q^{k-1}\i (Y)\bL_k$, for each choice of conformal scale
$\si$, is natural.  Upon restriction to the null space of $L_k$,
$ G^\si_k$ is conformally invariant {\rm(}and so we
omit the argument $\si${\rm)}. The composition $ G_k d:\ce^{k-1}\to
\ce_{k-1}$ is {\rm(}up to a non-zero scale factor{\rm)} $ L_{k-1}$.

\IT{iv} For $k=1,\cdots ,n/2$, $G_k=\d \tilde{M}_{k-1}$ where
$\tilde{M}_{k-1} : {\Cal N}(L_k)\to \ce_{k}/{\Cal N}(\d)$ is
conformally invariant. $G^\si_0=0$ and, up to a non-vanishing constant
multiple, $G^\si_{n/2}$ is $\d$ (and so is conformally invariant on $\ce^{n/2}$).

\IT{v}  The operators $L_0^\ell:\ce[\ell-n/2]\to \ce[-\ell-n/2]$
are {\rm(}up to a non-zero constant multiple{\rm)} the GJMS operators. 
The operators $\bL^{-1}_k$, 
$\obL^{-1}_k$, $L^{-1}_k$, and $G^\si_{n/2+1}$ all vanish. 
$L^0_k$ is a multiple of the identity, and $L_{n/2}=0$. 
\end{theorem}
\noindent{\bf Proof:}
\noindent Part \IT{i}.  For $u\in \ce^k[w]$ we have $\bL^\ell_ku:=
\bK^\ell_kq_k u$. Now from the definition of $q_k$ and \nn{parts} it
follows that the composition $ \i(X)q_k$ vanishes on $ \ce^k[w']$ (for
any weight $w'$).  So, making an arbitrary choice of scale, we have
$\bK^\ell_kq_k u= \bK^\ell_k \i(X)\e(Y)q_k u $. From Proposition
\ref{key} we have immediately that $\bL^\ell_k = \dt \bN^\ell_k$ with
$ \bN^\ell_k =2(\ell+1)\bK^{\ell-1}_{k+1}\e(Y)q_k$. From this we
obtain $\obL^\ell_k = \obN^\ell_k \td $, with $ \obN^\ell_k$ as given,
by taking formal adjoints.

Next we recall that for $u\in \ce^k$ we have $L_k u=q^k\bL^\ell_k u$
where $\ell=n/2-k$. So from our results just above we have $L_k u =
2(\ell+1) q^k \dt \bK^{\ell-1}_{k+1} \e(Y)q_k u$. Now on $\cg_{k+1}$
we have $q^k \dt=\d q^{k+1} $.  Using that $q^{k+1}\e(X)$ vanishes we
obtain $L_k u= 2(\ell+1) \d
q^{k+1}\i(Y)\e(X)\bK^{\ell-1}_{k+1}\e(Y)q_k u$. Calling on Proposition
\ref{key} then brings us to $4\ell(\ell+1) \d
q^{k+1}\i(Y)\bK^{\ell-2}_{k+2}\td \e(Y)q_{k} u$. Now from the
definition of $\td$ in terms of the ambient exterior derivative and
the relationship of $Y$ in \nn{Yform} to $\Y$ (or alternatively from
\nn{Yform}, Proposition \ref{tracefdd}, and \nn{Dform}) it is
straightforward to verify that as an operator on $\cg^k[w']$ (for any
weight $w'$) we have $\e(Y) = \e(\tilde Y)$ where
$\tilde{Y}:=\si^{-1}\td \si$. (Here $\si$ is the conformal scale
determining $Y$.) It follows immediately that on $\cg^k[w']$ we have
$\{\td, \e(Y) \}=\{\td, \e(\tilde{Y}) \} =0$. From this 
and using that as operators on  $\ce^k$ we have
$\td q_k =q_{k+1} d$ brings us to $L_k u= \d M_k d u $ where
$M^\si_k = -4\ell(\ell+1) q^{k+1}\i(Y)\bK^{\ell-2}_{k+2}\e(Y)q_{k+1}$, 
as claimed.

\noindent Part \IT{ii}.
Since $ \bL^\ell_k= \bK^\ell_kq_{k}$, $ \obL^\ell_k= q^k \bK^\ell_k$
and $ L^\ell_k= q^k\bK^\ell_k q_k$, from Proposition \ref{key} it is
immediate that these operators are natural for $\ell=-1$ and for 
$\ell$ in the range of
$m$ as in  Proposition \ref{powerslap}.  On the other
hand from part (i) above we also have $ \bL^\ell_k
=2(\ell+1)\dt \bK^{\ell-1}_{k+1}\e(Y)q_k $, $ \obL^\ell_k
=2(\ell+1)q^k \i(Y)\bK^{\ell-1}_{k+1}\td$ and $ L^\ell_k
=q^k\bL^\ell_k$. This shows these operators are natural for $
\ell=m+1$ except for the cases $k=0$ and $k=1$.  This exactly yields
the claimed result.  
(Note that by their definitions above each of
these is conformally invariant).

\noindent Part \IT{iii}. Since $ k\leq n/2+1 $ we have $n/2-k=\ell\geq
-1$. Thus from part \IT{ii}, and by construction, the operator $G_k^\si$ is
differential, natural and takes values in $\ce_{k-1}$. Consider $\bL_k
\phi$ for $\phi\in \cN(L_k)$.  Note that $q^k\bL_k \phi= L_k
\phi=0$. Thus $\e(X)\bL_k\phi=0$. 
 Using that $\{\e(X),\i(Y)\}$ is conformally invariant and the
identity on $\cg_k$ we have that the conformally invariant section
$\bL_k\phi\in \cg_k$ is equal to $\e(X)\i(Y)\bL_k\phi $. It follows
immediately that any conformal variation of $\i(Y)\bL_k\phi $ has the
form $\e(X)F$ and so is annihilated by $q^{k-1}$. 
Thus $G_k\phi =q^{k-1}\i(Y)\bL_k \phi$ is conformally
invariant.

Recall that $\bL_k= \bK^\ell_kq_k$ (with $\ell=n/2-k$).
So acting on $\ce^{k-1}$ we have
$2(\ell+2)G_k d = 2(\ell+2)q^{k-1}\i(Y)\bK^\ell_{k}q_k d$. 
Now since, for $v\in \ce^{k-1} $, we have
$q_k d v =\td q_{k-1} v$, Proposition \ref{key} gives
$2(\ell+2)G_k d v = q^{k-1}\i(Y)\e(X)\bK^{\ell+1}_{k-1}q_{k-1} v$. 
Recall that the composition
$q^{k-1}\e(X)$ vanishes on $ \cg_k[w']$ for any weight $w'$, so we have
$$
2(\ell+2)G_k d v = q^{k-1}\bK^{\ell+1}_{k-1}q_{k-1} v = L_{k-1}v .
$$ 

\noindent Part \IT{iv}. From part \IT{iii} we have that $G^\si_k:=
q^{k-1}\i (Y)\bL_k$.  
Since $q^{k-1}\i (Y)$
exactly recovers the coefficient of $\bX^{k}$ in $\cg_{k}$, 
the result is immediate from part \IT{i} and
the expression \nn{dtform} for $\dt$.

\noindent Part ({\it v}). Fix $\ell\in {\Bbb N}$. 
Since $q^0$, $q_0$ are both identity maps we have
$$
L^{\ell}_0 f= \bK^\ell_0 f
$$ for $f\in \ce[\ell-n/2]$. Now via the relationship of $\i(\fD)$ and
$\e(X)$ with the ambient operators $\i(\afD)$ and $\e(\X)$, or via
\nn{Dform} with \nn{trconnform}, it is easily shown that $
\i(\fD)\e(X)f= (\ell+1)(n+2\ell) f $.  Now $\bK^{\ell}_0f=
\fb_\ell\i(\fD)\e(X)f $, and
so
$$
L^{\ell}_0 f =(\ell+1)(n+2\ell) \fb_\ell f.
$$ 
But from \cite{GrZ} (see also \cite{FGrQ}) the operator $\fl_\ell$ 
(of Proposition \ref{powerslap})
is formally self-adjoint on $\ce[\ell-n/2]$ and so $\fb_\ell
f=\fl_\ell f$ is the GJMS operator of order $2\ell$.

Next recall that by definition $\bL^{-1}_k=\bK^{-1}_kq_k$, while
$\obL^{-1}_k=q^k\bK^{-1}_k$ and $\bK^{-1}_k =\i(X)\e(X)$. But
$\{\i(X),\e(X)\}=0$ and $q^k\e(X)=0=\i(X)q_k$.  So $\bL^{-1}_k$ vanishes, and
thus its formal adjoint $\obL^{-1}_k$ must also vanish.  As a result,
$L^{-1}_k=q^k\bL^{-1}_k=0$, and finally
$G^\si_{n/2+1}=q^{n/2}\bL_{n/2+1}=0$ as $\bL_{n/2+1}=q^{n/2}\bL_{n/2+1}^{-1}$.
 
That $L^{0}_k$ is a multiple of the identity, and that this multiple is
zero when $k=n/2$, is shown in Proposition \ref{splitop}, since
$L^0_k=q^k \i(\fD)\e(X)q_k$.  \quad $ \Box$

\begin{proposition} \label{bLelliptic}
For $\ell\in \bN$ and $k\in\{0,1, \cdots ,n\}$, $k\neq \ell+n/2$, the operator
$\bL^\ell_k:\ce^k[w]\to \cg_k[-w]$ is quasi-Laplacian. In particular
in Riemannian signatures it is injectively elliptic.
\end{proposition}

\noindent{\bf Proof:} Recall that $q^k
\bL^\ell_k=L^\ell_k$. 
On one hand, this 
implies that for $w\neq 0$ the result is immediate from Proposition
\ref{Lelliptic}. On the other hand, for the cases $w=0$, using  $q^k
\bL^\ell_k=L^\ell_k$ with $G^\si_k= q^{k-1}\i (Y)\bL_k$, we have that
$$
[\bL_k u]_\si = \left(\begin{array}{l}
                            L_k u \\
                            G^\si_k u
\end{array}
\right) ,
$$ 
in the splitting $[\cg_k[-w]]_\si =\ce_k[-w]\oplus \ce_{k-1}[-w]$ of
$\cg_k[-w]$ determined by a choice of scale $\si$.  
Next we have
already observed in the proof of Proposition \ref{Lelliptic} that, in
a choice of scale, $L_k$ is of the form $(\d d)^\ell +\LOT$ up to a
non-zero constant multiple, while from Theorem \ref{Main} part \IT{iii} we
have that $G_k d$ is $L_{k-1}$ up to a non-zero constant multiple. Using
this and considering possible leading symbols for $G_k$ it follows
that $G^\si_k = a\d (d\d)^\ell+b (d \d)^\ell d +\LOT$ where $a,b\in {\Bbb
R}$ with $a\neq 0$. Thus there is a pair $a',b'\in {\Bbb R}$ giving
$(a' \d d~,~b' d)[\bL_k u]_\si =\fl^{\ell+1} +\LOT$, showing that
$\bL_k$ is quasi-Laplacian.  \quad $\Box$ 

\noindent{\bf Remark:} In the proof we 
have observed that at leading order, $G^\si_k$ has the form 
$a\d (d\d)^\ell+b (d \d)^\ell d$ with $a\neq 0$. From Theorem \ref{Main} we have 
$ \bL^\ell_k = \dt \bN^\ell_k$. Considering also the explicit formula \nn{dtform}
for $\dt$ in a scale 
it follows that $\d$ is a left factor of 
$G^\si_k$, and so $b=0$ and $G^\si_k$ has the form
$$
G^\si_k =a\d ((d\d)^\ell +\LOT ).\qquad\endrk
$$

\subsection{Operators generalising Q-curvature} \label{Qsect}

Let us write $Y_\si$ for the section of $\ct[-1]$ given by $Y=I_\si-
\frac{1}{2}(I_\si\bullet I_\si) X$ where
$I_\si:=\frac{1}{n}\si^{-1}\e(\fD)\si$ and $ \si\in \ce[1]$.  Thus
$Y_\si$ is null and we have $X\bullet Y_\si=1$. By \nn{Yform}, if
$\si$ is a choice of conformal scale then $Y_\si$ is just $Y$ as above
but we want allow the possibility that $\si$ is not (necessarily) a
choice of conformal scale.  Note that the canonical surjection
$\ct^1[-1]\to \cg^1$ maps $Y_\si $ to a section of $\cg^1$. Explicitly
this image is $\tilde{Y}_\si=\si^{-1}\td \si$ (cf.\ the similar
observation for $\tilde{Y}$ in the proof of part \IT{i} of Theorem
\ref{Main}).  Note that as operators on $\cg^k[w']$ we have
$\e(\tilde{Y}_\si)=\e(Y_\si)$, and on $\cg_k[w']$ we have
$\i(\tilde{Y}_\si)=\i(Y_{\si})$.  Thus we shall normally omit the {\em
tilde} and write simply $Y_\si$ for the section in $\cg^1$ given by
$\si^{-1}\td \si$.  We now consider the differential operator $ q^k
\i(Y_\si)\bK^{\ell-1}_{k+1}$ for $\ell\geq 1$.  Apparently this
depends on $ \si$. For any weight $w\in {\Bbb R}$, let us denote by
$\cK^k[w]$ the subspace of $\cg^k[w]$ consisting of $U\in \cg^k[w]$
such that $\td U=0$.

\begin{lemma}\label{lemop}
For each $\si\in \ce[1]$, the composition
$$
q^{k-1} \i(Y_\si)\bK^{\ell}_{k}:\cg^{k}[w]\to \ce_{k-1}[-w] \quad w=k+\ell-n/2
$$ 
is a conformally invariant differential
operator {\rm(}natural for the range of $\ell$ as in
Theorem \ref{Main} part \IT{ii}{\rm)}.
Restricted to 
$\cN(\obL^\ell_k:\cg^{k}[w]\to \ce_k[w])$, it is independent of $\si$.
Thus in particular restricted to 
$ \cK^{k}[w]\subset\cN(\obL^\ell_k) $, $q^{k-1}
\i(Y_\si)\bK^{\ell}_{k} $ is independent of $\si$.
\end{lemma} 

Note that in the first statement here we mean that the operator is
conformally invariant with the choice of $\si\in \ce[1]$ fixed; that
is, we are not linking $\si$ to conformal scale. This point of view
will be continued below.  In addition, by `natural' here we mean, 
natural as an
operator on $E[1]\otimes \bG^k[w]$ (i.e.\ viewing $q^{k-1}
\i(Y_\si)\bK^{\ell}_{k}$ as an operator on $\si$ as well as the
section of $\cg^{k}[w]$.)

\noindent {\bf Proof:} The first statement is clear by construction
and the results above.  Next, recall that the conformally invariant
operator $\bK^\ell_k$ takes values in $\cg_k[-w]$, which has the
composition series $\ce_k[-w]\lpl\ce_{k-1}[-w]$.
The operator $q_k$ is the natural
surjection $\cg_k[-w]\to\ce_k[-w] $, while $q^{k-1}\i(Y_\si)$ is a
splitting of the natural injection $\ce_{k-1}[-w]\to \cg_k[-w]$. Thus
since $\obL^\ell_k:=q^k\bK^\ell_k$, it is immediate that if $U\in
\cN(\obL^\ell_k:\cg^{k}[w]\to \ce_k[w])$, then $q^{k-1}
\i(Y_\si)\bK^{\ell}_{k}U$ is independent of the choice of splitting, 
i.e.\ independent of $\sigma$.

For the final statement observe that when $\ell \geq 0$ we have
$\obL^\ell_k = \obN^\ell_k \td$ (see Theorem \ref{Main} part \IT{i}), 
and so $ \cK^{k}[w]\subseteq\cN(\obL^\ell_k) $.  On the other
hand $\obL^{-1}_k=0$ (Theorem \ref{Main}, part \IT{v}) so the final
statement follows trivially in this case. \quad $\Box$ 

\noindent{\bf Remark:} Note that the operator $q^{k-1}
\i(Y_\si)\bK^{\ell}_{k}:\cg^{k}[w]\to \ce_{k-1}[-w]$ is essentially a
generalisation of $G_k^\si$ and has many properties which reflect
this.  In particular note that on $\cg^{k-1}[w]$ we have
$2(\ell+2)q^{k-1} \i(Y_\si) \bK^\ell_k \td = q^{k-1}
\bK^{\ell+1}_{k-1} =\obL^{\ell+1}_{k-1} $ (cf.\ part \IT{iii} of
Theorem \ref{Main}). This uses the result $2(\ell+2)\bK^\ell_k \td =
\e(X)\bK^{\ell+1}_{k-1} $ of Theorem \ref{key}. From the latter it is clear that 
$\e(X)$, as well as $\i(X)$, annihilates $\bK^\ell_k \td $. This in turn 
implies that 
$q^k \bK^\ell_k \td =0$ and so 
\begin{equation}\label{Kdform}
2(\ell+2)\bK^\ell_k \td = \bX^{k}\cdot \obL^{\ell+1}_{k-1}.
\end{equation}
This is useful in the next section. $\qquad\endrk$

The next proposition constructs a family of operators with an
interesting conformal transformation property. 
In this, as in the lemma above, `natural' 
means natural as an operator on $\si$ as well as
the section in $\cg^{k}[w]$.

\begin{proposition} \label{nuQ}
For each choice of $\sigma\in \ce[1]$ the differential operator
$$
(\bQ^{\ell,\sigma}_k:=
{-2(\ell+1)}q^{k}\i(Y_\si)\bK^{\ell-1}_{k+1}\e(Y_\si)):\cg^k[w]\to
\ce_{k}[-w], \quad w=k+\ell-n/2
$$ 
is conformally invariant {\rm(}and 
natural for $\ell-1$ as in the range of $\ell$ in
Theorem \ref{Main} part \IT{ii}{\rm)}. Acting on $U\in \cK^k[w]={\Cal
N}(\td:\cg^k[w]\to \cg^{k+1}[w])$, $\bQ^{\ell,\sigma}_k$ has the
transformation law
$$
\bQ^{\ell,\hat\sigma}_kU=\bQ^{\ell,\sigma}_kU + \obL^\ell_k (\Up U)
$$
where $\hat\sigma=e^{-\Up}\sigma$. 
\end{proposition}

\noindent{\bf Proof:} The first statement is clear from the definition of
$\bQ^{\ell,\sigma}_k $.

Let us pick sections $\si_1,\si_2\in \ce[1]$.  Viewing $ Y_{\si_2}$ as
a section of $\cg^1$, we have $Y_{\si_2}=\si_2^{-1}\td \si_2$ and so
it is clear that $\td Y_{\si_2}=0$. 
Thus if $U\in \cK^k[w]$ then
$\e(Y_{\si})U\in \cK^{k+1}[w]$ and it follows at once from Lemma
\ref{lemop} that $ q^{k}\i(Y_{\si_1})\bK^{\ell-1}_{k+1}\e(Y_{\si_2})U$
is independent of $\si_1$.

Now let $\hat\si_2=e^{-\Up}\si_2$ for some smooth function $\Up$.
Viewing $Y_{\si_2}$ and $Y_{\hat\si_2} $ as sections of $\cg^1$ we have
$Y_{\hat\si_2}=Y_{\si_2}- \td \Up$. Thus
$$
q^{k}\i(Y_{\si_1})\bK^{\ell-1}_{k+1}\e(Y_{\hat\si_2})U 
- q^{k}\i(Y_{\si_1})\bK^{\ell-1}_{k+1}\e(Y_{\si_2})U 
= -q^{k}\i(Y_{\si_1})\bK^{\ell-1}_{k+1}\e(\td \Up)U .
$$
Since by assumption $ \td U=0$ we have $ \e(\td \Up)U=\td(\Up U)$
and so by Proposition \ref{key},
$$
\begin{array}{rl}
2(\ell+1)q^{k}\i(Y_{\si_1})\bK^{\ell-1}_{k+1}\e(\td \Up)U
&=q^{k}\i(Y_{\si_1})\e(X)\bK^\ell_{k} (\Up U) \\
&=q^{k}\bK^\ell_{k} (\Up U) = \obL^\ell_k (\Up U).
\end{array}
$$
Here we have used the operator equality  $\{ \i(Y_\si),\e(X)\}=X\cdot Y_\si =1$.
\quad $ \Box$

We now return to the convention that $\si$ denotes a conformal scale.

\noindent{\bf Definition:} For each choice of {\em conformal scale} 
$\si\in \ce[1]$ on even
dimensional conformal manifolds and for each $k\le n/2$, we
let
$(Q^\si_k:=\bQ^{\ell,\si}_kq_k):\ce^k\to \ce_k$ be given by
\begin{equation}\label{Qdef}
 Q^\si_k=-2(\ell+1)q^{k}\i(Y)\bK^{\ell-1}_{k+1}\e(Y) q_k.
\end{equation}  

The operator $Q^\si_k$ 
has the properties described in Theorem \ref{Qthm}, and in particular
 is a generalisation of Branson's $ Q$-curvature.

\noindent{\bf Proof of Theorem \ref{Qthm}:} Note that by construction $ 2(\ell+2)
Q_k^\si =M^\si_{k-1}$, where by the right hand side we mean the operator
given in \nn{Qtoo} above (viewed as an operator $\ce^k\to
\ce_k$). Thus Part \IT{iii} is already contained in Theorem
\ref{Main}.

Part \IT{i} is immediate from formula \nn{Qdef}, since
$\bK^{\ell-1}_{k+1}$ is formally self-adjoint by Theorem \ref{key}.

Part \IT{ii}. We have
$$
\d Q^\si_k =
-2(\ell+1)\d  q^{k}\i(Y)\bK^{\ell-1}_{k+1}\e(Y)q_k.
$$ 
Once again recall that
on $\cg^{k}$, $\d q^k=q^{k-1}\dt$ while as operators on
$\cg^{k+1}[w']$ for any weight $w'$,
we have $\{\dt,\i(Y)\}=0$.
Thus using Proposition \ref{key} and $\i(X)q_k=0$, we get
$$
\d Q^\si_k = q^{k-1}\i(Y)\bK^\ell_{k}q_k = G^\si_k ,
$$ 
where $G^\si_k$ is as defined in Theorem \ref{Main} (and its
restriction to $\cN(L_k)$ is denoted $G_k$).

Part \IT{iv}. If $u\in \cC^k $ then $q_k u\in \cK^k[0]$, since $ \td
q_k=q_{k-1} d$. Thus the transformation law is immediate from the
definition of $ Q^\si_k$ and Proposition \ref{nuQ}, since $L_k=\obL_k
q_k$, and $q_k$ commutes with the multiplication operator
$\e(\Upsilon)$, for any function $\Upsilon$.

Part \IT{v}. We have $L_k^\ell=q^k \bK^\ell_k q_k =
\obL^\ell_kq_k$.  Specialising to the case of densities $ \ce[w]$, observe
that $q^0$ and $q_0$ are both simply identity maps.  So for $f\in
\ce[w]$ we have $L^\ell:= L_0^\ell= \obL^\ell_0$, where
$w=\ell-n/2$. From Theorem \ref{Main} part \IT{v}, this is a GJMS operator.
{}From part \IT{i} of that theorem we also have 
$L^\ell= 2(\ell+1)\i(Y)\bK^{\ell-1}_1
\td f$.  Let us choose a conformal scale $\si$. Then on densities of 
weight $w$ we
have $\td =Z\cdot\e(\nd) +w \e(Y)$, from \nn{tdform}.  So
for a function $ f\in \ce[0]$ we have the operator
$$
\begin{array}{rl}
\si^{-w} L^\ell \si^w f &= 2(\ell+1)\si^{-w}\i(Y)\bK^{\ell-1}_1 
Z\cdot \e(\nd) \si^w f \\
&\qquad
+ w
2(\ell+1)\si^{-w}\i(Y)\bK^{\ell-1}_1 \e(Y) \si^w f .
\end{array}
$$
Note that the first term on the right-hand side annihilates
constant functions and so setting $f=1$, 
taking the coefficient of $w$ and setting in
this $\ell=n/2$ (i.e.\ $w=0$) yields (by definition) Branson's
Q-curvature. Thus Branson's $Q$ is given, in dimension $2\ell$, by
the operator $ 2(\ell+1)\i(Y)\bK^{\ell-1}_1 \e(Y) 1 $. 
But from the definition \nn{Qdef} 
this is
exactly $ -Q^{\si}_0 1$.
\quad $\Box$

\subsection{Other constructions and operators of order $n$}
In even dimension $n$, Theorem \ref{Main} constructs natural
conformally invariant differential 
operators $L^\ell_k$ up to order $n-2$. For
$k=0$ we have that $L^{n/2}_0$ is natural but Theorem \ref{long}
asserts the existence of curved generalisations of the conformally
flat operators at order $n$ for other $k$ values.  We obtain these by
a variation on our general construction. Note that for the operators
of order 4, the observation that one needs, and that there exist, such
alternative constructions is detailed in \cite{Grnon}.

Note that by the formula \nn{tdform} for $\td$ we have that, acting on
$\ce^k[w]$, $\i(X)\td q_k=wq_k$.  (Alternatively observe that on
$\aM$, if $U\in \act^k(w-k)$ has the property $\i(\X) U=O(Q)$, then
$\i(\X)\ad U=\LX U=wU$; the result follows.) Thus with
$w=k+\ell-n/2$, we have $w\bL^\ell_k= \bK^\ell_k \i(X)\td q_k$, as an
operator on $\ce^k[w]$.  So by Proposition \ref{key} we have $
w\bL^\ell_k= 2(\ell+1)\dt \bK^{\ell-1}_{k+1} \td q_k $.
(Integrating this by parts gives an alternative formula along these lines for
the formal adjoint $\obL^\ell_k$.)  
 For the
cases where $w\neq 0$ this provides an alternative construction of
$\bL^\ell_k$, and thus also of $L^\ell_k$:
$$
wL^\ell_k= 2(\ell+1)q^k\dt \bK^{\ell-1}_{k+1} \td q_k .
$$ 

Next from \nn{Kdform} we have that $2(\ell+1)\bK^{\ell-1}_{k+1} \td q_k =
\bX^{k+1}\cdot L^{\ell}_{k}$. 
So, from the explicit formula \nn{dtform} for $\dt$, 
it follows that when $w\neq 
0$, the action of $q^k\dt$ here is the same as some non-zero multiple
of $\bbY^{k+1} \bullet$.

To further re-express  $L^\ell_k$ we need the following lemma.
\begin{lemma}
The conformally invariant differential operator
$$
\i(\fD)\e(X)\td q_k: \ce^k[k+\ell-n/2]\to \ct^{k+1}[\ell-n/2-1]
$$
is a differential splitting of the canonical conformally invariant surjection 
$$
\bbX^{k+1} \bullet : \ct^{k+1}[\ell-n/2-1] \to \ce^k[k+\ell-n/2]
$$
for values of $k$ and $\ell$ such that $\ell\neq \pm 1$ and $k\pm \ell \neq n/2$.
\end{lemma} 

\noindent{\bf Proof:} First note that it is clear from \nn{parts}
that
$\bbX^{k+1}\bullet\i(\fD)\e(X)\td
q_k$ is same as $\bbZ^k\bullet\i(X)\i(\fD)\e(X)\td q_k$.  A
straightforward calculation shows that on form tractors of weight
$\tilde{w}$ we have
$(n+2\tilde{w}-2)\i(\fD)\i(X)+(n+2\tilde{w}+2)\i(X)\i(\fD)=0$. Thus
acting on $\ce^k[w]$ we have $2(\ell+1) \i(X)\i(\fD)\e(X)\td
q_k=-2(\ell-1)\i(\fD)\i(X)\e(X)\td q_k$.  Now from
$\e(X)\i(X)+\i(X)\e(X)=0$ and our observation above that $\i(X)\td
q_k=wq_k$ we obtain
$$
2(\ell+1) \i(X)\i(\fD)\e(X)\td q_k= 2(k+\ell-n/2)(\ell-1)\i(\fD)\e(X)q_k
$$ 
and the result follows immediately from Proposition \ref{splitop}
since on $\cf^k[w-k]$, $q^k$ is a non-zero constant multiple of
$\bZ^k\bullet$.  \quad $\Box$

Let us suppose that the integers $k,\ell$ are as in the lemma
above. From the lemma and the fact that $ \bK^{\ell-1}_{k+1} \td q_k$ is a
non-zero multiple of $\bX^{k+1}\cdot L^\ell_k $, we see that if $u,v\in
\ce^k[w]$ then
$$
( \i(\fD)\e(X)\td q_k v) \bullet( \bK^{\ell-1}_{k+1} \td q_k u)
$$ is a non-zero multiple of $v\cdot L^\ell_k u $. Integrating by
parts, it follows immediately that 
$$
q^k \dt \i(X)\e(\fD)\bK^{\ell-1}_{k+1} \td q_k =  
q^k \dt \i(X)\e(\fD)\fb_{\ell-1}\i(\fD) \e(X) \td q_k
$$ 
is a non-zero multiple of $L^\ell_k$ on $\ce^k[w]$.  
>From Proposition \ref{key} this is
natural for $\ell$ exactly as in Theorem \ref{Main} part (ii). The
importance of this expression, for our current purposes, is that this
formula is easily modified. In the constructions above we have used that 
$\fb_{\ell-1}\i(\fD) \e(X) \td q_k$ takes values in $\cg_{k+1}[-w]$.
Since $q^k \dt \i(X)\e(\fD)$
acts invariantly  on general form tractors we can, in the 
right-hand side of the last 
display, replace $\fb_{\ell-1}$ with 
$$
\Box_{\ell-1} =D^{A_1}\cdots D^{A_{\ell-2}}\Box  D_{A_{\ell-2}}\cdots D_{A_1} .
$$ 
For any integer $\ell\geq 2$ this is a natural, conformally invariant and
formally self-adjoint differential operator
on any tractor bundle
of weight $\ell-1-n/2$.  If $2\leq \ell \leq n/2$, then $\Box_{\ell-1}$ 
has leading term a
non-zero multiple of $\Delta^{\ell-1}$ (see
e.g. \cite{GoPet,gosrni}).  It follows that on $\ce^k[w]$,
$$
\tilde{L}^\ell_k:= q^k \dt \i(X)\e(\fD)\Box_{\ell-1}\i(\fD) \e(X) \td q_k
$$ 
gives an invariant operator $\tilde{L}^\ell_k: \ce^k[w]\to
\ce_k[-w]$ which at leading order agrees with $L^\ell_k$ provided
$0\neq w=k+\ell-n/2$, $\ell\leq n/2$, $\ell\neq \pm 1$, and
$k \neq \ell+n/2$.  In particular, since the splitting operator
$\i(\fD) \e(X) \td q_k$ and its formal adjoint $q^k \dt \i(X)\e(\fD)$
are natural, we have the following result.

\begin{theorem} \label{ordern}
For each $0<k<n$, the operator 
$$
\tilde{L}^{n/2}_k:\ce^k[k]\to \ce_k[-k]
$$ is  natural, conformally invariant, formally self-adjoint, and of
order $n$.  It is quasi-Laplacian. In the Riemannian signature
case the operator $ \tilde{L}_k^{n/2}$ is elliptic.
\end{theorem}

\section{Nontriviality of the cohomology maps} \label{nont}

In Theorem \ref{pairings} we have shown that the operators $Q_k$ give
conformally invariant maps $ Q_k: {\Cal H}^k\to H_{k}^L(M)$. Here we
demonstrate that these maps are not
trivial in general. Clearly it is sufficient to show the stronger claim
that the $Q_k: {\Cal H}^k\to H_{k}(M)$ of Corollary \ref{corpairings}
are non-trivial. For $k=0$ this result is already well known. It boils
down to checking the same question for the Q-curvature, but from
\cite{tomsharp} this integrates to $(n-1)!{\rm vol}(S^n)/2$ times
the Euler characteristic for 
conformally flat structures.
Here we show non-triviality of $Q_{p}: {\Cal H}^p\to H_{p}(M) $ for 
$M=S^p\times S^q$, where $p=n/2-1$, $q=n/2+1$, with
the standard Riemannian conformal structure.

A straightforward expansion of \nn{Qdef} shows that up to a non-zero
constant multiple, $Q^\si_p$ is given by $\frac12d\d+\J-2\V\sharp$.
(Recall that $\J$ is the trace of the Schouten tensor $\V$.)  For
our purposes here let us write $Q:
=\frac12d\d+\J-2\V\sharp$.  First we study $\cH^p$.

\begin{proposition} \label{Hp}
Let $M=S^p\times S^q$, where $p=n/2-1$, $q=n/2+1$, with
the standard Riemannian structure, and let
$Q=\frac12d\d+\J-2\V\sharp$.  Then $\f\in\cE^p(M)$ is in the joint
null space of $\d d$ and $\d Q$ if and only if $\f$ is harmonic.
\end{proposition}

{\bf Proof}: First note that by compactness, 
$\d d\f=0\implies d\f=0$. 

Let $\f$ be in the joint null space described above; then 
$d\f=0$, $\d Q\f=0$.
Let $\U:=\J-2\V\sharp$.  Then 
\begin{equation}\label{estimate} 
0=4\|\d Q\f\|^2=\|\d d\d\f\|^2+4\langle \d d\d\f,\d \U\f\rangle+4\|\d \U\f\|^2. 
\end{equation} 
Here $\langle \cdot,\cdot \rangle$ and $\|\cdot\|$ are the 
$L^2$ inner product and norm. 
The first term on the right in \nn{estimate} is 
$\langle \Delta\f,\Delta^2\f\rangle $, 
since $d\f=0$. 
The second term may be written 
$$ 
\langle d\d\f,d\d \U\f\rangle=\langle d\d\f,\Delta \U\f\rangle 
=\langle \Delta\f,\Delta \U\f\rangle. 
$$ 
The form Laplacian $\Delta$ commutes with the projections from 
the decomposition of $\cE^p(M)$ as $\oplus_{r+s=p}\cE^r(S^p)\boxtimes 
\cE^s(S^q)$, where $\boxtimes$ is the external tensor product. 
Furthermore, $\U$ is a linear combination of these projections, since 
$4\U$ takes the eigenvalue 
\begin{equation}\label{numb} 
\dfrac{16r}{n-2}-\dfrac{2(n-4)}{n-1} 
\end{equation} 
on $(r,s)$-forms.  As a result, $\U$ commutes with $\Delta$. 
Thus the second term in \nn{estimate} may be 
written 
$(\Delta\f,\U\Delta\f)$. 
Since the third term in \nn{estimate} is nonnegative, we have 
\begin{equation}\label{deceight}
0=4\|\d Q\f\|^2\ge\langle\Delta\f,(\Delta+4\U)\Delta\f\rangle. 
\end{equation}

If $n=4$, \nn{numb} shows that $\U$ is a nonnegative operator.
Since $\Delta$ is also nonnegative, \nn{deceight} shows that
both $\Delta$ and $\U$ kill $\Delta\f$.
Since $\Delta$ is a positive operator on its range, we have 
$\Delta\f=0$; i.e.\ $\f$ is harmonic.
In general even dimension $n\ge 4$,
\nn{numb} shows that the eigenvalues 
of $4\U$ are $>-2$, so that $\Delta+4\U$ can be 
nonpositive only on eigenspaces of 
the Laplacian with eigenvalue $<2$. 
But the non-zero 
eigenvalues of the form Laplacian on standard $S^m$ (see e.g.\ \cite{GM,IT})
are integers $\ge m$. 
Now the eigenvalues of the form Laplacian on $M$ are sums of 
form Laplacian eigenvalues on $S^p$ and $S^q$.  Thus non-zero 
eigenvalues less than 2 of $\Delta$ can only arise for $n/2-1\le 1$. 
This 
means that either we are in
the 4-dimensional case treated above, or else $\Delta\f$ is harmonic
(and thus vanishes), so that $\f$ is harmonic. 
This shows that the joint null space of $\delta d$ and $\delta Q$
is contained in the harmonics.

For the opposite inclusion, note that if $h$ is a harmonic $(r,s)$-form, then
$$
2\delta Qh=\delta(d\delta+{\rm const})h=
(\delta d+{\rm const})\delta h=0.
$$
But the $(r,s)$ components of a harmonic form $\f$ are harmonic, so
$\delta Q$ kills harmonics.$\qquad\Box$
 
We now use the above setup to give an example of a situation
in which the cohomology map $Q_p$ 
is non-trivial.  The cohomology of $M=S^p\times S^q$ is 
1-dimensional in the orders $0,p,q,n$, with $\omega$, the pullback of 
the $S^p$ volume form under projection onto the $S^p$ factor, generating 
the harmonics $H^p$
(as well as $H_p$). 
By the above, $\cH^p=H^p$, so
$\omega$ generates $\cH^p$. 
Since $\nabla\omega=0$, we have $Q\omega=\U\omega$, and \nn{numb} shows 
that $\U$ takes the eigenvalue $3n/2(n-1)$ on $(p,0)$-forms. 
Thus $Q\omega=3n\omega/2(n-1)$.
In particular, we have:

\begin{theorem}\label{Qp}
Let $M=S^p\times S^q$, where $p=n/2-1$, $q=n/2+1$, with
the standard  Riemannian conformal structure. Then 
$Q_p:\cH^p\to H_p(M)$ is non-trivial. 
\end{theorem}

\section{Variations on the theme of $Q$} \label{var}

Proposition 2.8 of \cite{GoPet} described one way to proliferate
natural scalar fields with transformation properties similar to the
Q-curvature in the sense that their conformal variation is by a
linear conformally invariant differential operator acting on the
variation function $\Up$.
Our first objective here is to observe that this generalises. We
follow this by making some connections with other recent constructions
of the Q-curvature.
 
\subsection{Semi-invariant operators}\label{semi}
Recall that as an operator on $\cg^{k}[w]$ we have that $\e(Y)$ is
the same as $\e(Y_\si)$ where $Y_{\si}:=\si^{-1}\td \si $.  Under a
conformal transformation given by $\si\mapsto \hat{\si}=e^{-\Up}\si$,
we thus have $\e(Y_{\hat{\si}})=\e(Y_\si)-\e(\td \Up)$. 
Now for each natural
conformally invariant operator $ S:\ct^{k+1}[w-k-1] \to \ct^{k'+1}[w'+k'-n+1]$ and
choice of conformal scale $\si$ there is an invariant operator
$$
(\bS^\si :=q^{k'}\dt \i(X)\e(\fD) S \i(\fD)\e(X)\e(Y_\si)):\cg^k[w]\to
\ce_{k'}[w'] .
$$ 
{}From the transformation law for $\e(Y_\si)$ it follows that upon
restriction to $\cK^k[w]$, our operator has the conformal
transformation (cf.\ Proposition \ref{nuQ})
$$
\bS^{\hat{\si}}=
\bS^\si - 
q^{k'}\dt \i(X)\e(\fD) S \i(\fD)\e(X) \td \e(\Up) ,
$$ 
where (as above) $\e(\Up)$ is $\Up$ viewed as a multiplication
operator.  Note that $ q^{k'}\dt \i(X)\e(\fD) S \i(\fD)\e(X) \td $ is
a composition of conformally invariant operators giving an operator
$\cg^k[w]\to \ce_{k'}[w']$. Composing $\bS^\si$ with $q_k$ and
restricting to $\ce^k$ we obtain the following result.

\begin{proposition} \label{Qopthings}
For each natural conformally invariant operator 
$$ 
S:\ct^{k+1}[-k-1] \to
\ct^{k'+1}[w'+k'-n+1]
$$ 
and choice of conformal scale $\si$, there is a natural invariant
operator
$$
\bS^\si q_{k}:\ce^k\to
\ce_{k'}[w'] .
$$
Upon restriction to $\cC^k$ this has the conformal transformation
$$
\bS^{\hat{\si}} q_k=
\bS^\si q_k - 
q^{k'}\dt \i(X)\e(\fD) S \i(\fD)\e(X) q_k d \e(\Up) .
$$
Acting  between $\ce^k$ and $\ce_{k'}[w']$, the natural differential operator
$$
q^{k'}\dt \i(X)\e(\fD) S \i(\fD)\e(X) q_k d
:\ce^k\to \ce_{k'}[w']
$$
is conformally invariant. If $k=k'$ and  $w'=0$ then we may re-express 
this by the formula
$\d q^{k} \i(X)\e(\fD) S \i(\fD)\e(X) q_k d$, and this is
is formally self-adjoint if $S$ is.
\end{proposition}
The final statement on formal self-adjointness is clear from the
symmetry of the formula and our earlier observations (identifying
$\i(\fD)$ as the formal adjoint of $\e(\fD)$ and so forth).  

Of course
the operators $\bS^\si$ are most interesting in the cases where $k=k'$
and $-w'=w=k+\ell-n/2$ for some $\ell\in {\Bbb N}$, since then they may
be added to $\bQ^{\ell,\si}_k$ without altering its properties
significantly.  If in addition 
$n$ is even and $w=w'=0$, then $\bS^\si q_k$ operates between
$\ce^k$ and $\ce_k$, and so similarly provides a possible modification
to the operator $Q_k^\si$. 
Note that if we take such an
operator, with $S$ also formally self-adjoint, and form the new
``Q-operator'' $Q^\si_k+\bS^\si q_k $, then this satisfies parts
\IT{i}, \IT{ii} and \IT{iii} of Theorem \ref{Qthm}. Part \IT{iv} of
that theorem also holds with the qualification that  
the invariant operator in the conformal
variation formula is a modification of $L_k$ by the 
addition of
a constant multiple of
the operator $-\d q^{k} \i(X)\e(\fD) S \i(\fD)\e(X) q_k d$. (Note that
$L_k-\d q^{k} \i(X)\e(\fD) S \i(\fD)\e(X) q_k d $ has the general form
\nn{mk}.) Finally (as pointed out in \cite{GoPet})  $(Q^\si_0+\bS^\si
q_0 )1$ gives an alternative to Branson's Q-curvature.

It is easy to construct non-trivial
examples. For example one can take $S$ to be $|C|^2$ (where $C$ is the
Weyl curvature), viewed as multiplication operator, to obtain
$$
q^{k'}\dt \i(X)\e(\fD)|C|^2 \i(\fD)\e(X)\e(Y_\si)):\cg^k[w]\to
\ce_{k}[w'],
$$ where $w'=w-2k+n-6$. Specialising to even dimensions, $w=0$,
$k=n/2-3$ and composing with $q_k$ we obtain
\begin{equation*}
\d q^{n/2-3} \i(X)\e(\fD)|C|^2 \i(\fD)\e(X)\e(Y_\si))q_{n/2-3}: 
\ce^{n/2-3} \to \ce_{n/2-3}.
\end{equation*}
Using the explicit formulae for $\e(\fD)$, etc.\ it is easy
expand this and verify that in (even) dimensions $n\geq 6$ this is
non-trivial. (In dimension 6 this boils down to a case treated this
way in \cite{GoPet}.)

Finally on this point we should remark that in constructing $\bS^\si$
we have have made no attempt to produce the most general object with a
transformation law similar to the $Q_k$ operators. 
Since, on $\cg^{k}[w]$,
 $\e(Y_\si)$ has the conformal transformation
$\e(Y_{\hat{\si}})=\e(Y_\si)-\e(\td \Up)$ it follows that {\em any}
conformally invariant operator $P$ which acts on $\cg^{k+1}[w]$ (for
any weight $w$) may be composed with $\e(Y_\si)$ to yield an operator
with a similar transformation law to the $\bQ$ operators. In the case
$w=0$ we may form the composition $P\e(Y_\si)q_k$, which has a
transformation law similar to that of the $Q_k^\si$ operators. 
However we envisage that the
main interest should be the $\bS^\si q_k$ that operate between $\ce^k$
and $\ce_k$, since these may play a role in understanding the nature of
the $Q$ operators and the Q-curvature.

\subsection{Other recent constructions of Q-curvature}
Up to a scale our construction here gives Branson's curvature $Q$ (as
a multiplication operator) by $\i(Y_\si)\bK^{n/2-1}_{k+1}\e(Y_\si)1=
\i(Y_\si)\fb_{n/2-1}\i(\fD)\e(X)\e(Y_\si)1$. Since $\fl_\ell$ is
formally self-adjoint on densities \cite{GrZ,FGrQ} the argument above
leading to this conclusion (in the proof of part \IT{v} of Theorem
\ref{Qthm}) works equally if we start with $\fl_\ell$ in place
$\fb_\ell$ in the formulae for the $L^\ell$. Thus up to a multiple,
$Q^\si_0$ is given by $ \i(Y_\si)\fl_{n/2-1}\i(\fD)\e(X)\e(Y_\si)$.
Applying this to the constant function 1, we have
$\i(Y_\si)\fl_{n/2-1}\i(\fD)\e(X)\e(Y_\si)1=
\i(Y_\si)\fl_{n/2-1}\i(\fD)\e(X)Y_\si$.  
Suppose we define $I^g$ by
$nI^g:=\i(\fD)\e(X)Y_\si$, where $g$ is the metric corresponding to
the conformal scale $\si$.
Then
an elementary calculation
using $Y_\si=\si^{-1}\td \si$ and the formulae in Section
\ref{tractor} shows that 
$I^g$ is exactly the scale dependent
standard tractor defined and used in Section 2.3 of \cite{GoPet} to
give a new construction of the Q-curvature.  So in this notation the
Q-curvature is given, up to a non-zero constant multiple, by
\begin{equation}\label{GPQ}
\i(Y_\si)\fl_{n/2-1}I^g.
\end{equation}
This is essentially the formula for $Q$ given by Gover
and Peterson in \cite{GoPet} (see Proposition 2.7). In fact, in that
paper, $Q$
is given by $\i(Y_\si) F I^g$, where $F:\ct[-1]\to \ct[1]$ is an
invariant operator derived from the ambient powers of the Laplacian
in a similar, but not identical way to $\fl_{n/2-1}$.  It is possible
that $F$ and $\fl_{n/2-1}$  differ as operators on $\ct[-1]$, but
up to scale they agree on $I^g$.

Next we observe that re-expressing our formula for the
Q-curvature recovers a formula given recently by Fefferman and
Hirachi. An ambient expression naturally corresponding to
$\i(Y_\si)\fl_{n/2-1}I^g=\i(Y_\si)\fl_{n/2-1}\i(\fD)\e(X)Y_\si$ is
$\i(\Y_{\!\!\si})\afl^{n/2-1}\i(\afD)\e(\X)\Y_{\!\!\si}$, 
where, with $\tilde{\si}\in
\cce(1)$ a homogeneous function on $\aM$ corresponding to $\si\in
\ce[1]$, we define $\Y_{\!\!\si}: =\tilde{\si}^{-1}\ad \tilde{\si}$.  
Now using Lemma
\ref{old-domino} we can re-express this as
$\i(\Y_{\!\!\si})\i(\X)\e(\afD)\afl^{n/2-1}\Y_{\!\!\si}$. Next observe that
$\Y_{\!\!\si}=\ad \log \tilde{\si}$ and that $[\afl^{n/2-1},\ad]=0$. So we
obtain 
$$
\i(\Y_{\!\!\si})\i(\X)\e(\afD)\ad \afl^{n/2-1} \log \tilde{\si}.
$$ 
But, from the formula
for $\e(\afD)$ in Proposition \ref{efdd}, we have 
$\e(\afD)\ad=\e(\X)\ad \afl$. Also note that  $\{\i(\X),\e(\X)\}$ 
vanishes modulo $O(Q)$, while  
$$\{\i(\Y_{\!\!\si}),\e(\X)\}=h(\X,\Y)=1.
$$ 
So, modulo $O(Q)$ terms, we get to 
$-\i(\X)\ad \afl^{n/2} \log \tilde{\si}
=-\cL_{\sX}\afl^{n/2} \log \tilde{\si} $. Since 
$\afl^{n/2} \log \tilde{\si} $ is homogeneous of degree $-n$, we obtain
finally
\begin{equation}\label{FHQ}
n \al^{n/2} \log \tilde{\si},
\end{equation}
modulo $O(Q)$ terms. 
Up to a non-zero constant 
multiple, this is the ambient expression for the Q-curvature
given in \cite{FeffHir}.  In summary,  we see that under the 
identification of
tractor sections with appropriate homogeneous ambient quantities, and
with the use of standard identities from exterior calculus, \nn{FHQ} and
\nn{GPQ} are really identical
formulae for the Q-curvature, and both are
generalised by $Q^\si_k$ as in Theorem \ref{Qthm}.

There is a corresponding ambient expression for
the cases $k\geq 1$. Recall that $\bQ^{\ell,\sigma}_k$ is given by
$-2(\ell+1) q^{k}\i(Y_\si)\fb_{\ell-1}\i(\fD)\e(X)\e(Y_\si)$, where
$\fb_{\ell-1}=\frac{1}{2}(\fl_{\ell-1}+(\fl_{\ell-1})^*)$. 
Following the case above, if we
replace $ \fb_{\ell-1}$ by $\fl_{\ell-1}$, then we obtain an alternative
operator 
$$
\tilde{\bQ}^{\ell,\sigma}_k:=-2(\ell+1)
q^{k}\i(Y_\si)\fl_{\ell-1}\i(\fD)\e(X)\e(Y_\si) 
$$ 
that agrees with
$\bQ^{\ell,\sigma}_k$ at leading order and has a conformal
transformation law very similar to that of
$\bQ^{\ell,\sigma}_k$. Up to a multiple the ambient
expression naturally corresponding to $\tilde{\bQ}^{\ell,\sigma}_k$ is
$\qq^k\i(\Y_{\!\!\si})\afl^{\ell-1}\i(\afD)\e(\X)\e(\Y_{\!\!\si})$.
Here $\qq^k$ is an algebraic map on $\acg_k$ corresponding to 
$q^k$ on $\cg_k$ (so that along $\cq$ these maps are equivalent).
Now viewing this as
an operator on ambient forms $U$ of degree $k$ and such that
$\e(\X)\ad U=O(Q)$ along $\cq$ we obtain, by a very similar argument to 
that above, a re-expression of this as 
$$
-\qq^k \i(\X)\ad \afl^\ell \e( \log \tilde{\si}).
$$

Finally we should say that, in an obvious way, the constructions of
Section \ref{semi} above may be carried out on the ambient manifold,
and in that setting, the observation $\Y_{\!\!\si}=\ad \log \tilde{\si}$ 
may be used
to express the $\bS^\si$ operators in terms of ambient operators acting on
$\log \tilde{\si} $. 
It follows that the constructions there may also
be viewed as a generalisation of Theorem 2.2 of \cite{FeffHir}.

\end{document}